\documentclass[12pt]{amsart}
\usepackage{amssymb, amsmath, amsfonts}
\usepackage[raiselinks=false,colorlinks=true,
citecolor=blue,urlcolor=blue,linkcolor=blue,
bookmarksopen=true,pdftex]{hyperref}
\newtheorem{theorem}{Theorem}[section]
\newtheorem{proposition}[theorem]{Proposition}
\newtheorem{lemma}[theorem]{Lemma}
\newtheorem{remark}[theorem]{Remark}

\newtheorem{corollary}[theorem]{Corollary}

\usepackage{graphicx}%
\usepackage[left=1in,right=1in,top=1in,bottom=1in]{geometry}
\usepackage{enumerate}
\usepackage{tikz}
\usetikzlibrary{positioning}
\usetikzlibrary{patterns}
\usepackage{float}

\renewcommand{\hom}{\textrm{Hom}}

\newcommand{\n}{\textrm{n}}
\newcommand{\SL}{\textrm{SL}}
\newcommand{\Sp}{\textrm{Sp}}
\newcommand{\SO}{\textrm{SO}}

\newcommand{\SU}{\textrm{SU}}
\newcommand{\Fix}{\text{Fix}}

\newcommand{\Res}{\textrm{Res}}

\def\r{\mathbb{R}}

\def\f{\mathfrak}

\begin{document}
\baselineskip=15.5pt
\title[Geometric cycles and automorphic representations]{Geometric cycles in compact locally Hermitian symmetric spaces and 
automorphic representations}  
\author[A. Mondal]{Arghya Mondal}
\address{Indian Statistical Institute, 8th Mile, Mysore Road, Bengaluru 560059, India.}
\author[P. Sankaran]{Parameswaran Sankaran}
\address{The Institute of Mathematical Sciences, (HBNI), Chennai 600113, India.}
\email{arghya\_vs@isibang.ac.in}
\email{sankaran@imsc.res.in}
\subjclass[2010]{22E40, 22E46\\
Keywords and phrases: Geometric cycles, automorphic representation.}

\date{}

\begin{abstract}   
Let $G$ be a linear connected non-compact real simple Lie group 
and let  
$K\subset G$ be a maximal compact subgroup of $G$.  Suppose that the centre of $K$ isomorphic to $\mathbb{S}^1$ 
so that $G/K$ is a global Hermitian symmetric space.
Let $\theta$ be the Cartan involution of $G$ that fixes $K$. 
Let $\Lambda$ be a uniform lattice in $G$ such that $\theta(\Lambda)=\Lambda.$
Suppose that $G$ is one of the groups $\SU(p,q), p<q-1, q\ge 5, \SO_0(2,q)$, $\Sp(n,\mathbb{R}), n\ne 4, \SO^*(2n), n\ge 9.$ 
Then there exists a unique irreducible unitary representation 
$\mathcal{A}_\mathfrak{q}$ associated to a proper $\theta$-stable parabolic subalgebra $\mathfrak{q}$ with $R_+(\mathfrak{q})=R_-(\mathfrak{q})$ such that if $H^{s,s}(\mathfrak{g},K;A_{\mathfrak{q}',K})\ne 0$ for some $0<s\le R_+(\mathfrak{q})$, then $\mathcal{A}_{\mathfrak{q}'}$ is unitarily equivalent to either the trivial representation or to $ \mathcal{A}_{\mathfrak{q}}$.    
As a consequence, under suitable hypotheses on $\Lambda,$ we show that the multiplicity of $\mathcal{A}_\mathfrak{q}$ 
occurring in $L^2(\Gamma\backslash G)$ is positive for {\it any} torsionless lattice $\Gamma\subset G$ commensurable with $\Lambda$.   

\end{abstract}
\maketitle
\section{Introduction}
Let $G$ be a linear connected real semisimple Lie group which is non-compact and let $K$ 
denote a maximal compact subgroup 
of $G$.    
Let $\theta:G\to G$ be the Cartan involution of $G$ which fixes $K$.   
We denote by $\mathfrak{g}_0, \mathfrak{k}_0$ the Lie algebras of $G$ and $K$ respectively and by 
the same symbol $\theta$ the involution of the Lie algebra $\mathfrak{g}_0$, which is the differential 
of $\theta:G\to G.$  The complexification of $\mathfrak{g}_0, \mathfrak{k}_0$ will be denoted 
by $\mathfrak{g}, \mathfrak{k}$, etc.  One has the Cartan decomposition $\mathfrak{g}_0
=\mathfrak{k}_0\oplus \mathfrak{p}_0$ where $\mathfrak{p}_0$ is the $(-1)$-eigenspace of $\theta$.
We denote by $\mathfrak{p}$ the complex vector space $\mathfrak{p}_0\otimes_\mathbb{R}\mathbb{C}$.  Thus 
$\mathfrak{g}=\mathfrak{k}\oplus \mathfrak{p}$.  

Since $G$ is semisimple, the Killing form of $\mathfrak{g}_0$ restricted to $\mathfrak{p}_0$ is positive definite 
and induces a $G$-invariant Riemannian metric on $X:=G/K$ with respect to which $X$ is a symmetric space.   
We will assume that $G$ is not a complex Lie group so that $X$ is a product of irreducible symmetric spaces 
of type III (see \cite{helgason}).   This condition automatically holds when $X$ is a Hermitian symmetric space.  

Let $\Gamma\subset G$ be a torsionless uniform lattice in $G$ so that $\Gamma \backslash 
G/K=\Gamma\backslash X=:X_\Gamma$ is a compact locally symmetric space which is an 
Eilenberg-MacLane space $K(\Gamma,1)$.   
The cohomology algebra $H^*(X_\Gamma;\mathbb{C})=H^*(\Gamma;\mathbb{C})$ is an important 
object of study and is of interest not only in topology but also in number theory and representation theory. 
Our aim here is to construct so-called special cycles whose Poincar\'e duals are non-zero cohomology classes 
in $H^*(\Gamma;\mathbb{C})$ when $X$ is an irreducible {\it Hermitian} (non-compact) symmetric space.  
Our results have implications to occurrence with non-zero multiplicity of the irreducible unitary representations $(\mathcal{A}_\mathfrak{q}, A_\mathfrak{q})$ 
of $G$ associated to certain $\theta$-stable parabolic algebras $\mathfrak{q}\subset \mathfrak{g}$ of $\mathfrak{g}_0$ 
in the Hilbert space $L^2(\Gamma\backslash G)$.  

Special cycles, which are closed oriented totally geodesic submanifolds of $X_\Gamma$, whose Poincar\'e duals
are non-zero cohomology classes, were first constructed by Millson and Raghunathan \cite{mr} 
when $G=\SU(p,q), \SO_0(p,q), \Sp(p,q)$.   Schwermer and Waldner \cite{sw} dealt with the case $G=\SU^*(2n)$ 
and Waldner, when $G$ is the non-compact real form of the exceptional complex Lie group $G_2$.  More recently, 
the cases $G=\SL(n,\mathbb{R}), \SL(n,\mathbb{C})$ were considered by Susanne Schimpf \cite{schimpf} 
and the case $G=\SO^*(2n)$ by Arghya Mondal \cite{mondal}. (See also \cite{mondal-sankaran}.) 

Millson and Raghunathan's construction yields a {\it pair} of special cycles $C, C'\subset X_\Gamma$  whose dimensions 
add up to the dimension of 
$X_\Gamma$.   In fact $C$ and $C'$ are sub locally symmetric space $X(\sigma)_{\Gamma(\sigma)}$ and $X(\sigma\circ \theta)_{\Gamma(\sigma\circ\theta)}$, where $\sigma$ arises from an algebraically defined involutive automorphism that commutes 
with the Cartan involution, $X(\sigma)=G(\sigma)/K(\sigma)$, $G(\sigma)\subset G$ is the subgroup of fixed points of 
$\sigma$,  $K(\sigma)=G(\sigma)\cap K$, and 
$\Gamma(\sigma)=G(\sigma)\cap \Gamma$.  The involutions $\sigma, \theta$ are required to stabilize $\Gamma$ so that 
$\Gamma(\sigma), \Gamma(\sigma\circ \theta)$ are lattices in $G(\sigma), G(\sigma\circ \theta)$.  
Under certain additional hypotheses on the special cycles which ensure that their intersection is transverse,   
and, if necessary, replacing $\Gamma$, which is assumed to be arithmetic, by a suitable finite index subgroup, they showed that 
the cup-product of the Poincar\'e duals $[C],[C']$ of $C$ and $C'$ is a non-zero class in the top cohomology 
of $X_\Gamma$.   They deduced that the Poincar\'e dual of such special cycles cannot arise from a $G$-invariant form.  (See \cite[Theorem 2.1]{mr}.)   Equivalently, their dual cohomology classes are  
not in the image of the Matsushima homomorphism $H^*(X_u;\mathbb{C})\to H^*(X_\Gamma;\mathbb{C})$.  Here $X_u$ denotes the compact dual of the non-compact symmetric space $X$.   Rohlfs and Schwermer \cite{rs} obtain an excess intersection formula leading to a criterion for the non-vanishing of the cup-product of special cycles in a more general setting.  

In order to state our main results, we first recall some well-known results concerning the Hilbert space of square-integrable 
functions on $\Gamma\backslash G$, where $G$ is any non-compact semisimple Lie group with finite centre and 
$\Gamma$ a lattice in $G$.  To a Haar measure on $G$ is associated $G$-invariant measure on 
$\Gamma\backslash G$ with finite volume.   The Hilbert space $L^2(\Gamma\backslash G)$ affords a unitary representation 
of $G$ via the translation action of $G$ on $\Gamma\backslash G$. 
When $\Gamma$ is a uniform lattice, 
Gelfand and Pyatetskii-Shapiro \cite{ggp},\cite{gp} proved that $L^2(\Gamma\backslash G)$ decomposes into a Hilbert direct sum of 
irreducible unitary representations $(\pi, H_\pi)$ of $G$ each occurring with {\it finite} multiplicity $m(\pi, \Gamma)$.  Those 
unitary representations $\pi$ such that $m(\pi,\Gamma)$ are positive are referred to as {\it automorphic} representations. 
Let $K\subset G$ be a maximal compact subgroup of $G$.   If $V$ is any $G$-representation 
on a Hilbert space, we denote by $V_K$ the space of all smooth $K$-finite vectors of $V$.  Recall that $V_K$ 
is a $(\mathfrak{g},K)$-module, known as the Harish-Chandra module of $V$. 

The cohomology of $X_\Gamma=\Gamma\backslash G/K$ is described in terms of the relative Lie algebra cohomology 
by the Matsushima isomorphism \cite{matsushima}:
\[H^*(X_\Gamma;\mathbb{C})\cong H^*(\mathfrak{g},K;L^2(\Gamma\backslash G)_K)=\oplus_{\pi}m(\pi,\Gamma) H^*(\mathfrak{g},K;H_{\pi, K}).\]   
A theorem of D. Wigner says that if $(\pi, H_\pi)$ is 
an irreducible unitary representation of $G$, then $H^*(\mathfrak{g},K;H_{\pi,K})$ is non-zero only when its infinitesimal 
character  $\chi_\pi$ is trivial (that is, $\chi_\pi=\chi_0,$ the infinitesimal character of the (trivial) representation $\mathbb{C}$). 
The irreducible unitary representations which have trivial infinitesimal characters have been classified in 
terms of $\theta$-stable parabolic subalgebras $\mathfrak{q}\subset \mathfrak{g}$ of $\mathfrak{g}_0$.   If $\mathfrak{q}$ is a 
$\theta$-stable parabolic subalgebra of $\mathfrak{g}_0$, we will denote the corresponding irreducible unitary 
representation of $G$ by $(\mathcal{A}_\mathfrak{q},A_\mathfrak{q})$ and set $m(\mathfrak{q},\Gamma)
:=m(\mathcal{A}_\mathfrak{q},\Gamma)$.   One has the equivalence relation $\sim$ on the set of all $\theta$-stable parabolic subalgebras of $\mathfrak{g}_0$ where $\mathfrak{q} \sim \mathfrak{q}'$ if $\mathcal{A}_{\mathfrak{q}}$ is unitarily equivalent to $\mathcal{A}_{\mathfrak{q}'}$.   The set $\mathfrak{Q}$ of equivalence classes of $\theta$-stable parabolic subalgebras of $\mathfrak{g}_0$ is finite.  

Suppose that $X=G/K$ is a Hermitian symmetric space.  Then $X_\Gamma$ is a compact K\"ahler manifold and we have the Hodge decomposition 
$H^r(X_\Gamma;\mathbb{C})\cong \oplus_{p+q=r} H^{p,q}(X_\Gamma;\mathbb{C})$.   Also one has a Hodge decomposition $H^r(\mathfrak{g},K;
H_{\pi, K})=\oplus_{p+q=r}H^{p,q}(\mathfrak{g},K;H_{\pi,K})$ for any unitary representation $\pi$ of $G$.  
See \cite[Ch. II,\S4]{borel-wallach}.
When $\pi=\mathcal{A}_{\mathfrak{q}}$, there exists a pair of integers $R_+(\mathfrak{q}),R_-(\mathfrak{q})$ such that 
$H^{p,q}(\mathfrak{g},K;A_{\mathfrak{q},K})=0$ unless $p\ge R_+(\mathfrak{q}), q\ge R_-(\mathfrak{q})$ and $p-q
=R_+(\mathfrak{q})-R_-(\mathfrak{q})$.  Moreover $H^{r+R_+(\mathfrak{q}), r+R_-(\mathfrak{q})}(\mathfrak{g},K;A_{\mathfrak{q}
,K})\cong H^{r,r}(Y_\mathfrak{q};\mathbb{C})$ for a certain compact globally Hermitian symmetric space $Y_\mathfrak{q}$.   
We refer to $(R_+(\mathfrak{q}), R_-(\mathfrak{q}))$ as the {\it Hodge type} of $\mathfrak{q}$ (or of $\mathcal{A}_\mathfrak{q}$) and we set 
$R(\mathfrak{q})=R_+(\mathfrak{q})+R_-(\mathfrak{q})$.   Note that the Hodge type of $\mathfrak{q}$ depends only on its class in $\mathfrak{Q}$.  
The Matsushima 
isomorphism preserves the Hodge decomposition, that is, its inverse maps $H^{p,q}(\mathfrak{g},K;A_{\mathfrak{q},K})$ into 
$H^{p,q}(X_\Gamma;\mathbb{C})$ for all $p,q$.  

Suppose that $X=G/K$ is an irreducible Hermitian symmetric space. 
Denote by $r(\mathfrak{g}_0)$ the smallest positive integer $r$ such that there exists a $\theta$-stable parabolic subalgebra $\mathfrak{q}$ of 
$\mathfrak{g}_0$ with $R_+(\mathfrak{q})=R_-(\mathfrak{q})=r$.  See Table \ref{valuesofr0} for the values of $r(\mathfrak{g}_0)$.    
 
Let $F$ be a totally real number field $F\ne \mathbb{Q}$ and let 
$u\in F_{>0}$ be an element all of whose conjugates, except $u$ itself, are negative. Then one obtains, via 
Weil's restriction of scalars, 
a {\it uniform} lattice $\Gamma(F,u)\subset G$ arising from an $F$-structure on $\mathfrak{g}_0$ 
using a suitable Chevalley basis of $\mathfrak{g}$.  (The Chevalley basis is assumed to be adapted to $\mathfrak{t}$ where $T\subset K$ is a compact Cartan subgroup of $G$ and to the compact form 
$\mathfrak{k}\oplus i\mathfrak{p}=\mathfrak{g}_\textrm{u}$ in the sense of \cite{borel63}.)  
This construction is due to Borel \cite{borel63}.   Let 
$\mathcal{L}(G)$ be the family consisting of all torsionless lattices $\Lambda\subset G$ which are commensurable 
to $\Gamma(F,u)$ for some pair $(F,u)$.

\begin{theorem}  \label{main1} 
We keep the above notations.  
Let $G$ be one of the groups $\SU(p,q), 1\le p<q-1, q\ge 5$, $\SO_0(2,p), p\ge 3, \Sp(n,\mathbb{R}), n\ne 4$, and $\SO^*(2n), n\ge 9$.  
Then there exists a unique irreducible unitary representation $\mathcal{A}_\mathfrak{q}$ of $G$ where $R_\pm(\mathfrak{q})=r(\mathfrak{g}_0).$  Moreover, 
$\mathcal{A}_\mathfrak{q}$ occurs with non-zero multiplicity in $L^2(\Gamma \backslash G)$ 
for every $\Gamma\in \mathcal{L}(G)$. 
\end{theorem}

The above theorem leaves out the infinite families $G=\SU(p,p), \SU(p,p+1)$, the exceptional groups with Lie algebras 
$\mathfrak{e}_{6,(-14)}, \mathfrak{e}_{7,(-25)} $ and a few classical groups for small complex rank.   
We do consider these cases also and obtain, but a weaker result.  See \S \ref{proofs}.

The above theorem will be obtained as an application of the following theorem.  
We associate to each irreducible Hermitian symmetric space $X=G/K$ a number, denoted $c(X)$, as follows: 
$c(\SU(p,q)/K)=p$ where $p\le q$, $c(\SO_0(2,p)/K)=1, c(\Sp(n,\mathbb{R})/K)=n-1, c(\SO^*(2n)/K)=n-1, 
c(E_{6,(-14)}/K)=6, c(E_{7,(-25)}/K)=11$.   Here $E_{6, (-14)}, E_{7,(-25)}$ are exceptional Lie groups 
with Lie algebras $\mathfrak{ e}_{6,(-14)}, \mathfrak{e}_{7,(-25)}$ respectively. 
The significance of $c(X)$ is that, as we shall see, there 
exists a complex analytic special cycle in $X_\Gamma$ of complex dimension $c(X)$ for $\Gamma \in \mathcal{L}(G)$.

\begin{theorem} \label{main2}
We keep the above notations.  
Let $\Gamma\in \mathcal{L}(G)$. 
For any integer $r$ such that $c(X)\le r\le \dim_\mathbb{C} X-c(X)$, there exist a non-zero cohomology 
class in $H^{r,r}(X_\Gamma;\mathbb{C})$ which is not in the image of the Matsushima homomorphism 
$H^{2r}(X_\textrm{u};\mathbb{C})\to H^{2r}(X_\Gamma;\mathbb{C})$.   
\end{theorem}

Theorem \ref{main1} seems to be a new addition to the vast literature on the non-vanishing and the asymptotic behaviour 
of the multiplicity of automorphic representations in $L^2(\Gamma\backslash G)$ in various settings, including, 
the work of Anderson \cite{anderson}, Clozel \cite{clozel}, DeGeorge and Wallach \cite{degeorge-wallach}, and 
Li \cite{li}.  See also \cite[\S6]{parthasarathy80} and \cite[Ch. VIII]{borel-wallach}.  
It should be pointed out that the work of Li \cite{li} establishes non-vanishing results for $m(\pi,\Gamma)$ 
in a general setting using entirely different (and rather deep) techniques, but it does not cover the case of the representation $\mathcal{A}_\mathfrak{q}$ as in Theorem \ref{main1} when $G=\Sp(n, \mathbb{R}), \SO^*(2n)$ or $\SO_0(2,p)$ with $p$ odd. 
This is because, in these cases the group 
$L\subset G$ corresponding to the Lie algebra $\mathfrak{l}_0= 
\mathfrak{q}\cap \bar{\mathfrak{q}}\cap \mathfrak{g}_0$ has more than one non-compact simple factor and so 
does not satisfy the hypotheses of \cite[Prop. 6.1]{li}.
(See \S \ref{pptypeparabolics}.)   
When $G=\SU(p,q)$ or $\SO_0(2,p), p$ even, Theorem \ref{main1} does follow from the work of 
Li, at least when $\Gamma $ is sufficiently `deep'.  

This paper was inspired by the work of Schwermer and Waldner \cite{sw}.  

We will prove both the theorems simultaneously. 
Our proofs are quite elementary and involves Lie theory in identifying elements of $\mathfrak{Q}$ having Hodge 
type $(r,r)$, especially when $r\le c(X),$ and exploits well-known cohomological consequences resulting from 
the existence of complex analytic cycles in a compact K\"ahler manifold.   The construction of the 
lattices in $\mathcal{L}(G)$ is recalled in \S \ref{commutinginvolutions}.   
The group of commuting involutions obtained in Proposition 
\ref{involutions}, which is applicable in greater generality, is used in the construction of {\it analytic} special cycles.  In 
\S\ref{pptypeparabolics} we determine all $\theta$-stable parabolic subalgebras of Hodge type $(r,r)$ for $r\le c(X)$.   
The main theorems are proved in the last section.
\newpage
\begin{center}
{\bf List of notations}
\end{center}
\begin{tabular}{ll}
$G\supset K\supset T$ & linear connected (semi) simple non-compact Lie group, a maximal compact \\
& subgroup of $G$, a maximal torus.\\
$\mathfrak{g}_0, \mathfrak{k}_0, \mathfrak{t}_0$& $\text{Lie}(G), \text{Lie}(K), \text{Lie}(T)$\\
$\mathfrak{g},\mathfrak{k},\mathfrak{k}$ & complexifications of $\mathfrak{g}_0,\mathfrak{k}_0, \mathfrak{t}_0$\\
$X, X_\textrm{u}$ & the symmetric space $G/K$, compact dual of $X$\\
$\Gamma$, $\mathcal{L}(G)$ & a (uniform) lattice in $G$, a certain collection of torsionless uniform lattices in $G$.\\
$G(\sigma), X(\sigma)$ & fixed points of an automorphism $\sigma$ of $G$.\\
& respectively an isometry $\sigma$ of $X$.\\
$X_\Gamma$ & $\Gamma\backslash X$\\
$C(\sigma, \Gamma)$ & $\Gamma(\sigma)\backslash X(\sigma)$ \\
$\theta$  &  Cartan involution of $G$ that fixes $K$, the induced Lie algebra \\
&automorphism of $\mathfrak{g}_0$ or its complexification.\\
$\mathfrak{p}_0, \mathfrak{p}$ & the $(-1)$-eigenspace of $\theta:\mathfrak{g}_0\to \mathfrak{g}_0$, its complexification.\\
$\mathfrak{p}_+, \mathfrak{p}_-$ & the $K$-submodules of $\mathfrak{p}$, holomorphic and anti-holomorphic tangent space at the origin of $X$.\\
$\mathfrak{g}_\textrm{u}$ & the compact form $\mathfrak{k}_0\oplus i\mathfrak{p}_0\subset \mathfrak{g}$ of 
$\mathfrak{g}$.  \\
$\Phi_\mathfrak{g}, \Delta_\mathfrak{g}$ & the roots of $\mathfrak{g}$, the set of simple roots of $\mathfrak{g}.$\\ 
$\Phi^+; \Phi^+_\n, \Phi^-_\n$ & the set of positive roots; the set of positive (resp. negative) non-compact roots.\\
$\alpha_0; \psi, \psi_j,$ & the highest root; simple roots.\\
$\varpi_\psi, \varpi_j$ & fundamental weights of $\mathfrak{g}$.\\
$\lambda, \mu$ & elements of $i\mathfrak{t}^*.$\\
$(.,.)$ & the innerproduct on $i\mathfrak{t}$ or on $i\mathfrak{t}_0^*$ induced by the Killing form on $\mathfrak{g}_0$.\\
$h_\lambda;  H_\lambda\in i\mathfrak{t}_0$ & the element of $i\mathfrak{t}_0$ such that $\lambda(h)=(h_\lambda, h); 2h_\lambda/||h_\lambda||^2$\\
$\# A$ & cardinality of a set $A$.\\
$|A|$ & sum of elements of $A$ when $A$ is a finite set of vectors in a vector space.\\
$\mathfrak{q}, \mathfrak{q}_x,\mathfrak{q}_\lambda$ &  $\theta$-stable parabolic subalgebras of $\mathfrak{g}_0$\\
$\mathfrak{l}, \mathfrak{u}$ & a Levi subalgebra and the unipotent radical of $\mathfrak{q}$ respectively.\\
$\mathfrak{m}_F, \widetilde{\mathfrak{m}}_\mathbb{R}$ & an $F$-form of $\mathfrak{g}_0$, $\mathfrak{m}_F\otimes_\mathbb{Q}\mathbb{R}$\\  
$\mathbf{M}, \mathcal{M}, \mathcal{M}_\mathbb{R}$ & an $F$-algebraic group, $\Res_{F|\mathbb{Q}}\mathbf{M}$, $\mathbb{R}$-points of $\mathcal{M}$.\\
$(\pi, H_\pi)$ & an irreducible unitary representation of $G$ on a Hilbert space $H_\pi$.\\
$H_{\pi,K}$ & the $(\mathfrak{g},K)$-module of smooth $K$-finite vectors in $ H_\pi$.\\
$(\mathcal{A}_\mathfrak{q}, A_\mathfrak{q})$ & a certain irreducible unitary representation associated to $\mathfrak{q}$.\\
\end{tabular}

\section{Geometric cycles} 
Let $ G$ be a connected real semi simple linear Lie group without compact factors and let $K$ be a maximal 
compact subgroup of $G$.  We have the Cartan decomposition $\mathfrak{g}_0=
\mathfrak{k}_0\oplus \mathfrak{p}_0$; thus $\mathfrak{p}_0$ is the $(-1)$-eigenspace of $\theta$.  Let $X=G/K$.  We shall refer 
to the trivial coset as the origin of $X$. 
 
For any lattice $\Gamma$ in $G$, let $X_\Gamma:=\Gamma\backslash X=\Gamma\backslash G/K$.   
We shall assume that $\Gamma$ is torsionless, 
irreducible, and uniform.  Thus $X_\Gamma$ is a smooth compact manifold.  The $G$-invariant metric on $X$ descends to yield a Riemannian metric on $X_\Gamma$.  Our aim in this section is to exhibit, for suitable lattices $\Gamma\subset G$, pairs of {\it geometric cycles} $C_1, C_2$ in $X_\Gamma$ of complementary dimensions (i.e., $\dim C_1+\dim C_2=\dim X$) 
such that the cup-product of their Poincar\'e duals $[C_1],[C_2]$ is a non-zero element of $H^*(X_\Gamma;\mathbb{C})$ when $X$ is an 
irreducible Hermitian symmetric space (of non-compact type).   The significance of such pairs is that the cohomology classes $[C_j], j=1,2,$ are then {\it not} representable by $G$-invariant forms on $X$; this is a result due to Millson and Raghunathan \cite[Theorem 2.1]{mr}.  
Equivalently, the Poincar\'e duals $[C_1], [C_2]$ are not in the image of the Matsushima homomorphism $H^*(X_\textrm{u};\mathbb{C})\to H^*(X_\Gamma;\mathbb{C})$.   Millson and Raghunathan constructed pairs of geometric cycles which intersect transversally 
when $G=\SU(p,q), \SO_0(p,q), \Sp(p,q)$ for certain lattices $\Gamma$.  Working in a more general setup where the 
geometric cycles are allowed to intersect along positive dimensional submanifolds,  Rohlfs and Schwermer \cite{rs} obtained 
a formula for the cup-product $[C_1].[C_2]$ when $C_1,C_2$ satisfy a certain orientation condition called Or (to be explained below).   
This has been applied, for suitable uniform arithmetic lattices, when $G=\SU^*(2n)$ by Schwermer and 
Waldner, for the (non-compact) exceptional group of type $\textrm{G}_2$ by Waldner \cite{waldner} and for the groups 
$\SL(n,\mathbb{R}), \SL(n,\mathbb{C})$ by Schimpf \cite{schimpf}.   More recently, Mondal \cite{mondal} has considered the case $G=\SO^*(2n).$    (See also \cite{mondal-sankaran}.) In this paper we use the term \textit{special cycles} interchangeably with the term geometric cycles, although the special cycles considered by Rohlfs and Schwermer are more general.

Suppose that $G$ is simple so that $X$ is irreducible.    Let $\sigma_1$ be an involutive automorphism of $G$ that 
stabilizes $K$ and that $\sigma_1$ commutes with the Cartan involution $\theta$.   
Set $\sigma_2:=\sigma_1\circ \theta$.   
It is known (\cite{borel63}, \cite[\S2]{rohlfs-speh}) that there are arithmetic lattices $\Gamma\subset G$ such that $\theta(\Gamma)=\Gamma$.  
We assume that $\sigma_1(\Gamma)$ and $\Gamma$ are commensurable so that, by passing to a finite index subgroup 
of $\Gamma$ if necessary, we have $\sigma_j(\Gamma)=\Gamma, j=1,2.$  
Let $G_j=\Fix(\sigma_j)$ and let $K_j=K\cap G_j$.  The group $G_j$ is in general a reductive Lie subgroup, not necessarily 
semi simple.   In any case $X_j:=G_j/K_j$ is a Riemannian symmetric space that naturally embeds in $X=G/K$ as a totally 
geodesic submanifold.     
Denote by $C_j$ the image of $X_j$ under the projection $X\to X_\Gamma$.   
Setting $\Gamma_j:=G_j\cap \Gamma$, the $C_j=\Gamma_j\backslash G_j/K_j, j=1,2,$ are closed submanifolds of $X_\Gamma$ of complementary dimensions.  
Following \cite{rs} one says that $C_j$ satisfies condition Or if the action of $G_j$ on (the left of) $X_j$ is orientation preserving.   
This requirement is trivially valid when $G_j$ is connected.     It is also valid when $X$ is Hermitian symmetric and 
$\sigma_j:X\to X$ commutes with translation by elements of the centre of 
$K$; see \cite{rs}.   Moreover, in this case the $X_j$ are also Hermitian symmetric and the inclusions $X_j\hookrightarrow X$ and $C_j\hookrightarrow X_\Gamma$  
are complex analytic.  As $X_\Gamma$ is K\"ahler---in fact it is a smooth complex projective variety by a theorem of Kodaira \cite[Theorem 6]{kodaira}---so are the $C_j$.

Our aim is to show the existence of a uniform lattice in $G$ stabilized by $\sigma, \theta,$ 
where $\sigma\in Aut(G)$ is such that $\sigma(K)=\theta(K)=K$.   
We will achieve this by choosing an appropriate $F$-algebraic group $\mathbf{M}$, with $\mathcal{M}$ the $\mathbb{Q}$-algebraic group obtained by applying Weil's restriction of scalars functor $\Res_{F|\mathbb{Q}}$  to $\mathbf{M}$, such that (i) $G/Z(G)$ equals the identity component of the Lie group given by the 
$\mathbb{R}$-points $\mathcal{M}_\mathbb{R}$ modulo its maximal compact connected normal subgroup   
(ii) $\sigma$ and $\theta$ are induced by $F$-rational involutions $\sigma_F, \theta_F$ of $\mathbf{M}$.  
The existence of $F$-rational structures and an $F$-rational Cartan involution $\theta_F$ that induces $\theta $ on 
$G$ are well-known \cite{borel63}, \cite{rohlfs-speh}.  We shall proceed as in \cite{rohlfs-speh} to show 
the existence of $\sigma_F$ that commutes with $\theta_F$.   It suffices to do this at the level of Lie algebras (as in \cite{rohlfs-speh}). 

\subsection{Commuting family of involutions}\label{commutinginvolutions}
Throughout this section we suppose that $G$ is a connected semisimple linear Lie group without 
compact factors.   Let $K$ be a maximal compact subgroup of $G$ and denote by $\theta$ the Cartan involution of $G$ that 
fixes $K$.    Let $T\subset K$ be a maximal torus 
in $K$. We assume that $\mathfrak{t}$ is a Cartan subalgebra of $\mathfrak{g}$, although $G/K$ is not required to be Hermitian symmetric. This hypothesis simplifies the exposition of the Chevalley basis of 
$\mathfrak{g}$ needed in the construction 
of $\theta$-stable uniform lattices to be described below, 
although Borel obtained his results in complete generality.

Denote the Killing form on $\mathfrak{g}$ by $(.,.)$.  Its
restriction to $\mathfrak{t}$ is non-degenerate and hence yields an isomorphism $\mathfrak{t}\cong \mathfrak{t}^*$ 
and an induced bilinear form on $\mathfrak{t}^*$ denoted by the same symbol.  It is an innerproduct on 
$i\mathfrak{t}_0^*$. 
For any non-zero $\lambda\in \mathfrak{t}^*$, 
we denote by $h_\lambda\in\mathfrak{t}$ the unique element so that 
$\lambda(H)=(H,h_\lambda)$ and set $H_\lambda:=2h_\lambda/||\lambda||^2$.  
Note that 
$(\lambda,\mu)=(h_\lambda,h_\mu)=\mu(h_\lambda)$ and that if $\lambda\in i\mathfrak{t}_0^*$, then 
$h_\lambda\in i\mathfrak{t}_0$.   
Let $\Phi=\Phi(\mathfrak{g},\mathfrak{t})$ be the set of roots.  Let $\Phi^+\subset \Phi$ be a positive root system and let 
$\Delta_\mathfrak{g}\subset \Phi^+$ be the set of simple roots.

We choose a Chevalley basis $\{H_\gamma\}_{\gamma\in \Delta_{\mathfrak{g}}},  X_\alpha, \alpha\in \Phi_\mathfrak{g}=\Phi,$ of $\mathfrak{g}$ 
adapted to $\mathfrak{t}$ and 
the compact form $\mathfrak{g}_\textrm{u}=\mathfrak{k}_0\oplus i\mathfrak{p}_0$ so that the structure constants are all rationals, that is:
\[\sum_{\gamma\in \Delta_{\mathfrak{g}}} \mathbb{R}H_\gamma=i\mathfrak{t}_0,\]
\[  [H, X_\alpha]=\alpha(H)X_\alpha, ~\forall H\in \mathfrak{t},\forall \alpha\in \Phi,\]
\[ [X_\alpha, X_\beta]=\left\{ \begin{array}{ll} 
N_{\alpha,\beta}X_{\alpha+\beta} & \textrm{if~} \alpha,\beta,\alpha+\beta\in \Phi, \alpha+\beta\ne 0 \\
 H_\alpha & \textrm{if~} \alpha+\beta=0.
 \end{array}\right.         \] 
where $N_{\alpha,\beta}=-N_{-\alpha,-\beta}=\pm(p+1)\in \mathbb{Z}$ where $p\ge 0$ is the greatest integer such that $\alpha-p\beta\in \Phi$.  Set 
$U_\alpha:=(X_\alpha-X_{-\alpha}), V_\alpha:=i(X_\alpha+X_{-\alpha}), \alpha\in \Phi^+$.  
Then $iH_\gamma, \gamma\in \Delta_{\mathfrak{g}}, U_\alpha, V_\alpha, \alpha\in \Phi_{\mathfrak{k}}^+, 
iU_\beta, iV_\beta, \beta\in \Phi_\n^+$  
form a basis for $\mathfrak{g}_0$ with rational structure constants.  For any real number field $F$ let 
$\mathfrak{g}_F$ denote 
the $F$-vector space spanned by these elements.  
Then $\mathfrak{g}_F$ is a Lie algebra over $F$ and is an $F$-form of $\mathfrak{g}_0$, that is, $\mathfrak{g}_F\otimes_F\mathbb{R}=\mathfrak{g}_0$.    Note that the set $\mathcal{B}_{\mathfrak{k}}$ consisting of  $iH_\gamma, 
\gamma\in \Delta_{\mathfrak{g}}, U_\alpha, V_\alpha, \alpha\in \Phi^+_{\mathfrak{k}}$ is an $F$-basis for an $F$-form 
$\mathfrak{k}_F$, 
of $\mathfrak{k}_0$.  We denote by $\mathfrak{p}_F$ the $F$-span the set $\mathcal{B}_\mathfrak{p}=\{iU_\alpha, iV_\alpha\mid \alpha\in \Phi^+_\n\}$.  

Suppose that $F$ is a totally real number field.  Choose an element $u\in F, u>0$, such that $s(u)<0$ for all 
$s\in S=S(F)$, the set of 
all embeddings $s:F\to\mathbb{R}$ other than the inclusion $\iota:F\hookrightarrow \mathbb{R}$. 
Let $E=F(\sqrt{u})$ and let 
\[\mathfrak{m}_F=\mathfrak{k}_F\oplus \sqrt{u} \mathfrak{p}_F.\]
In view of the fact that $[\mathfrak{k}_0,\mathfrak{p}_0]\subset \mathfrak{p}_0, 
[\mathfrak{p}_0,\mathfrak{p}_0]\subset \mathfrak{k}_0$ we see that $\mathfrak{m}_F$ is a 
Lie algebra over $F$.   
It has an $F$-basis $\{b_j\}:=\mathcal{B}_\mathfrak{k}\cup \sqrt{u}\mathcal{B}_\mathfrak{p}$.  
Choose a primitive element  $v\in F$ over $\mathbb{Q}$ and let $d=\deg_\mathbb{Q} F$.  
Then $\{v^l\otimes_\mathbb{Q} b_j\mid 0\le l<d, 1\le j\le \dim G\}$ yields a $\mathbb{Q}$-structure on  
$\mathfrak{m}_F$.

Also, for any $s\in S, \mathfrak{m}^s_F:=
\mathfrak{k}_F\oplus i\sqrt{-s(u)}\mathfrak{p}_F$ is a Lie algebra over $F$ which is an $F$-form 
of $\mathfrak{u}:=\mathfrak{k}_0+i\mathfrak{p}_0$, a 
maximal compact Lie subalgebra of $\mathfrak{g}$.   Again $\mathfrak{m}^s_F$ has a $\mathbb{Q}$-structure 
given by $\{v^l\otimes_\mathbb{Q} s(b_j)\}$ where $s(b_j)=b_j$ if $b_j\in \mathfrak{k}_F$ and $s(b_j)
=i\sqrt{-s(u)} b'_j$ if $b_j=\sqrt{u} b'_j, b'_j\in \mathfrak{p}_F$.    
Since $F$ is totally real, we have an isomorphism of real Lie algebras   
\[\widetilde{\mathfrak{m}}_\mathbb{R}:=\mathfrak{m}_F\otimes_\mathbb{Q} \mathbb{R}
\cong \mathfrak{g}_0\bigoplus (\oplus_{s\in S} \mathfrak{m}^s_\mathbb{R}). \] 
In particular $\widetilde{\mathfrak{m}}_\mathbb{R}$ is a semi simple Lie algebra in which all simple ideals 
{\it not} contained in $\mathfrak{g}_0$ are compact Lie algebras.

Let $\mathbf{M}$ denote the adjoint type 
$F$-algebraic group corresponding to the $F$-Lie algebra $\mathfrak{m}_F$ and let 
$\mathcal{M}=\Res_{F|\mathbb{Q}}(\mathbf{M})$ be the $\mathbb{Q}$-algebraic group obtained from 
$\mathbf{M}$ 
by Weil's restriction of scalars from $F$ to $\mathbb{Q}$.  Let $\textrm{Lie}(\mathcal{M}_\mathbb{Q})$ be the Lie algebra of $\mathcal{M}_\mathbb{Q}$.  Then $\textrm{Lie}(\mathcal{M}_\mathbb{Q})\otimes_\mathbb{Q}\mathbb{R}\cong \widetilde{\mathfrak{m}}_\mathbb{R}$ is the Lie algebra of 
the real Lie group $\mathcal{M}_\mathbb{R}.$   
Denote by $\mathcal{M}_\mathbb{R}^0$ the identity component of $\mathcal{M}_\mathbb{R}$.  
It follows that $\mathcal{M}_\mathbb{R}^0$ has exactly one non-compact 
factor, namely,  $G/Z(G)$.  

The $F$-rational involution $\theta^\iota_F:\mathfrak{m}_F\to \mathfrak{m}_F$ 
defined by $X+\sqrt{u}Y\mapsto X-\sqrt{u}Y$ induces the Cartan involution $\theta=\theta^\iota_F\otimes _F
\mathbb{R}$ on $\mathfrak{g}_0$ and an involution $\theta^s$ on $\mathfrak{m}^s_\mathbb{R}$ given by  conjugation,  $X+i\sqrt{-s(u)}Y\mapsto X-i\sqrt{-s(u)}Y$.  The $\mathbb{Q}$-rational involution $\Res_{F|\mathbb{Q}}(\theta^\iota_F):\mathcal{M}\to \mathcal{M}$ yields an involution 
$\widetilde{\theta}_\mathbb{R}$ on $\mathcal{M}_\mathbb{R}$ which induces the product 
$\theta\times (\prod_{s\in S} \theta^s)$ on $\textrm{Lie}(\mathcal{M}_\mathbb{R})=\widetilde{\mathfrak{m}}_\mathbb{R}$.  
This shows that $\theta:\mathfrak{g}_0\to \mathfrak{g}_0$ arises from a $\mathbb{Q}$-rational involution 
of the $\mathbb{Q}$-group $\mathcal{M}$.   

 In view of the fact that $G/Z(G)$ is the only non-compact factor of $\mathcal{M}^0_\mathbb{R}$, we see that the $\mathbb{Z}$-points $\mathcal{M}_\mathbb{Z}\cap \mathcal{M}^0_\mathbb{R}$ of $\mathcal{M}^0_\mathbb{R}$ projects to a {\it uniform} arithmetic lattice $\bar{\Gamma}$ in  $G/Z(G)$. 
(Here we need $F\ne \mathbb{Q}$
so that there is at least one non-trivial compact factor in $\mathcal{M}_\mathbb{R}^0$.) 

If $p:G\to G/Z(G)$ is the canonical projection, then $\Gamma=\Gamma(F,u):=p^{-1}(\bar{\Gamma})$ is a uniform arithmetic lattice in $G$.  Since $\widetilde{\theta}_\mathbb{R}(\mathcal{M}_\mathbb{Z})=\mathcal{M}_\mathbb{Z},$  it follows that 
$\theta(\Gamma)=\Gamma$.   

{\it We denote by $\mathcal{L}(G)$ the family of all torsionless uniform lattices which are commensurable with 
$\Gamma(F,u)$ as  $F$ varies over all 
totally real number fields and $u$ over $ F_{>0}$ all whose conjugates, other than itself, are negative.  }

\begin{proposition}\label{involutions}
We keep the above notation.  Let $\varpi=\varpi_\psi$ be a fundamental weight corresponding to 
a compact simple root $\psi\in \Delta_\mathfrak{g}$.   Set $t_0:=\pi ||\varpi||^2/||\psi||^2$. 
Then:\\
(i)  The automorphism $\sigma=\sigma_\psi:=e^{\emph{ad} it_0H_\varpi}:\mathfrak{g}_0\to \mathfrak{g}_0$ 
is a $\mathbb{Q}$-rational involution for the $\mathbb{Q}$-structure given by the Chevalley basis. \\
 (ii)  The involution $\sigma_\mathbb{Q}$ commutes with $\theta_\mathbb{Q}$  and defines an $F$-involution $\sigma^\iota_F$ on $\mathfrak{m}_F$ which commutes with $\theta^\iota_F$.  \\
(iii) The involution $\widetilde{\sigma}:=\Res_{F|\mathbb{Q}} (\sigma^\iota_F)$ of $\mathbb{Q}$-group $\mathcal{M}$ commutes with $\widetilde{\theta}=\Res_{F|\mathbb{Q}}(\theta^\iota_F)$.  In particular
$\widetilde{\sigma}_\mathbb{R}$ commutes with $\widetilde{\theta}_\mathbb{R}$.  \\
(iv)  The uniform lattice $\mathcal{M}_\mathbb{Z}\cap \mathcal{M}_\mathbb{R}^0$ 
is preserved by $\widetilde{\sigma}_\mathbb{R}$. \\
(v)  Any two involutions $\sigma_\psi, \psi\in \Delta_\mathfrak{k},$ commute.
\end{proposition}
\begin{proof}   
Evidently $[iH_\varpi,H_\gamma]=0$ for any $\gamma\in \Delta_{\mathfrak{g}}$. Let $\alpha\in \Phi^+$ and write $\alpha=\sum_{\gamma\in\Delta_{\mathfrak{g}}} a_{\alpha,\gamma} \gamma$ as an (integral) 
linear combination of simple roots.  
$[iH_\varpi, X_\alpha]=i\alpha(H_\varpi)X_\alpha= 2i(\alpha,\varpi)||\varpi||^{-2} X_\alpha=i a_{\alpha,\psi} c X_\alpha$ where $c=c(\psi):=(||\psi||/||\varpi||)^{2}$.    
It follows that $[iH_\varpi, U_\alpha]=ca_{\alpha,\psi} V_\alpha, [iH_\varpi, V_\alpha]
=-ca_{\alpha,\psi} U_\alpha$.  Therefore the automorphism $e^{\textrm{ad}i tH_\psi}$ of $\mathfrak{g}_0$ restricts 
to the identity on $\mathfrak{t}_0$, preserves the planes $\mathfrak{k}_{0,\alpha}\subset \mathfrak{k}_0$ spanned by $U_\alpha, V_\alpha; \alpha\in \Phi^+_{\mathfrak{k}},$ and the planes $\mathfrak{p}_{0,\alpha}\subset \mathfrak{p}_0$ 
spanned by $iU_\alpha, iV_\alpha, \alpha\in \Phi^+_\n$.  The matrix $E_\alpha$ of the operator $e^{\textrm{ad}itH_\psi}$ restricted to $\mathfrak{k}_{0,\alpha}$ or to $\mathfrak{p}_{0,\alpha}$
with respect to their respective chosen basis is $e^{tA_\alpha}$ where 
$A_\alpha:=\left(\begin{matrix}  0 & -ca_{\alpha,\psi}  \\
ca_{\alpha,\psi} & 0\end{matrix} \right)$. 
Taking $t_0:=\pi /c$ we see that $E_\alpha=I$ or $-I$ according as $a_{\alpha,\psi}=(1/c)\alpha(H_\varpi)$ is even or odd. 
 As the value of $t_0$ is independent of $\alpha$ we see that $\sigma$ sends each Chevalley basis 
 element either to itself or to its negative.  Hence it preserves the $\mathbb{Q}$-structure on $\mathfrak{g}_0$ 
and is an involution. 
 
 Parts (ii)-(iv) of the proposition follow easily from the observation that the matrix of $\sigma$ with respect to the Chevalley basis is diagonal with eigenvalues $\pm 1$.    Part (v) is trivial since all the $H_\varpi$ 
 belong to the abelian subalgebra $i\mathfrak{t}$.   
\end{proof}

The involution $\sigma_\psi:\mathfrak{g}_0\to \mathfrak{g}_0$ induces an involution of 
the universal cover $\widetilde{G}$ of $G$ which leaves fixed the centre of $Z(\widetilde{G})$.  Hence it 
induces an involution of $G$, which is again denoted by the same symbol $\sigma_\psi$. 
Since $\sigma_\psi(\mathfrak{k}_0)=\mathfrak{k}_0$, we have $\sigma_\psi(K)=K$.  
Hence $\sigma_\psi$ induces an isometry of $X=G/K$ which is also denoted $\sigma_\psi$. 

The group $\Sigma=\Sigma(F,u)$ generated by $\sigma_\psi, \psi\in \Delta_{\mathfrak{k}}$ and $\theta$ is an 
elementary abelian group of order $2^n$ where $n=\dim \mathfrak{t}_0$.  Then the assertions 
(ii), (iii), and (iv) hold for any $\sigma\in \Sigma$.  
We regard $\Sigma$ also 
as a group of isometries of $X$.

From now on, we assume that $X=G/K$ be hermitian symmetric. Since $\sigma_\psi$ commutes  
with the isometries of $X$ defined by elements of $Z(K)$, 
it follows that $\sigma_\psi:X\to X$ is complex analytic.   
Let $\Gamma\subset G$ be a torsionless uniform lattice stabilized by 
a $\sigma\in \Sigma$ and the Cartan involution $\theta$ and set $\tau:=\sigma\circ \theta$.  
Note that if $\Lambda\in \mathcal{L}(G)$, then $\Gamma:=\Lambda\cap \theta(\Lambda)\cap 
\sigma(\Gamma)\cap \theta\circ \sigma(\Lambda)$ is such a lattice.
Denote by $G(\sigma), 
K(\sigma), \Gamma(\sigma)$ the fixed points $\Fix(\sigma), \Fix(\sigma|_K), \Fix(\sigma|_\Gamma)$, etc.  
Then $K(\sigma)$ is a maximal compact subgroup of the reductive group $G(\sigma)$ (which may not be 
connected) and $\Gamma(\sigma)$ is a uniform lattice in $G(\sigma)$. 
We obtain a pair of complementary dimensional special cycles 
$C(\sigma,\Gamma)=X(\sigma)_{\Gamma(\sigma)}, C(\tau,\Gamma)=X(\tau)_{\Gamma(\tau)} \subset X_\Gamma$ which are complex analytic submanifolds of the locally Hermitian symmetric space $X_\Gamma$.  
The group $G(\sigma)$ acts on $X(\sigma)$ preserving the orientation by \cite[Remark 4.8(ii)]{rs}.  The same is true of $G(\tau)$ as well and so condition {\it Or} of \cite[Theorem 4.11]{rs} is met and we have the following 
corollary (cf. \cite[Theorem 2.1]{mr}).  

\begin{corollary} \label{geometriccycles}
Let $\Gamma\subset G$ be a torsionless uniform arithmetic group stabilized by $\sigma\in \Sigma$ 
and by the Cartan involution $\theta$.   Let $\textrm{codim}_{X_\Gamma}C(\sigma,\Gamma)=p$.  
Then 
 $[C(\sigma,\Gamma)]\in H^{p,p}(X_\Gamma;\mathbb{C})$ is {\em not} in the image of the Matsushima homomorphism 
$H^*(X_u;\mathbb{C})\to H^*(X_\Gamma;\mathbb{C})$.  \hfill $\Box$
 \end{corollary}

In view of the fact that $C(\sigma,\Gamma)$ is a complex analytic submanifold of the K\"ahler manifold 
$X_\Gamma$, we have that $[C(\sigma,\Gamma)]\ne 0$.  In fact, $[V]\ne 0$ in $H^*(X_\Gamma;\mathbb{C})$ for any complex analytic subvariety $V\subset C$; see \cite{gh}.  Evidently $[V]$ is of Hodge type $(p,p)$ where 
$p$ is the complex codimension of $V$ in $X_\Gamma$. 

We will be concerned with special cycles $C(\sigma,\Gamma)$ as in the above corollary having {\it minimum} codimension. 
The following corollary, whose proof is immediate from the proof of Proposition \ref{involutions}(i), is a useful tool in 
determining $X(\sigma)$ (the universal cover of $C(\sigma,\Gamma)$) and its compact dual $X(\sigma)_u$, in particular the dimension of the special cycles.  

\begin{corollary} 
With notations as above, let $\psi$ be a compact simple root.  
The Lie algebra of $G(\sigma_\psi)$ 
is $\mathfrak{g}_{0}(\sigma_\psi)=
\mathfrak{t}_0\oplus  \sum \mathfrak{k}_{0,\alpha}\oplus \sum \mathfrak{p}_{0,\beta}$ where the sum is over all 
$\alpha\in \Phi^+_\mathfrak{k}$, (resp. $\beta\in \Phi^+_\n$) such that $(||\varpi_\psi||^2/||\psi||^2)\alpha(H_{\varpi_\psi})$ 
(resp. $(||\varpi_\psi||^2/||\psi||^2)\beta (H_{\varpi_\psi}))$  
is even.   \hfill $\Box$ 
\end{corollary}

\begin{remark}\label{codim}{\em 
(i) We remark that $\theta|_{G(\sigma_\psi)}$ is the Cartan involution of $G(\sigma_\psi)$ that fixes $K(\sigma_\psi):=K\cap G(\sigma_\psi)$.    The expression for $\mathfrak{g}_0(\sigma_\psi)$ in the above corollary is the corresponding 
Cartan decomposition where the last summand equals $\mathfrak{p}_0(\sigma_\psi):=\mathfrak{p}_0\cap \mathfrak{g}_0(\sigma_\psi)$.
 It follows that $T_0(X(\sigma_\psi))=\mathfrak{p}_0(\sigma_\psi)$ and so we obtain the following formula 
 for the (complex) codimension of $X(\sigma)=G(\sigma_\psi)/K(\sigma_\psi)$ in $X$:
 \[c(\sigma_\psi):= \textrm{codim}_X X(\sigma)=\#\{\alpha\in \Phi^+_\n\mid a_{\alpha,\psi}\equiv 1 \mod 2\}
 =\#\Phi^+_\n-\#\Phi^+(\sigma_\psi)_n.\]

(ii) Although $\Delta_{\mathfrak{g}}\cap \Phi(\mathfrak{g}(\sigma_\psi), \mathfrak{t})=\Delta_\mathfrak{g}\setminus\{\psi\}$, 
it is not, in general, the set of simple roots for the positive system $\Phi^+(\sigma_\psi)=
\Phi^+(\mathfrak{g}(\sigma_\psi),\mathfrak{t})=
\Phi^+\cap\Phi(\sigma_\psi)$. }
\end{remark}

\subsection{Outer automorphisms commuting with $\theta$}\label{outerinvolutions}
The involutions of $G$ commuting with $\theta$ arising from the involutions of the Lie algebra $\mathfrak{g}_0$ given by the above proposition are all inner automorphism of $G$.   There are also involutive outer automorphisms which commute 
with $\theta$ when the Dynkin diagram of $\mathfrak{g}$ admits a non-trivial symmetry.    We consider the case
 $G=\SO_0(2,2n-2), K=\SO(2)\times \SO(2n-2)$ in some detail as it will be used later.  We take $T=(\SO(2))^n$ embedded block diagonally in $K$.
 We 
label the simple roots of $\mathfrak{g}$ as in \cite[Planche IV]{bourbaki}.   

Let $H_\psi, \psi\in \Delta_\mathfrak{g}, X_\gamma, \gamma\in \Phi=\Phi(\mathfrak{g},\mathfrak{t})=\{\pm(\epsilon_i\pm\epsilon_j)
\mid 1\le i<j\le n\}$ be a Chevalley 
basis for a $\mathbb{Q}$-form $\mathfrak{g}_\mathbb{Q}$ of $\mathfrak{g}=\mathfrak{so}(2n,\mathbb{C})$ adapted to 
$\mathfrak{t}$ 
and 
$\mathfrak{g}_u$.  Then the Lie algebra automorphism $\tau_\mathbb{Q}:\mathfrak{g}_\mathbb{Q}\to \mathfrak{g}_\mathbb{Q}$ defined by $H_\psi\mapsto H_\psi, \psi\in \Delta_{\mathfrak{g}}, \psi\ne \epsilon_{n-1}\pm \epsilon_n$, 
$H_{\epsilon_{n-1}\pm\epsilon_n}\mapsto H_{\epsilon_{n-1}\mp\epsilon_n}$, $X_\gamma\mapsto X_\gamma, \gamma=\pm 
(\epsilon_i\pm \epsilon_j), 1\le i<j\le n-1$, $X_{\pm(\epsilon_j\pm\epsilon_n)}\mapsto X_{\pm(\epsilon_j\mp\epsilon_n)}$ 
induces an involutive {\it outer} automorphism $\tau:\mathfrak{g}_0\to \mathfrak{g}_0$ that fixes $\mathfrak{so}(2,2n-3)$. 
Denote by $\tau_\mathbb{C}:\mathfrak{g}\to \mathfrak{g}$ the complex linear extension of $\tau$.

 It is evident that $\tau_\mathbb{Q}$ commutes with 
$\theta_\mathbb{Q}$.    If $F$ is any totally real number field and an element $u\in F\cap  \mathbb{R}_{>0}$ such that 
$s(u)<0$ for all embeddings $s:F\to \mathbb{R}$ other than the inclusion $\iota:F\to \mathbb{R}$, we obtain 
from $\tau_\mathbb{Q}$ an involution 
$\tau_F^\iota:\mathfrak{m}_F\to \mathfrak{m}_F$ 
which commutes with $\theta_F^\iota$ (with notations as in \S \ref{commutinginvolutions}).  
We may apply the restriction of 
scalar functor to obtain a $\mathbb{Q}$-automorphism $\tilde{\tau}$ of $\mathcal{M}$ that commutes with $\theta^\iota$.   
This allows us to conclude that the arithmetic group $\mathcal{M}_\mathbb{Z}\cap \mathcal{M}_\mathbb{R}^0$ is stable by $\tilde{\tau}_\mathbb{R}$.  
It follows that the uniform lattice $\Gamma=\Gamma(F,u)\subset \SO_0(2,2n-2)$ (as in \S\ref{commutinginvolutions}) is preserved by $\tau$ and that $\Gamma\cap \SO_0(2,2n-3)$ is a lattice in $\SO_0(2,2n-3)$.  
By passing to a finite index subgroup, we may (and do) 
assume that $\Gamma$ is torsionless.   

Let $\Lambda\in \mathcal{L}(G)$ be a lattice commensurable with $\Gamma=\Gamma(F,u)$ and let $\Lambda(\tau)=\Lambda\cap G(\tau)$.   Then  
$C(\tau,\Lambda)=\Lambda(\tau) \backslash G(\tau) /K(\tau) $ is a special cycle which is complex analytic   
since $\tau$ commutes with the centre $SO(2)\hookrightarrow \SO(2)\times\SO(2n-2)$.  Note that $C(\tau,\Gamma)$ is a divisor in $X_\Gamma$.   
As in Corollary \ref{geometriccycles}, its Poincar\'e dual $[C(\tau,\Lambda)]$ is a non-vanishing 
 cohomology class in $H^{1,1}(X_\Gamma;\mathbb{C})$   
 which is {\it not} in the image of the Matsushima homomorphism 
 $H^*(X_\textrm{u};\mathbb{C})\to H^*(X_\Lambda;\mathbb{C})$.

We summarise the above discussion as a proposition.

\begin{proposition}   Let $G=\SO_0(2,2n-2), K=\SO(2)\times \SO(2n-2)$. 
Let $\Lambda\in \mathcal{L}(G)$.    Then there exists an involutive automorphism 
$\tau:G\to G$ with $\Fix(\tau)=G(\tau)=\SO_0(2,2n-3)$  
such that (i)  $\tau(K)=K$, (ii) $\tau$ commutes with conjugation by any central element of 
$K$; in particular, $\tau\circ \theta=\theta\circ \tau$, (iii) $\tau(\Lambda)\in \mathcal{L}(G)$, and 
(iv) $\Lambda\cap \Fix(\tau)$ is a uniform lattice in $\SO_0(2,2n-3)$, and, (v) The Poincar\'e dual of the special cycle 
$[C(\tau,\Lambda)]\subset \Lambda$ is not in the image of the the Matsushima homomorphism $H^*(X_\textrm{u};\mathbb{C}) \to H^*(X_\Lambda;\mathbb{C})$.  \hfill $\Box$
\end{proposition}

 \begin{remark}{\em 
Millson and Raghunathan \cite{mr} constructed an involutive automorphism $\sigma$ of $\SO_0(2,n)$ so that 
$X(\sigma)\cong \SO_0(2,n-1)/\SO(2)\times \SO(n-1)$ irrespective of the parity of $n$.  In fact they considered  
the more general case when $G=\SO_0(m,n)$ and construct involutions $\sigma_k$ so that $X(\sigma_k)\cong \SO_0(m,k)/(\SO(m)\times \SO(n-k))$, for $1\le k<n$.  They further show that $\sigma$ arises from an $F$-algebraic 
automorphism (for a suitable number field $F\subset \mathbb{R}$) and hence leads to construction of special cycles of complementary dimensions in $\Gamma\backslash X$ for appropriate uniform lattices.  Our approach to the construction 
of $C(\tau,\Gamma)$ is different from theirs.  Although our approach is applicable to the more general case of $\SO_0(p,q)$, we 
will have no need for it for our present purposes. }
\end{remark}


\subsection{Dimensions of Hermitian special cycles}\label{codimension}
We shall describe the tangent space at the origin of $X(\sigma)\subset X=G/K$ for certain involutive automorphisms $\sigma$ for  
which $X(\sigma)$ is Hermitian symmetric and $\textrm{codim}_XX(\sigma)$ is minimum.   
In all but one case $\sigma$ is an element of the group $\Sigma$.     
In the case when $G$ is an exceptional group we shall merely compute the codimension as this is the only information 
that will be needed for our puposes. 
Again we shall use the conventions of the Planches in \cite{bourbaki} for labelling of the 
simple roots of $\Phi(\mathfrak{g}, \mathfrak{t})$. We shall also use the formula for 
$\text{codim} (X(\sigma))$ given in Remark \ref{codim}(i).

{\bf Type AIII}.  $\mathfrak{g}_0=\mathfrak{su}(p,n-p), 2p\le n$.    We have $\dim_\mathbb{C} X=p(n-p)$.  The simple non-compact root is $\epsilon_p-\epsilon_{p+1}$.  
Taking $\psi=\epsilon_{n-1}-\epsilon_n$, we see that for any positive root $\alpha=\sum a_{\alpha, \gamma}\gamma,
\gamma\in \Delta_{\mathfrak{g}}$, $a_{\alpha,\gamma}$ is either $0$ or $1$ and that $a_{\alpha, \psi}=1$ if and only 
if $\alpha=\epsilon_i-\epsilon_n$.  
Such a root $\alpha$ is non-compact if and only if $i\le p$. Thus there are 
exactly $p$ such roots and so the complex codimension of $T_0X(\sigma)\subset T_0X$ equals $p$.  
In this case $G(\sigma)=\SU(p,n-1-p)$ which embeds in $\SU(p,n-p)$ by fixing the standard basis element 
$\epsilon_n\in \mathbb{C}^n$. 
 
It is easily verified that if we take $\psi$ to be any other compact simple root, the (complex) codimension of 
$X(\sigma_\psi)\subset X$ is at least $2p$, except when $2p=n$ and $\psi=\epsilon_1-\epsilon_2,$
 in which $\textrm{codim}_XX(\sigma_\psi)=p$ again.  

{\bf Type BDI (rank=2)}. 
$\mathfrak{g}_0=\mathfrak{so}(2,n)$ with non-compact simple root $\psi_1=\epsilon_1-\epsilon_2$.\\ There are two cases 
to consider depending on the parity of $n$.

{\it Case 1}. Let $n=2p-1$ be odd, $p\ge 3$.  Consider $\psi:=\psi_p=\epsilon_p$. Then $\Phi^+_\n=
\{\epsilon_1\pm \epsilon_j, 1<j\le p\}\cup \{\epsilon_1\}$ and $\Phi^+(\sigma_\psi)=\{\epsilon_1\pm \epsilon_j, 
1<j\le p\}$.  Therefore 
there is exactly one non-compact root $\alpha$ for which the coefficient of $\psi$ is odd, namely, $\epsilon_1.$  
So $\textrm{codim}_XX(\sigma_\psi)=1.$   The Lie algebra $\mathfrak{g}(\sigma_\psi)$ modulo its radical is 
isomorphic to $\mathfrak{so}(2,n-1)$ with positive roots  $\{\epsilon_i\pm \epsilon_j\mid 1\le i< j\le p\}$.

{\it Case 2.}  Let $n=2p-2$ be even, $p\ge 4$.  Consider $\psi=\psi_2=\epsilon_2-\epsilon_3$.  The only 
non-compact roots 
in which the coefficient of $\psi_2$ is even are $\psi_1=\epsilon_1-\epsilon_2$ and the highest root $\alpha_0:=\epsilon_1+\epsilon_2$. 
It follows that 
$\textrm{dim}_XX(\sigma_\psi)=2$ and so $\textrm{codim}_XX(\sigma_\psi\circ\theta)=2$.  

This is the minimum codimension of a geometric cycle as we vary $\sigma \in \Sigma$.  However, with notations as in 
\S\ref{outerinvolutions}, we note that $\textrm{codim}_XX(\tau)=1$.

{\bf Type CI}.
$\mathfrak{g}_0=\mathfrak{sp}(n,\mathbb{R})$ with non-compact simple root $2\epsilon_n$.  
Take $\psi=\psi_1=\epsilon_1-\epsilon_2$.  Then $\Phi^+_\n=\{\epsilon_i+\epsilon_j, 1\le i<j\le n; 2\epsilon_j, 1\le j\le n\}$ 
and $\Phi^+(\sigma_\psi)_\mathfrak{n}=\{2\epsilon_j, 1\le j\le n; \epsilon_i+\epsilon_j, 1<i<j\le n  \}$.   Therefore 
$\textrm{codim}_XX(\sigma_\psi)=\#\Phi^+_\n-\#\Phi^+(\sigma_\psi)_\n=n-1.$ 
In this case $\mathfrak{g}(\sigma)$ modulo its radical is isomorphic to $\mathfrak{sl}(2,\mathbb{C})\times \mathfrak{sp}(n-1,\mathbb{C})$, where the roots of the $\mathfrak{sl}(2,\mathbb{C})$ factor are $\pm 2\epsilon_1$.  (In fact 
$\mathfrak{g}_0\cong \mathfrak{su}(2)\oplus \mathfrak{sp}(n-1,\mathbb{R})$ since $\pm \epsilon_1$ are compact roots.) 

{\bf Type DIII}
$\mathfrak{g}_0=\mathfrak{so}^*(2n)$.  The non-compact root is $\psi_n=\epsilon_{n-1}+\epsilon_n$ and 
$\Phi_\n^+=\{\epsilon_i+\epsilon_j\mid 1\le i<j\le n\}$.
Let $\psi=\psi_1=\epsilon_1-\epsilon_2$.  Then $\Phi^+(\sigma_\psi)_\n=\{ \epsilon_i+\epsilon_j\mid 1<i<j\le n\}$. 
So $\#\Phi^+_\n-\Phi^+(\sigma_\psi)_\n=\#\{\epsilon_1+\epsilon_j\mid 2\le j\le n\}=n-1$. Thus $\textrm{codim}_XX(\sigma_\psi)=n-1$.

{\bf Type EIII}
$\mathfrak{g}_0=\mathfrak{e}_{6,(-14)}$, with non-compact simple root $\psi_1$.  In this case $\#\Phi^+_\n=16$. 
Take $\psi=\psi_3$.  Using \cite[Planche V]{bourbaki}, we observe that among the non-compact positive roots,  
the coefficient of $\psi_3$ is at most $2$.   Among them
one has coefficient $0$ and five have coefficient $2$.
So $\textrm{codim}_XX(\sigma_\psi)=10$.  Hence $\textrm{codim}_XX(\sigma_\psi\theta)=6$.  A routine calculation 
shows that for {\it any} $\tau\in \Sigma$, the codimension, $\textrm{codim}_XX(\tau)$ is at least six.  

{\bf Type EVII}
$\mathfrak{g}_0=\mathfrak{e}_{7,(-25)}$.  The non-compact simple root is $\psi_7$.  $\#\Phi^+_\n=27$. 
Take $\psi=\psi_6$.  Again using Planche-VI in \cite{bourbaki}, we see that the coefficient of $\psi_6$ in any non-compact root is at most $2$ and that the number of non-compact 
roots in which $\psi_6$ occurs with coefficient $0, 2$ are respectively, $1, 10$.  
Thus $\#\Phi^+(\sigma_\psi)_\mathfrak{n}=
11$ and we have $\textrm{codim}_XX(\sigma_\psi)=16$.  Hence $\textrm{codim}_XX(\theta\sigma_\psi)=11$.  
This is the least possible codimension of $X(\tau)$ as $\tau$ varies in $\Sigma$ by a routine verification.  
The verification was made by direct computation as well as by using Python.  

Table \ref{c(X)} summarises the results obtained above.
The number $c(X)$ denotes the smallest number that arises 
as complex codimension of $X(\sigma)\subset X$  that has been constructed 
above.  In the table we have also indicated a $\sigma$ which 
realizes this number.  Except when $G=\SO_0(2,2p-2)$, $\sigma$ belongs to $\Sigma$.

\begin{table}[h]
\centering
\begin{tabular}{|c|c|c|c|}
\hline 
Type & $\mathfrak{g}_0$  & $\sigma$&$c(X)$\\
\hline
AIII& $\mathfrak{su}(p,n-p)$ & $\sigma_{\epsilon_{n-1}-\epsilon_n}$   &$p$\\
&$1\le p\le n/2$ &   &\\
\hline
BD I (rank ~2) & $\mathfrak{so}(2,n)$ &   &  \\
& $n=2p-1$& $\sigma_{\epsilon_p}$ & $1$\\
& $n=2p-2$& $\tau (\textrm{see~} \S\ref{outerinvolutions})$ &$1$\\
\hline 
CI & $\mathfrak{sp}(n,\mathbb{R})$ &$\sigma_{\epsilon_1-\epsilon_2}$ &  $n-1$ \\
\hline
DIII&$\mathfrak{so}^*(2n)$ & $\sigma_{\epsilon_1-\epsilon_2}$ &$n-1$\\
\hline
EIII & $\mathfrak{e}_{6,(-14)}$ &$\sigma_{\psi_3} \circ \theta$      &$6$ \\
\hline
EVII &$\mathfrak{e}_{7,(-25)}$ & $\sigma_{\psi_6}\circ \theta$ &$11$\\
\hline
\end{tabular}
~\\
\caption{The values of $c(X)$ and involutions yielding analytic special cycles of minimum codimension.}
\label{c(X)}
\end{table}


\section{$\theta$-stable parabolic subalgebras}\label{parabolics} \label{thetastableparabolics}

Let $G$ be a linear connected non-compact simple Lie group and $K\subset G$ a maximal compact subgroup of $G$. Denote by $\theta$ the Cartan involution that fixes $K$.  We denote its 
differential $\mathfrak{g}_0\to \mathfrak{g}_0$ and also its  complexification by the same symbol  
$\theta$. We have the Cartan decomposition $\mathfrak{g}_0=\mathfrak{k}_0\oplus \mathfrak{p}_0$ 
where $\mathfrak{p}_0$ is the $(-1)$-eigenspace of $\theta$. Complexification yields $\mathfrak{g}=\mathfrak{k}\oplus \mathfrak{p}$.  We assume that $G/K$ is Hermitian symmetric.  Thus the centre 
of $K$, denoted $Z(K)$, is isomorphic to $U(1)$ and the (complex) rank of $G$ equals the rank of $K$.  We fix a maximal 
torus $T\subset K$.   
  We denote by $\Phi$ the set of 
roots of $(\mathfrak{g},\mathfrak{t})$; $\Phi_\mathfrak{k}, \Phi_\n\subset \Phi$ denote the set of compact, 
respectively, non-compact roots.   We fix a positive root system for $(\mathfrak{g}, \mathfrak{t})$ 
such that the set of simple roots $\Delta_\mathfrak{g}$ has exactly one non-compact root; $\Phi^+, \Phi^+_\mathfrak{k}$ denote
the set of positive roots of $\mathfrak{g}, \mathfrak{k}$ respectively and $\Phi^+_\n, \Phi^-_\n$ the set of 
positive, resp. negative, non-compact roots.  Then $\Delta_\mathfrak{k}:=\Delta_\mathfrak{g}\cap \Phi^+_\mathfrak{k}$ is the simple 
roots for the positive system $\Phi_\mathfrak{k}^+$.  The fundamental weight corresponding to a simple root $\psi\in \Delta_\mathfrak{g}$ will be 
denoted $\varpi_\psi$.  

If $\mathfrak{s}$ affords a $\mathfrak{t}$-representation, especially when $\mathfrak{s}\subset \mathfrak{g}$, we shall denote by 
$\Phi(\mathfrak{s})$ the multiset of non-zero weights of $\mathfrak{s}$ and by $|\Phi(\mathfrak{s})|$ the sum of elements 
of $\Phi(\mathfrak{s})$ each appearing as many times as its multiplicity.

The Killing form on $\mathfrak{g}$ restricted to $i\mathfrak{t}_0\subset \mathfrak{t}$ is an 
inner product.  This in turn defines an inner product on $i\mathfrak{t}_0^*$, which will be denoted $(\cdot,\cdot)$.

The tangent space $\mathfrak{p}_0$ of $G/K$ at the origin is a complex vector space and so we have a 
decomposition $\mathfrak{p}:=\mathfrak{p}_0\otimes_\mathbb{R}\mathbb{C}=\mathfrak{p}_+\oplus \mathfrak{p}_-$ as the holomorphic and anti-holomorphic tangent spaces at the identity coset of $G/K$ where 
$\mathfrak{p}_+=  \oplus_{\alpha\in \Phi^+_\n}\mathfrak{g}_\alpha, 
\mathfrak{p}_-=  \oplus_{\alpha\in \Phi^+_\n}\mathfrak{g}_{-\alpha}$.

Denoting by $~\bar{}~$ the complex 
 conjugation of $\mathfrak{g}=\mathfrak{g}_0\oplus i\mathfrak{g}_0$ with respect to $\mathfrak{g}_0$, 
 we recall that a $\theta$-stable parabolic subalgebra $\mathfrak{q}$ of $\mathfrak{g}_0$ is a parabolic subalgebra contained in  $\mathfrak{g}$ such that (a) $\theta(\mathfrak{q})=\mathfrak{q}$, and, (b) $\mathfrak{q}\cap \overline{\mathfrak{q}}=:\mathfrak{l}$ is a Levi subalgebra of $\mathfrak{q}$.   If $k\in K$ and $\mathfrak{q}$ is a $\theta$-stable parabolic subalgebra, then so is $\text{Ad}(k)(\mathfrak{q})$.   
Then $\mathfrak{l}_0:=\mathfrak{l}\cap \mathfrak{g}_0$ contains $\text{Lie}(T_1)$ for some maximal torus $T_1\subset K$.
Conjugating by an element of $K$ if required, we may assume that $\mathfrak{t}_0\subset \mathfrak{l}_0$ so that $\mathfrak{t}\subset \mathfrak{q}$.    
Let  $\mathfrak{u}$ be the nilradical of $\mathfrak{q}$; thus we have $\mathfrak{q}=\mathfrak{l}\oplus \mathfrak{u}$.  

The $\theta$-stable parabolic subalgebras of $\mathfrak{g}_0$ containing $\mathfrak{t}$ are constructed 
as follows: Let $x\in i\mathfrak{t}_0$.  Note that the roots of $(\mathfrak{g},\mathfrak{t})$ take real values 
on $i\mathfrak{t}_0$.  Let $\mathfrak{q}_x:=\mathfrak{t}+\sum_{\alpha(x)\ge 0, \alpha\in \Phi} \mathfrak{g}_\alpha $. Then $\mathfrak{q}_x=\mathfrak{l}_x\oplus \mathfrak{u}_x$ is a $\theta$-stable subalgebra of $\mathfrak{g}_0$ where  the nilradical of $\mathfrak{q}_x$ equals $\mathfrak{u}_x=\oplus \sum_{\alpha(x)>0}\mathfrak{g}_\alpha$ and the Levi subalgebra $\mathfrak{l}_x$ equals  
$\mathfrak{t}+\sum_{\alpha(x)=0}\mathfrak{g}_\alpha$.  We denote by $\Phi_x$ the roots of $(\mathfrak{l},\mathfrak{t})$.  
Every $\theta$-stable subalgebra $\mathfrak{q}$ that contains $\mathfrak{t}$ arises as $\mathfrak{q}_x$ for some $x\in i\mathfrak{t}_0$.   Moreover, fixing a positive system for $(\mathfrak{k},\mathfrak{t})$ we may assume, without loss of 
generality that, $\alpha(x)\ge 0$ for all $\alpha\in\Phi^+_\mathfrak{k}$.    

Recall that if, $\lambda$ is any element of $i\mathfrak{t}_0^*$, then $h_\lambda\in i\mathfrak{t}_0$ is the  
unique element such that $\lambda(H)=(H,h_\lambda) ~\forall H\in \mathfrak{t}^*$ and we have $\alpha(h_\lambda)=(\lambda, \alpha)~\forall \alpha\in 
i\mathfrak{t}^*_0$.    
When $x=h_\lambda$ we often denote 
$\mathfrak{q}_x, \mathfrak{l}_x, \mathfrak{u}_x$ by $\mathfrak{q}_\lambda, \mathfrak{l}_\lambda,\mathfrak{u}_\lambda$ respectively.

We choose a positive root system for $(\mathfrak{l}\cap\mathfrak{k},\mathfrak{t})$ and extend it to a positive root system $\Phi^+_\mathfrak{k}(x)$ 
for $(\mathfrak{k},\mathfrak{t})$ such that the set $\Phi(\mathfrak{u}\cap \mathfrak{k})\subset \Phi^+_{\mathfrak{k}}(x)$ where $\Phi(\mathfrak{u}\cap\mathfrak{k})$ denotes the $\mathfrak{t}$-weights of $\mathfrak{u}\cap \mathfrak{k}$.   One has an irreducible unitary representation $(\mathcal{A}_{\mathfrak{q}}, A_\mathfrak{q})$ of $G$ with trivial infinitesimal character  
such that (i) the $(\mathfrak{g},K)$-module $A_{\mathfrak{q},K}$ has an irreducible $K$-submodule $V$ with highest weight 
(with respect to $\Phi^+_\mathfrak{k}(x)$) equal to the sum $|\Phi(\mathfrak{u}\cap \mathfrak{p})|=
\sum_{\alpha\in \Phi(\mathfrak{u\cap p})}\alpha$, (ii) the $K$-type of $V$ occurs in $A_{\mathfrak{q},K}$ with multiplicity one, i.e., 
$\hom_K(V,A_{\mathfrak{q},K})\cong \mathbb{C}$, and, (iii) any other $K$-type that occurs in $A_{\mathfrak{q},K}$ has 
highest weight of the form $|\Phi(\mathfrak{u}\cap \mathfrak{p})|+\sum_{\gamma\in \Phi(\mathfrak{u}\cap\mathfrak{p})}
a_\gamma\gamma$ with $a_\gamma\ge 0$.   If $\mathfrak{q}':=\mathfrak{q}_{x'}=\mathfrak{l}_{x'}\oplus 
\mathfrak{u}_{x'},x'\in i\mathfrak{t}_0,$ is another $\theta$-stable parabolic 
subalgebra of $\mathfrak{g}_0$,
then 
$A_{\mathfrak{q},K}$ is unitarily equivalent to $A_{\mathfrak{q}',K}$ 
if and only if $\mathfrak{u}_x\cap\mathfrak{p}=\mathfrak{u}_{x'}\cap \mathfrak{p}$.   See \cite{sriba} for a more general statement.  It is a result due to Harish-Chandra that two irreducible {\it unitary} representation of $G$ are unitarily equivalent 
if and only if their spaces of smooth $K$-finite vectors are isomorphic as $(\mathfrak{g},K)$-modules.  (See \cite[Ch. IX, \S 1]{knapp}.)
In particular, $|\Phi(\mathfrak{u}_\lambda\cap \mathfrak{p})|=|\Phi(\mathfrak{u}_\mu\cap\mathfrak{p})|$ for two $\theta$-stable parabolic subalgebras 
$\mathfrak{q}_\lambda, \mathfrak{q}_\mu$ of $G$ if and only if 
$\mathcal{A}_{\mathfrak{q}_\lambda}$ and $\mathcal{A}_{\mathfrak{q}_\mu}$ are unitarily equivalent $G$-representations.

The group $L=\{g\in G\mid \text{Ad}(g)(\mathfrak{q})=\mathfrak{q}\}$ is a connected reductive closed Lie subgroup of 
$G$ with Lie algebra $\mathfrak{l}_0=\mathfrak{l}\cap \mathfrak{g}_0.$ 
Moreover, $T\subset L$ and $[L,L]\cap K $ is a maximal compact subgroup of the 
semisimple Lie group $[L,L]$.  
Let $Y_\mathfrak{q}$ denote the compact dual of $[L,L]/(K\cap[L, L])$.  
It turns out that $Z(K)\subset [L,L]$ and so $Y_\mathfrak{q}$ is {\it Hermitian} symmetric.
(See Proposition \ref{hermitian} below.) 

It is known that if $(\pi,H_\pi)$ is an irreducible unitary representation of 
$G$ such that $H^*(\mathfrak{g},K;H_{\pi,K})$ is non-zero, then $\pi$ is unitarily equivalent to 
$\mathcal{A}_\mathfrak{q}$ for some $\theta$-stable parabolic subalgebra $\mathfrak{q}$.  
Also $H^r(\mathfrak{g},K;A_{\mathfrak{q},K})$ 
is isomorphic to $H^{r-R(\mathfrak{q})}(Y_\mathfrak{q};\mathbb{C})$, where $R(\mathfrak{q})=\dim_\mathbb{C}\mathfrak{u}\cap \mathfrak{p}$. In fact a more refined statement is valid as we shall now describe. 

Let $R_+(\mathfrak{q})=\dim_\mathbb{C} \mathfrak{u}\cap \mathfrak{p}_+, R_-(\mathfrak{q})=
\dim_\mathbb{C}\mathfrak{u}\cap \mathfrak{p}_-$ so that $R(\mathfrak{q})=R_+(\mathfrak{q})+R_-(\mathfrak{q})$.  Then 
$H^{p,q}(\mathfrak{g},K;A_{\mathfrak{q},K})\cong 
H^{p-R_+(\mathfrak{q}),q-R_-(\mathfrak{q})}(Y_\mathfrak{q};\mathbb{C})$.   (See \cite{vz}.)
Since $Y_\mathfrak{q}$ is Hermitian symmetric, we see that $H^{p,q}(\mathfrak{g},K;A_{\mathfrak{q},K})=0$ unless $p-q=R_+(\mathfrak{q})-R_-(\mathfrak{q})$. See \cite{gh}.
We refer to $(R_+(\mathfrak{q}),R_-(\mathfrak{q}))$ as the {\it Hodge type} of $\mathfrak{q}$ and to $R(\mathfrak{q})=R_+(\mathfrak{q})+R_-(\mathfrak{q})$ as the {\it degree} of $\mathfrak{q}$.

The $(\mathfrak{g},K)$-module $A_{\mathfrak{q},K}$ was first constructed by Parthasarathy \cite{parthasarathy}. 
Vogan and Zuckerman \cite{vz} and Vogan \cite{vogan} gave a construction via cohomological induction and proved 
that they are unitarizable.   We refer the reader to the paper \cite{vogan97} for a very readable account of 
the basic properties of $\mathcal{A}_\mathfrak{q}$.  For basic representation theory of semisimple Lie groups  
we refer the reader to Knapp's book \cite{knapp}.  For the cohomology of $(\mathfrak{g},K)$-modules 
and its relation to the cohomology of lattices in Lie groups, see \cite{borel-wallach}.


\subsection{The Levi subalgebras of $\theta$-stable subalgebras of $\mathfrak{g}_0$.} \label{levi}
Having fixed a positive system of roots for $(\mathfrak{g},\mathfrak{t})$,  we have a partial order on the set of roots where $\alpha\ge \beta$ if $\alpha-\beta$ is a non-negative linear combination of simple roots.  Let $\mathfrak{q}=\mathfrak{q}_x$ where $x\in i\mathfrak{t}_0$ is such that $\gamma(x)\ge 0$ for all compact roots $\gamma\in \Phi^+_\mathfrak{k}.$   Write $\Phi_x=\{\alpha\in \Phi\mid \alpha(x)=0\}$.  It is clear that $\Phi_x$ is the set of roots of $(\mathfrak{l},\mathfrak{t})$.  
Our assumption on $x$ that $\psi(x)\ge 0$ for all $\psi\in \Phi^+_\mathfrak{k}$ implies that if $\alpha+\beta\in \Phi_x, \alpha,\beta\in \Phi^+_\mathfrak{k}$, then $\alpha,\beta\in \Phi_x$.  
 If $\Phi_x\subset \Phi_\mathfrak{k}$, then $[\mathfrak{l},\mathfrak{l}]\subset \mathfrak{k}$ and $Y_\mathfrak{q}$ 
is reduced to a point.  So assume that $\Phi_x\cap \Phi_\n\ne \emptyset.$

Let $\mathcal{C}=\mathcal{C}(x)\subset \Phi_x\cap \Phi^+_\n$ be the set of all positive non-compact roots such that 
there is no positive non-compact root $\beta\in\Phi_x$ where $\beta<\alpha$.  For each $\alpha\in \mathcal{C}$, 
let $\Delta_\alpha\subset \Phi^+$ denote the set consisting of $\alpha$ 
and all compact simple roots $\psi$ such that there exists a positive non-compact root $\beta\in \Phi_x$ where $\beta>\alpha$ and 
$(\varpi_\psi,\beta-\alpha)\ne 0.$   
Also let $\Phi_\alpha\subset \Phi_x$ be the set consisting of all roots $\beta \in \Phi_x$ 
in the span of $\Delta_\alpha$.   Then $(\Phi_\alpha,\Delta_\alpha)$ is a reduced root system.
From the definition of $\mathcal{C}$, it follows that if $\gamma\in \Phi_\alpha$ is a compact root then 
$\gamma$ is in the span of $\Delta_\alpha\cap \Phi_\mathfrak{k}$.  

  The Lie subalgebra of $\mathfrak{l}$ generated by the root spaces $\mathfrak{g}_\gamma, \gamma\in \Phi_\alpha,$ will be denoted by $\mathfrak{l}_\alpha$.  It is clear that $\mathfrak{l}_\alpha\subset 
[\mathfrak{l},\mathfrak{l}]$.  
Denote the subgroup of $[L,L]$ corresponding to $\mathfrak{l}_{\alpha,0}:=\mathfrak{l}_0\cap \mathfrak{l}_\alpha$  by $L_\alpha$ and the group 
$K\cap L_\alpha$ by $K_\alpha$ for $\alpha\in \mathcal{C}$.   The set $\Delta_\alpha$ has exactly one non-compact root, namely $\alpha$. Since the coefficient of 
this root in any non-compact root of $\Phi_\alpha$ is $\pm 1$,  $(\Phi_\alpha,\Delta_\alpha)$ is a Borel-de Siebenthal root system 
of a Hermitian symmetric pair $(L_\alpha,K_\alpha)$ of non-compact type. (See \cite{bds}.) 
Thus $L_\alpha/K_\alpha$ is an irreducible Hermitian symmetric space of non-compact type.

We let $\Phi_\textrm{c}=\Phi_x\setminus (\cup_{\alpha\in \mathcal{C}}\Phi_\alpha)$.  Then $\Phi_\textrm{c}$ 
consists entirely of compact roots.
The Lie subalgebra of $\mathfrak{l}$ generated by $\mathfrak{g}_\gamma, \gamma\in \Phi_\textrm{c}$ is an ideal 
$\mathfrak{l}_\textrm{c}$ in $\mathfrak{[l,l]}$ and the subgroup of $[L,L]$ corresponding to 
$\mathfrak{l}_0\cap\mathfrak{l}_\textrm{c}$ is a maximal compact normal subgroup, which we denote by $K_\textrm{c}$.  
A maximal compact subgroup of $[L,L]$ is the product 
$K_\textrm{c}.\prod_{\alpha\in \mathcal{C}}K_\alpha$.   We summarise below the above discussion.

\begin{proposition} \label{hermitian}
With the above notation, the simple ideals of $[\mathfrak{l,l}]$ are $\mathfrak{l}_\textrm{c}$ and $\mathfrak{l}_\alpha,\alpha\in \mathcal{C}$.  
The homogeneous space $L_\alpha/K_\alpha$ is an irreducible globally Hermitian symmetric space of non-compact type.  Hence $[L,L]/(K\cap [L,L])=\prod_{\alpha\in \mathcal{C}} L_\alpha/K_\alpha$ is Hermitian symmetric. \hfill $\Box$
\end{proposition}

We have the following lemma which implies, in particular, that $\#\mathcal{C}$ does not exceed the real rank of $G$.

\begin{lemma}
The set $\mathcal{C}\subset \Phi^+_\n$ is a set of strongly orthogonal roots.  
\end{lemma}
\begin{proof}  Let $\alpha,\beta\in \mathcal{C}$. 
Since the sum of two positive non-compact roots is never a root (as $\mathfrak{p}_+$ is an abelian subalgebra) it suffices 
to show that $\beta-\alpha$ is not a root.  Indeed if $\beta-\alpha=:\kappa$ is a root, it has to be a compact root, which we assume is positive.  Now $\beta,\alpha\in \Phi_x$ implies that $\kappa(x)=0$.   Therefore $\beta=\alpha+\kappa$ implies that $\beta\in \Phi_\alpha$ and hence $\beta\notin\mathcal{C}$, a contradiction.  Since $\beta\pm \alpha$ are not roots, we must have $(\alpha,\beta)=0$.  
\end{proof}

\begin{remark} \label{rootsofu}
{\em 
(i)   Since $\theta$ restricts to the identity on $\mathfrak{t}$, it is clear that $\theta(\mathfrak{l})=\mathfrak{l}$.  In fact 
$\theta|_\mathfrak{l}$ is the $\mathbb{C}$-linear extension of $\theta|_{\mathfrak{l}_0}$.  
Moreover, $\theta|_{[\mathfrak{l}_0,\mathfrak{l}_0]}$ is a Cartan involution of 
$[\mathfrak{l}_0,\mathfrak{l}_0]$ and $\theta|_{\mathfrak{l}_{\alpha,0}}$ is a 
Cartan involution of $\mathfrak{l}_{\alpha,0}$.

(ii)  Recall that $R(\mathfrak{q})=\dim_\mathbb{C}\mathfrak{u\cap p}$.  If $\beta\in \Phi_\n\setminus \Phi_x$, then $\beta(x)\ne 0$ and so exactly one of the roots $\beta, -\beta$ is a weight of $\mathfrak{u\cap p}$.  It follows that 
$R(\mathfrak{q})=(1/2)\#(\Phi_\n\setminus\Phi_x)=
\dim_\mathbb{C}G/K-\dim_\mathbb{C} [L,L]/(K\cap [L,L]).$ 
In particular, if $R_+(\mathfrak{q})=R_-(\mathfrak{q})$, then $R_+=(1/2) \#(\Phi^+_\n\setminus \Phi_x^+)$.

(iii) 
Let $\beta>\alpha$ where $\beta,\alpha\in \Phi^+_\textrm{n}$.  
Since $\psi(x)\ge 0$ for all $\psi\in \Delta_\mathfrak{k}$, we have $\beta\in \Phi(\mathfrak{u}_x\cap\mathfrak{p}_+)$ if $\alpha\in 
\Phi(\mathfrak{u}_x\cap \mathfrak{p}_+)$.
In particular there exists a unique set of pairwise non-comparable positive roots $\xi_1,\ldots, \xi_r\in \Phi^+_\textrm{n}$  (depending on $x$) such that 
$\Phi(\mathfrak{u}_x\cap \mathfrak{p}_+)=\cup_{1\le i\le r}\{\eta\in \Phi^+_\textrm{n}\mid \eta\ge \xi_i\}$.  
An analogous statement holds for $\Phi(\mathfrak{u}_x\cap\mathfrak{p}_-)$.   
Write $x=h_\lambda$.
If $\psi$ is a compact simple root  
such that $\xi_j-\psi$ is a root for some $j$ and $\xi_j-\psi\in\Phi(\mathfrak{l}_x\cap \mathfrak{p}_+)$,  then 
$(\lambda,\psi)=(\lambda,\xi_j)>0$.
If $\alpha\in\Phi^+_n$ and $\psi$ is a simple compact root such that $\alpha,\alpha+\psi\in\Phi(\mathfrak{l}_x\cap\mathfrak{p}_+)$, then $(\lambda,\psi)=0$. 
These elementary observations will be used in classifying $\theta$-stable parabolic subalgebras 
of $\mathfrak{g}_0$ with prescribed Hodge type, particularly in the case of exceptional Lie algebras of type EIII and EVII.
}
\end{remark}

The Weyl group $W(K,T)\cong W(\mathfrak{k},\mathfrak{t})=:W_\mathfrak{k}$ acts on $i\mathfrak{t}_0$ and is generated by the set $S$ of simple reflections $s_\gamma, \gamma \in \Delta_\mathfrak{k}$.   We have the length function defined on $W_\mathfrak{k}$ with respect to $S$.    We will denote by $w_0^\mathfrak{k}$ (or more briefly $w_0$), 
the longest element of $W_\mathfrak{k}$.   Recall that   
$w_0(\Delta_\mathfrak{k})=-\Delta_\mathfrak{k}$ and that $w_0^2=1$.  
We denote by $\iota$ the Weyl involution $-w_0$ on $\mathfrak{t}, i\mathfrak{t}_0$ or on their duals.

\begin{lemma}  \label{purity}
Suppose that $x\in i\mathfrak{t}_0$ satisfies the condition that $\gamma(x)\ge 0~\forall \gamma\in \Phi_\mathfrak{k}^+$, then $\iota(x)$ also satisfies this condition. Moreover 
$(R_+(\mathfrak{q}_{\iota(x)}),R_-(\mathfrak{q}_{\iota(x)}))=(R_-(\mathfrak{q}_x),R_+(\mathfrak{q}_x))$.  In particular, if 
$x=\iota (x)$, then $R_+(\mathfrak{q}_x)=R_-(\mathfrak{q}_x)$.  
\end{lemma}
\begin{proof} 
Note that $\iota$ yields a 
bijection of $\Phi^+_\mathfrak{k}$ onto itself.  So, if $\gamma\in \Phi^+_\mathfrak{k},$ then $\iota(\gamma)\in \Phi^+_\mathfrak{k}$ and we have $\gamma(\iota(x))=\iota(\gamma)(x)\ge 0$.  This proves the first assertion.  

Since $\mathfrak{p}_+$ and $\mathfrak{p}_-$ are irreducible representations of $K$ which are dual to each other, 
we have $\iota (\Phi^+_\n)=-\Phi^+_\n=\Phi^-_\n$.   
We need only show that $\alpha\in \Phi^+_\n$ is a weight of $\mathfrak{u}\cap \mathfrak{p}_+$ if and only if 
$\iota (\alpha)\in \Phi^-_\n$ is a weight of $\mathfrak{u}_{\iota(x)}\cap \mathfrak{p}_-$.    This is immediate from 
the observation $\iota (\alpha)(\iota(x))=w_0(\alpha)(w_0(x))=\alpha(w_0^{-1}w_0(x))=\alpha(x)$, and the lemma 
follows. 
\end{proof}


\subsection{The $\theta$-stable parabolic subalgebras of type $(p,p)$} \label{pptypeparabolics}
As at the beginning of \S \ref{thetastableparabolics},  $G$ is a linear connected simple Lie group with finite centre, $K\subset G$ is a maximal compact subgroup and $X=G/K$ is an irreducible Hermitian symmetric space of non-compact 
type.

We shall classify $\theta$-stable parabolic subalgebras $\mathfrak{q}$ of $\mathfrak{g}_0$ such that $R_+(\mathfrak{q})=R_-(\mathfrak{q})\le c(X)$, the (complex) 
codimension 
of a complex analytic geometric cycle $X(\sigma)\subset X=G/K$ constructed in \S\ref{involutions}; see Table \ref{c(X)}.  
(In most cases $c(X)$ is the smallest such 
positive integer.  However this property of $c(X)$ will not be needed.)

Recall that if $x=h_\lambda\in i\mathfrak{t}_0$,  then $\mathfrak{q}_\lambda$ stands for $\mathfrak{q}_x$.  Note that 
$\psi(h_\lambda)\ge 0$ for all $\psi\in \Phi^+(\mathfrak{k})$ if and only if $\lambda$ is $\mathfrak{k}$-dominant.

For any $p\ge 1$, let $N(p)$ be the number of 
$\theta$-stable parabolic subalgebras $\mathfrak{q}=\mathfrak{q}_\lambda$  of 
$\mathfrak{g}_0$ (where $\lambda\in i\mathfrak{t}_0^*$ is in the dominant Weyl chamber) such that $\mathfrak{q}$ is of 
Hodge type $(p,p)$, i.e., $R_\pm(\mathfrak{q})=p$.

In this section we shall determine $N(p)$ for $p\le c(X)$.   We denote by $r=r(\mathfrak{g}_0)$ the smallest positive 
integer such that 
$N(r)\ge 1$.  Our results are summarised in Table \ref{valuesofr0}. 

\begin{table}[h]
\centering
\begin{tabular}{|c|c|c|c|c|}
\hline
Type & $\mathfrak{g}_0$&$c(X)$& $r(\mathfrak{g}_0)$ & $\sum_{1\le p\le c(X)}N(p)$\\
\hline
AIII &$\mathfrak{su}(p,q), p<q-1$& $p$& $p$& $1$ \\
&~~~$p=q-1$&$p$&$p$&$3$\\
&~~~$p=q$&$p$&$p-1$& $4$\\
 \hline
 BDI& $\mathfrak{so}(2,2n-1)$& $1$ & $1$ &$1$\\
 &$\mathfrak{so}(2,2n-2)$ &$1$&$1$&$1$ \\
 \hline
 CI & $\mathfrak{sp}(n,\mathbb{R}), n\ne 4$ & $n-1$ &$n-1$&$1$\\
 \hline
 DIII & $\mathfrak{so}^*(2n),n\ge 9$ & $n-1$& $n-2$& $1$\\
 \hline
 EIII &$\mathfrak{e}_{6,(-14)}$ & $6$ & $4$& $4$\\
 \hline
 EVII & $\mathfrak{e}_{7,(-25)}$ & $11$ & $6$&$6$\\
 \hline
\end{tabular}
\caption{Values of $r(\mathfrak{g}_0).$}
\label{valuesofr0}
\end{table}

We shall proceed with the task of the classification in each type.  

\subsection*{Type AIII}
Let $G=\SU(p,q), p\le q$.  
The set $\Phi^+(\mathfrak{g})$ of positive roots of $\mathfrak{g}=\mathfrak{sl}(p+q,\mathbb{C})$ equals 
$\{\epsilon_i-\epsilon_j\mid 1\le i<j\le p+q\}$.  The non-compact simple root is $\epsilon_{p}-\epsilon_{p+1}$ and the set  
of positive non-compact roots equals $\Phi_\n^+=\{\epsilon_i-\epsilon_j\mid 1\le i\le p, p+1\le j\le p+q\}$.  Set $n:=p+q-1,$ the 
rank of $\mathfrak{g}$.   
We regard $i\mathfrak{t}_0^*$ as the 
subspace of the Euclidean space 
$\mathbb{R}^{p+q}$ where the sum of the coordinates is equal to zero.  (Cf. \cite[Planche I]{bourbaki}.)  
Thus $\lambda=\sum_{1\le i\le p+q}a_i\epsilon_i\in i\mathfrak{t}_0^*$ 
if and only if $\sum_{1\le i\le n} a_i=0$.  It is convenient to set 
$\epsilon_0=\sum_{1\le i\le p+q}\in \mathbb{R}^{p+q}$.  
(Of course $\epsilon_0\notin i \mathfrak{t}_0^*$.)

Let $\omega:=\epsilon_{p+1}-\epsilon_{p+q}$.    
For any $\gamma\in \Phi^+_\mathfrak{k}$ we see that $(\gamma,\omega) \ge 0$. 
We have $(\epsilon_k-\epsilon_{p+1},\omega) =-1$ and $(\epsilon_k-\epsilon_{p+q},\omega)=1$ for $k\notin\{p+1,p+q\}.$  Also $(\epsilon_i-\epsilon_j,\omega)=0$ if $i,j\notin\{p+1,p+q\}$. Hence, considering $\mathfrak{q}_\omega=\mathfrak{l}_\omega\oplus \mathfrak{u}_\omega$, we have 
$\Phi(\mathfrak{u}_\omega\cap\mathfrak{p}_-)
=\{\epsilon_{p+1}-\epsilon_k\mid 1\le k\le p\}$, $\Phi(\mathfrak{u}_\omega\cap\mathfrak{p}_+)
=\{\epsilon_{k}-\epsilon_{p+q}\mid 1\le k\le p\}$ and so 
$R_\pm(\mathfrak{q}_\omega)=p$.
Also $|\Phi(\mathfrak{u}_\omega\cap \mathfrak{p})|=p(\epsilon_{p+1}-\epsilon_{p+q})=p\omega$.      

Similarly, if $\mu:=\epsilon_1-\epsilon_p$, then $(\gamma,\mu)\ge 0$ for all $\gamma\in \Phi_\mathfrak{k}^+$, 
$R_+(\mathfrak{q}_\mu)=R_-(\mathfrak{q}_\mu)=q, \Phi(\mathfrak{u}_\mu\cap \mathfrak{p})=\{\epsilon_1-\epsilon_k, 
\epsilon_k-\epsilon_p\mid p+1\le k\le p+q\}$ and $|\Phi(\mathfrak{u}_\mu\cap \mathfrak{p})|=q(\epsilon_1-\epsilon_p)=q\mu$.  Since 
$|\Phi(\mathfrak{u}_\mu\cap \mathfrak{p})|\ne |\Phi(\mathfrak{u}_\omega\cap \mathfrak{p})|$, the highest weights of the lowest 
 $K$-type in $A_{\mathfrak{q}_\omega,K}$ and $A_{\mathfrak{q}_\mu,K}$ are not equal. Hence we conclude that 
 $\mathcal{A}_{\mathfrak{q}_\omega}, \mathcal{A}_{\mathfrak{q}_\mu}$ are not unitarily equivalent.  {\it In particular, 
 when $p=q, $ we have $N(p)\ge 2$.}

Let $q=p+1$ and let  $\kappa:=p\epsilon_1+q\epsilon_{p+1}-\epsilon_0\in i\mathfrak{k}^*_0, \nu=\epsilon_0-p\epsilon_p-q\epsilon_{p+q}$ which are $\mathfrak{k}$-dominant.  
By a straightforward computation $\Phi^+(\mathfrak{u}_{\kappa}\cap \mathfrak{p}_+)=\{\epsilon_1-\epsilon_j\mid p+1<j\le q\}$, 
so that $R_+(\mathfrak{q}_{\kappa})=q-1=p$.  Also $\Phi^+(\mathfrak{u}_{\kappa}\cap \mathfrak{p}_-)= 
\{\epsilon_{p+1}-\epsilon_j\mid 1\le j\le p\}$, and consequently $R_-(\mathfrak{q}_\kappa)=p$.  Observe that 
$|\Phi(\mathfrak{u}_\kappa\cap\mathfrak{p})|
=p\epsilon_1+(p+1)\epsilon_{p+1}-\epsilon_0=\kappa$.    Similarly $R_\pm(\mathfrak{q}_{\nu})=p$ and 
$|\Phi(\mathfrak{q}_{\nu}\cap \mathfrak{p})|=\nu$.
Since $p\omega, \kappa,\nu$ are pairwise distinct, the corresponding 
representations $\mathcal{A}_{\mathfrak{q}_\lambda}, \lambda=\omega,\kappa, \nu,$ of $G$ are pairwise inequivalent.  Hence $N(p)\ge 3$ in this case.

Suppose that $p=q$.   We now show that $N(p-1)\ge 2$.  Indeed, the 
$\xi:=p(\epsilon_1+\epsilon_{p+1})-\epsilon_0$ and $\eta:=\epsilon_0-p(\epsilon_p+\epsilon_{p+q})$ are $\mathfrak{k}$-dominant weights.   A straightforward computation shows that  
$R_\pm(\mathfrak{q}_\xi)=R_\pm(\mathfrak{q}_\eta) 
=p-1$ and $|\Phi(\mathfrak{u}_\xi\cap \mathfrak{p})|=\xi, |\Phi(\mathfrak{u}_\eta\cap \mathfrak{p})|= \eta$ and hence 
$\mathcal{A}_{\mathfrak{q}_\xi}, \mathcal{A}_{\mathfrak{q}_\eta}$ are  
inequivalent unitary representations.    Thus $N(p-1)\ge 2.$
 
We have the following theorem:  

\begin{theorem} \label{aiii}    We keep the above notations.  Let $G=\SU(p,q),1\le p\le q, q\ge 5$.  
Any unitary representation $\mathcal{A}_\mathfrak{q}$ having Hodge type $(r,r)$ where $1\le r\le p$ is 
unitarily equivalent to one of the representations $\mathcal{A}_{\mathfrak{q}_\lambda}$ where 
$\lambda\in \{\omega, \mu, \nu, \xi,\eta,\kappa\}$.   In particular, 
the number $N(r)$ of such unitary representations   
up to unitary equivalence, is as follows:  (i) if $p<q-1, N(p)=1$; if $p=q-1$, 
$N(p)=3$; if $p=q, N(p)=2$. (ii)  If  $1\le r<p$, then $N(r)=0$, except when $r+1=p=q$ in which case $N(p-1)
=2$.  
\end{theorem}
\begin{proof}    
When $p=1$ the statement is easily seen to be true.  So we shall assume that $p\ge 2$.  

 Let $\lambda=\sum_{1\le i\le p+q} a_j\epsilon_j\in i\mathfrak{t}^*_0$ be a non-zero 
 $\mathfrak{k}$-dominant weight.  Thus $(\gamma,\lambda)\ge 0~\forall \gamma\in \Phi_\mathfrak{k}^+,$ which implies 
 that $a_1\ge \cdots \ge a_p$ and $a_{p+1}\ge \cdots \ge a_{p+q}$.  
Also $\sum_{1\le j\le p+q} a_j=0$ as $\lambda\in i\mathfrak{t}_0^*$.  
Suppose that $0<R_+(\mathfrak{q}_\lambda)=R_-(\mathfrak{q}_\lambda)\le p$.   We will show that $R_\pm(\mathfrak{q}_\lambda)=p$ if $p<q$ and $R_\pm(\mathfrak{q})\in \{p,p-1\}$ when $p=q$  and that 
$\Phi(\mathfrak{u}_\lambda\cap\mathfrak{p})$ equals $\Phi(\mathfrak{u}\cap\mathfrak{p})$ where 
$\mathfrak{q}\in \{\mathfrak{q}_\omega,\mathfrak{q}_\mu,\mathfrak{q}_\kappa, \mathfrak{q}_\xi,\mathfrak{q}_\eta,  \mathfrak{q}_\nu\}$ considered in the discussion prior to the statement 
of the theorem.   This is decisive for the proof.

Since $\#\Phi(\mathfrak{u}_\lambda\cap\mathfrak{p})\le 2p$, we have $\#\Phi(\mathfrak{l}_\lambda\cap \mathfrak{p}_+)
\ge\#\Phi^+_\textrm{n}-2p=p(q-2)$.  Therefore there exists an $s\le p$ so that there are at least $q-2$ numbers among 
$a_j, p+1\le j\le p+q$, such that $a_s=a_j$.  
Suppose that $m$ is the cardinality of the set $C=\{ j\in [p+1,p+q]\mid a_j=a_s\}$.  Then $m\in \{q-2,q-1,q\}$. 
We break up the rest of the proof into three cases depending on the value of $m$.

Case (1): $m=q-2$.   We claim that $a_s=a_i$ for all $1\le i\le p$.  Otherwise there exists an $i$ such that     
$a_i\ne a_s$ for some $ i\le p.$ 
Any such $a_i$ equals {\it at most two} of the numbers $a_j, j\ge p+1$.  If there are $t\ge 1$ such numbers, we have $2t+(p-t)(q-2)\ge \#\Phi(\mathfrak{l}_\lambda\cap \mathfrak{p})\ge p(q-2).$  This implies that $2\ge (q-2)$ contrary to our assumption 
that $q\ge 5$.  
Therefore $a_s=a_i$ for all $i\le p$.   
Since $R_+(\mathfrak{q}_\lambda)>0$, 
$a_1-a_{p+q}=(\epsilon_1-\epsilon_{p+q},\lambda)>0$.  Similarly $R_-(\mathfrak{q}_\lambda)>0$ implies that $a_p-a_{p+1}<0$. 
Since $a_1=a_p$, it follows that $a_{p+1}>a_{p+q}$ and $p+1,p+q\notin C$. 
Now $\lambda=a_1(\sum_{1\le j\le p+q} \epsilon_j)+(a_{p+1}-a_1)\epsilon_{p+1} +(a_{p+q}-a_1)\epsilon_{p+q}=a_1\epsilon_0+
(a_{p+1}-a_1)\epsilon_{p+1}+(a_{p+q}-a_1)\epsilon_{p+q}$.  (Here $\epsilon_0=\sum_{1\le i\le p+q}\epsilon_i$.)  
Using $a_{p+1}>a_1>a_{p+q}$, a direct computation shows that $\Phi(\mathfrak{u}_\lambda
\cap \mathfrak{p}_\pm)=\Phi(\mathfrak{u}_\omega\cap \mathfrak{p}_\pm)$ where 
$\omega=\epsilon_{p+1}-\epsilon_{p+q}$.   

Case (2): $m=q-1$.  In this case $p+1\in C$ or $p+q\in C$. 
As in case (1) above, we have $(\epsilon_1-\epsilon_{p+q},\lambda)>0, (\epsilon_p-\epsilon_{p+1},\lambda)<0$, 
that is $a_1>a_{p+q}, a_p<a_{p+1}$.  
If $p+1\in C$, then $a_p-a_j<0$ for $p+1\le j\le p+q-1$ and so $-\epsilon_p+\epsilon_j\in \Phi(\mathfrak{u}_\lambda\cap \mathfrak{p})$ which implies that $q-1\le R_-(\mathfrak{q}_\lambda)\le p,$ that is, $q-1\le p$.  
We arrive at the same conclusion if $p+q\in C$ using $R_+(\mathfrak{q}_\lambda)\le p$.  

Let $A=\{i\le p\mid a_i\ne a_s\}, B=\{i\le p\mid a_i=a_s\}$.  
For each $i\in A$, $(\epsilon_i-\epsilon_j)$ or $(\epsilon_j-\epsilon_i) $ belongs to $\Phi(\mathfrak{u}_\lambda\cap \mathfrak{p}_+)$ for $(q-1)$ distinct values 
of $j, p+1\le j\le p+q$. For each $i\in B$, either $\epsilon_i-\epsilon_{p+1}$ or $-\epsilon_i+\epsilon_{p+q}$ belongs to $\Phi(\mathfrak{u}_\lambda\cap \mathfrak{p})$.  
It follows that, setting $a:=\#A$,  we have $\#B=p-a$ and 
$2p \ge R(\mathfrak{q}_\lambda)\ge a(q-1)+(p-a)=p+aq-2a$.  Therefore $a(q-2)\le p$.  Using the observation that $p\le q, q\ge 5, $ the only possibilities are $a=0,1$. 
 
Suppose that $a=0$. Then all the $a_i, 1\le i\le p$, are equal and so each $a_i=a_s, i\le p,$ equals $a_j$ 
for every $j\in C$. Thus $\epsilon_i-\epsilon_j\in \Phi(\mathfrak{l}_\lambda\cap \mathfrak{p}) ~\forall i\le p, j\in C$ 
and $\epsilon_i-\epsilon_{p+1}\in \Phi(\mathfrak{u}_\lambda\cap\mathfrak{p}_+)$ for all $i\le p$. So $\#\Phi(\mathfrak{l}_\lambda\cap \mathfrak{p}_+)=p(q-1)$ and 
$R_+(\mathfrak{q})=p,$ whence $R_-(\mathfrak{q}_\lambda)=0$, 
a contradiction.  So we are left with the possibility that 
$a=1$, in which case $A=\{1\}$ or $\{p\}$, in view of the monotonicity of $a_1,\ldots,a_p$.  
As observed earlier $p+1\in C $ or $p+q\in C$. There are four possibilities to consider, one for each choice of $A$ and $C$:

(a)  If $A=\{1\}, C=[p+1,p+q-1]\cap \mathbb{N}$, then $a_p=a_{p+1}$ which implies that $R_-(\mathfrak{q}_\lambda)=0$, 
a contradiction.  

(b) If $A=\{p\}, C=[p+2,p+q]\cap \mathbb{N}$, then $a_1=a_{p+q}$ which implies $R_+(\mathfrak{q}_\lambda)=0$, a contradiction. 

(c) Let $A=\{1\}, C=[p+2, p+q]\cap \mathbb{N}$.   Then $\epsilon_1-\epsilon_j\in \Phi(\mathfrak{q}_\lambda\cap\mathfrak{p}_+)$ for $p+1<j\le p+q$, $\epsilon_{p+1}-\epsilon_i\in \Phi(\mathfrak{q}_\lambda\cap\mathfrak{p}_-)$ for $2\le i\le p$.  
Therefore $R_+(\mathfrak{q}_\lambda)\ge q-1$, $R_-(\mathfrak{q}_\lambda)\ge p-1$.  If $a_1>a_{p+1}$, then $R_+(\mathfrak{q}_\lambda)=q$, $R_-(\mathfrak{q}_\lambda)=p-1$, which is impossible as $p\le q$ and $R_+(\mathfrak{q}_\lambda)=R_-(\mathfrak{q}_\lambda)$.  So $a_1\le a_{p+1}$.  If equality holds, then $R_+(\mathfrak{q}_\lambda)=q-1, R_-(\mathfrak{q}_\lambda)=p-1$.  This forces $p=q$ and we have $\lambda=a_p\epsilon_0 +(a_1-a_p)\epsilon_1+(a_1-a_p)\epsilon_{p+1}$. 
In this case $\Phi(\mathfrak{u}_\lambda\cap\mathfrak{p}_\pm)=\Phi(\mathfrak{u}_\xi\cap\mathfrak{p}_\pm)$ where $\xi=p(\epsilon_1+\epsilon_{p+1})
-\epsilon_0$.   On the other hand, if $a_1<a_{p+1}$, then $\epsilon_{p+1}-\epsilon_1\in \Phi(\mathfrak{q}_\lambda\cap 
\mathfrak{p}_-)$ and $R_-(\mathfrak{q}_\lambda)=p$.  Furthermore we must have $R_+(\mathfrak{q}_\lambda)=q-1$ whence 
$p=q-1$.  In this case $\lambda=(a_1-a_p)\epsilon_1+(a_{p+1}-a_p)\epsilon_{p+1}+a_p\epsilon_0$ and $\Phi(\mathfrak{u}_\lambda\cap\mathfrak{p}_\pm)=\Phi(\mathfrak{u}_\kappa\cap\mathfrak{p}_\pm)$ 
 where $\kappa=p\epsilon_1+q\epsilon_{p+1}-\epsilon_0$. 
 
(d)  Let $A=\{p\}, C=[p+1, p+q-1]\cap \mathbb{N}$.  This is similar to (c) above and we obtain, one of the two possibilities: Either 
$p=q-1, R_\pm(\mathfrak{q}_\lambda)=p$, in which case $\lambda=a_1\epsilon_0+(a_p-a_1)\epsilon_p+(a_{p+q}-a_1)\epsilon_{p+q}$, 
$\Phi(\mathfrak{q}_\lambda\cap \mathfrak{p}_\pm)
=\Phi(\mathfrak{q}_\nu \cap\mathfrak{p}_\pm)$ or $p=q, R_\pm(\mathfrak{q}_\lambda)=p-1$, $\lambda=a_1\epsilon_0+
(a_p-a_1)(\epsilon_p+\epsilon_{p+q})$, and $\Phi(\mathfrak{u}_\lambda\cap\mathfrak{p}_\pm)
=\Phi(\mathfrak{u}_\eta\cap\mathfrak{p}_\pm)$ where $\eta=\epsilon_0-p(\epsilon_p+\epsilon_{p+q})$.

Case (3) $m=q$.   Thus $C=\{p+1,\ldots, p+q\}$.  As $\epsilon_1-\epsilon_{p+q}\in 
\Phi(\mathfrak{u}_\lambda\cap\mathfrak{p}_+)$ 
we see that $\epsilon_1-\epsilon_j\in \Phi(\mathfrak{u}_\lambda\cap \mathfrak{p}_+)$ for all $j\ge p+1$.  
Thus $p\ge R_+(\mathfrak{q}_\lambda)\ge q$ and so $R_+(\mathfrak{q}_\lambda)=p=q$.   Similarly $-\epsilon_p+\epsilon_{p+1}\in 
\Phi(\mathfrak{u}_\lambda\cap \mathfrak{p}_-)$ implies that $R_-(\mathfrak{q}_\lambda)=p$. 
If $1<i<p$, then $a_i$ must equal $a_{p+1}$, for,  
otherwise $\epsilon_i-\epsilon_{p+1}$ or $-\epsilon_i+\epsilon_{p+1}\in \Phi(\mathfrak{u}_\lambda\cap \mathfrak{p})$ resulting in 
$R(\mathfrak{q}_\lambda)>2q=2p$.   It follows that $\lambda=(a_1-a_{p+1})\epsilon_1-(a_{p+1}-a_{p})\epsilon_p+a_{p+1}\epsilon_0$. In this case $\Phi(\mathfrak{q}_\lambda\cap\mathfrak{p}_\pm)=\Phi(\mathfrak{q}_\mu\cap\mathfrak{p}_\pm)$ where $\mu 
=\epsilon_1-\epsilon_p$ considered previously.  This completes the proof.
\end{proof}

\begin{remark}{\em 
Suppose that $2\le p\le q\le 4$, 
Apart from the representations $\mathcal{A}_\mathfrak{q}$ given by the theorem, which are 
valid also for $2\le p\le q\le 4$, there are a few exceptional ones with $1\le R_+(\mathfrak{q})=R_-(\mathfrak{q})\le p$ 
which we list below:\\ 
(i) when $p=2, q=3$,  $\lambda= \epsilon_1-\epsilon_2+\epsilon_3-\epsilon_5$ is the only exceptional case 
and we have $R_\pm(\mathfrak{q}_\lambda)=2$. \\
(ii) when $p=2, q=4$, $\lambda=\epsilon_1-\epsilon_2+\epsilon_3+\epsilon_4-\epsilon_5-\epsilon_6$ 
yields $R_\pm(\mathfrak{q}_\lambda)=2$. This is the only exceptional case.\\
(iii) when $p=3,q=3,4$, $\lambda=\epsilon_1-\epsilon_3+\epsilon_4-\epsilon_6$ yields $R_\pm(\mathfrak{q}_\lambda)
=3$.\\
(iv) when $p=4=q$ and $\lambda=\epsilon_1+\epsilon_2-\epsilon_3-\epsilon_4+\epsilon_5+\epsilon_6-\epsilon_7-\epsilon_8$ 
yields $R_\pm(\mathfrak{q})=4$.  There are no other exceptional cases.
}
\end{remark}

\begin{remark}\label{supqyq}
{\em 
Suppose that $p<q-1, q\ge 5$.  By the above proposition, if $R_\pm(\mathfrak{q}_\lambda)=p$ where $\lambda$ is $\mathfrak{k}$-dominant, then $\mathfrak{q}_\lambda=\mathfrak{q}_\omega$.  
It follows that $\Phi(\mathfrak{l}_\lambda)=\{\pm(\epsilon_i-\epsilon_j)\mid 1\le i< j\le p+q, i\ne p+1, j\ne p+q\}$.  The compact roots are $\pm(\epsilon_i-\epsilon_j), 1\le i<j\le p$ or $p+1<i<j<p+q.$   It is readily seen that $\{\epsilon_i-\epsilon_{i+1}\mid 1\le i<p+q, i\ne p, p+1, p+q-1\}\cup\{\epsilon_p-\epsilon_{p+2}\}$ is the set of simple roots for the positive system 
$\Phi^+(\mathfrak{l}_\lambda)=\Phi(\mathfrak{l}_\lambda)\cap \Phi^+$.  The only non-compact simple root 
is $\epsilon_p-\epsilon_{p+2}$.    Since $\mathfrak{t}\subset \mathfrak{l}_\lambda$, the rank of $\mathfrak{l}$ equals $p+q-1$.   In fact the centre of $\mathfrak{l}$ is spanned by the vectors $H_\omega$ and $H_\epsilon$ where 
$\epsilon:=\epsilon_0-(p+q)(\epsilon_{p+1}+\epsilon_{p+q})/2$.  
It follows that the 
real reductive Lie algebra $\mathfrak{l}_{0,\lambda}$ is isomorphic to $\mathfrak{su}(p,q-2)\oplus  \mathbb{R}iH_\omega\oplus\mathbb{R}iH_\epsilon$.   The connected Lie subgroup $L\subset \SU(p,q)$ corresponding 
to $\mathfrak{l}_{0,\lambda}$ is locally isomorphic to $\SU(p, q-2)\times \mathbb{S}^1\times \mathbb{S}^1.$  We note that 
the compact globally Hermitian symmetric space $Y_\mathfrak{q}$ dual to the symmetric space $L/(K\cap L)$ is the 
Grassmann manifold $\textrm{U}(p+q-2)/(\textrm{U}(p)\times \textrm{U}(q-2))\cong G_p(\mathbb{C}^{p+q-2})$. 
 }  
\end{remark}

\subsection*{Type BDI}
Let $G=\SO_0(2,p)$,  $K=\SO(2)\times \SO(p), p\ge 3$.  Set $n:=\lfloor p/2\rfloor +1=\textrm{rank}(\mathfrak{g})$.   
We have $\Phi^+=\{\epsilon_i\pm\epsilon_j\mid 1\le i<j\le n\} $ if $p$ is even, $\Phi^+=\{\epsilon_i\pm \epsilon_j\mid 1\le i<j\le n\}
\cup\{\epsilon_j\mid 1\le j\le n\}$, if $p$ is odd;  $\Delta_\mathfrak{g}=\{\psi_{j}:=\epsilon_j-\epsilon_{j+1}\mid 1\le j<n\}\cup 
\{\psi_n\}$ where $\psi_n=\epsilon_n$ if $p$ is odd and $\psi_n=\epsilon_{n-1}+\epsilon_n$ if $p$ is even.   The simple 
non-compact 
root is $\psi_1=\epsilon_{1}-\epsilon_2$ for any parity of $p$.  We have $\Phi^+_\n:=\{\epsilon_1\pm\epsilon_j\mid 1< j\le n\}$ if $p$ is even and $\Phi^+_\n:=\{\epsilon_1\pm \epsilon_j\mid 1<j\le n\}\cup \{\epsilon_1\}$ if $p$ is odd.    In this case 
we shall classify {\it all} $\theta$-stable parabolic subalgebras of $\mathfrak{g}_0$ having Hodge type of the form $(r,r)$ although our main concern is to show that $N(1)=1$.

It is readily 
verified that when $\lambda_r=\epsilon_2+\cdots+\epsilon_r, 2\le r\le n,$ we have $(\epsilon_i\pm\epsilon_j,\lambda)\ge 0$ 
for $2\le i<j\le n$, $(\epsilon_k,\lambda_r)\ge 0$ for $2\le k\le n$. Thus $(\gamma, \lambda_r)\ge 0$ for all 
$\gamma\in \Phi^+_\mathfrak{k}$ for any parity of $p$.    
$\Phi(\mathfrak{u}_{\lambda_r}\cap 
\mathfrak{p}_\pm)=\{\pm\epsilon_1+\epsilon_j\mid 2\le j\le r\}$ and so 
$R_\pm(\mathfrak{q}_{\lambda_r})=r-1$.   

When $p$ is even, $\mu_n=\epsilon_2+\cdots+\epsilon_{n-1}-\epsilon_n$ is also $\mathfrak{k}$-dominant
and we have $R_\pm(\mathfrak{q}_{\mu_n})=n$.  However $\Phi(\mathfrak{u}_{\mu_n}\cap \mathfrak{p})=\{\pm \epsilon_1+\epsilon_j\mid j<n\}
\cup \{\pm\epsilon_1-\epsilon_n\}\ne \Phi(\mathfrak{u}_{\lambda_n}\cap \mathfrak{p})$.   Moreover $2\lambda_n=|\Phi(\mathfrak{u}_{\lambda_n}\cap \mathfrak{p})|\ne |\Phi(\mathfrak{u}_{\mu_n}\cap \mathfrak{p})|=2\mu_n$.  This 
shows that the two representations $\mathcal{A}_{\mathfrak{q}_\lambda}, \lambda=\lambda_n,\mu_n$,  
are {\it inequivalent} representations of $G$ whence $N(n)\ge 2$.     

More generally, suppose that $\lambda=\sum_{1\le j\le n} a_j\epsilon_j\ne 0$.  Then $(\gamma,\lambda)\ge 0$ 
for all $\gamma\in \Phi^+_\mathfrak{k}$ if and only if the following condition holds depending on the parity of $p$: 
(i) when $p$ is odd,  $a_2\ge \cdots\ge a_{n}\ge 0$ and, 
(ii) when $p$ is even, $a_2\ge \cdots\ge a_{n-1}\ge |a_n|$.  Assume that $\lambda \ne 0$ satisfies this condition. 

\begin{lemma}  We keep the above notation.
Suppose that $\lambda=\sum_{1\le i\le n}a_i\epsilon_i$ is $\mathfrak{k}$-dominant and that 
$R_+(\mathfrak{q}_\lambda)= R_-(\mathfrak{q}_\lambda)$.  Then 
(i) $a_1=0$ if $p$ is odd, and, (ii) $|a_1|< |a_n|$ if $p$ is even.  
\end{lemma}
\begin{proof}
We will assume that $a_1>0$ the case $a_1<0$ being analogous---one merely has to interchange $\epsilon_1$ and $-\epsilon_1$ throughout. 
Also we assume that $a_n\ge 0$ when $p$ is even.   
The case when $a_n<0$ is similar---one merely has to interchange 
$\epsilon_n$ with $-\epsilon_n$ throughout.  

We shall pair the positive non-compact root $\epsilon_1+\epsilon_j\in \Phi^+_\n$ with the negative non-compact
root $-\epsilon_1+\epsilon_j\in \Phi^-_\n$. 
and similarly the root $\epsilon_1-\epsilon_j$ with $-\epsilon_1-\epsilon_j$.   In addition, when $p$ is odd, the non-compact root $\epsilon_1$ is paired with $-\epsilon_1$.

When $p$ is odd, as $a_1>0$ we note that $\epsilon_1\in \Phi(\mathfrak{u}_\lambda\cap \mathfrak{p})$ but not $-\epsilon_1$.

Let $i_0\le n$ be the largest integer such that $a_{i_0}>a_1$.   Similarly, let $i_1$ be the smallest integer such that 
$a_1> a_{i_1}$. 
If there is no such integer we put $i_1=n+1$.   
Since $R_-(\mathfrak{q}_\lambda)>0$, we have $-\epsilon_1+\epsilon_2\in \Phi(\mathfrak{u}_\lambda\cap \mathfrak{p}_-)$ and so $a_1<a_2$.  Thus  $2\le i_0<i_1$.   

If $1< j\le i_0$, then both the roots $\pm\epsilon_1+\epsilon_j$ belong to 
$\Phi(\mathfrak{u}_\lambda\cap \mathfrak{p})$ and neither of the roots $\pm\epsilon_1
-\epsilon_j $ belong to $\Phi(\mathfrak{u}_\lambda\cap \mathfrak{p})$.  
 
Let $i_0<j<i_1$.   Then $a_j=a_1$ and none of the pairs of roots $\pm\epsilon_1-\epsilon_j$ belongs to $\Phi(\mathfrak{u}_\lambda\cap\mathfrak{p})$.   However $\epsilon_1+\epsilon_j\in \Phi(\mathfrak{u}_\lambda\cap \mathfrak{p})$ {\it but not its pair}  $-\epsilon_1+\epsilon_j$. 

If $i_1\le j\le  n$, then $\epsilon_1\pm \epsilon_j$ is in 
$\Phi(\mathfrak{u}_\lambda \cap \mathfrak{p})$ but neither of their paired  
roots $-\epsilon_1\pm \epsilon_j$ is in $\Phi(\mathfrak{u}_\lambda\cap\mathfrak{p})$.  In this case, 
for each such $j$ there are two non-compact positive roots which 
belong to $\Phi(\mathfrak{u}\cap\mathfrak{p})$ but {\it there are no matching negative non-compact roots} in $\Phi(\mathfrak{u}_\lambda\cap \mathfrak{p})$.   

The above observations, together with the equality $R_+(\mathfrak{q}_\lambda)=R_-(\mathfrak{q}_\lambda)$, imply that $a_1=0$ if $p$ is odd and that 
and $i_1=n+1$. Thus if $a_1>0$, we must have $p$ is even and $a_1<a_n$ and the lemma follows. 
\end{proof}

Finally, suppose that $R_+(\mathfrak{q}_\lambda)=R_-(\mathfrak{q}_\lambda)$.  
Now observe that when $p$ is even and $0<|a_1|<|a_n|$, $\mathfrak{q}_\lambda$ equals $\mathfrak{q}_{\mu}$ 
where $\mu=\sum_{2\le j\le n}a_j\epsilon_j$. Therefore, in view of the above lemma we may (and do) assume that $a_1=0$.    
We see that $\Phi(\mathfrak{q}_\lambda\cap \mathfrak{p}_\pm)=\Phi(\mathfrak{q}_{\mu}\cap\mathfrak{p}_\pm)$ where 
$\lambda=\lambda_r$ when $a_n\ge 0$ and 
$r\ge 1$ is the largest number such that $a_r>0$ and $\lambda=\mu_n$ (defined above) 
when $p$ is even and $a_n<0$.   We have proved 

\begin{proposition} Suppose that  $G=\SO_0(2,p)$ and $n=\lfloor p/2\rfloor+1$.  Then  
$N(r)=1$ for $1\le r<n$.  Also $N(n)=\left\{\begin{array}{cc} 1 & ~\textrm{if $p$ is odd}, \\
2 &\textrm{~if $p$ is even}.\\
\end{array}\right. $ \hfill $\Box$ 
\end{proposition}

\begin{remark}\label{bd1yq}  {\em When $\lambda=\epsilon_2$, we have $\Phi(\mathfrak{l}_\lambda)$ consists of 
the roots $\{\pm(\epsilon_i\pm \epsilon_j)\mid 1\le i<j\le n, i\ne 2\ne j\}$ if $p$ is even, and, when 
$p$ is odd, besides the above set of roots  we have also the roots $\epsilon_j, 1\le j\le n, j\ne 2$.   
Here $n=\lfloor p/2\rfloor +1$. 
Also $H_{\lambda}$ is in the centre of $\mathfrak{l}_\lambda.$ 
Thus, for any parity of $p$,  $\mathfrak{l}_\lambda
\cong \mathfrak{so}(p,\mathbb{C})\oplus \mathbb{C}H_\lambda$.   The set of simple roots 
of $\mathfrak{l}_\lambda$ for the positive system $\Phi(\mathfrak{l}_\lambda)\cap \Phi^+$ consists 
of $\epsilon_1-\epsilon_3, \epsilon_i-\epsilon_{i+1}, 3\le i\le n,$ when $p$ is even; we have one 
more simple root, namely $\epsilon_n$, when $p$ is odd.  The only non-compact simple root is $\epsilon_1-\epsilon_3$. 
It follows that, for 
any parity of $p$, $\mathfrak{l}_{\lambda,0}$ is isomorphic to $\mathfrak{so}(2,p-2)\oplus \mathbb{R}iH_\lambda$.  
It is readily seen that the connected Lie subgroup $L$ of $\SO_0(2,p)$ corresponding to $\mathfrak{l}_{\lambda,0}$ 
is (locally) isomorphic to $\SO(2,p-2)\times \mathbb{S}^1$.  
Hence the compact dual of the symmetric space $L/(K\cap L)$ is isomorphic to $\SO(p)/(\SO(2)\times \SO(p-2))$.  This 
homogeneous space may be identified with the 
non-singular complex quadric $Q_{p-2}$ defined by the vanishing of the equation $z_1^2+\cdots +z_p^2 $ 
in the complex projective $(p-1)$-space $\mathbb{C}P^{p-1}$.   
}
\end{remark}


\subsection*{Type CI}
Let $G=\Sp(n,\mathbb{R})$ so that $K\cong U(n)$.  We have $\Phi^+=\{\epsilon_i\pm \epsilon_j\mid 1\le i<j\le n\}\cup \{2\epsilon_j\mid 1\le j\le n\}, \Delta_\mathfrak{g}=\{\psi_j=\epsilon_j-\epsilon_{j+1}\mid 1\le j<n\}\cup\{\psi_n=2\epsilon_n\},$ with $\psi_n$ being the non-compact simple root.  Also $\Phi_\n^+=\{\epsilon_i+\epsilon_j\mid 1\le i\le j\le n\}$.  
Let $\lambda=\epsilon_1-\epsilon_n.$  Then $(\epsilon_i-\epsilon_j,\lambda)\ge 0$ for all $1\le i< j\le n$. That is, $\lambda$ is 
$\mathfrak{k}$-dominant.  Consider $\mathfrak{q}_\lambda$.  We have $\Phi(\mathfrak{u}_\lambda\cap \mathfrak{p}_+)=
\{\epsilon_1+\epsilon_j\mid 1\le j<n\}$, $\Phi(\mathfrak{u}_\lambda\cap \mathfrak{p}_-)=\{-\epsilon_j-\epsilon_n\mid 1< j\le n\}$. 
Thus $R_\pm(\mathfrak{q}_{\epsilon_1-\epsilon_n})=n-1$.   

More generally let $\lambda=\sum_{1\le i\le n} a_i\epsilon_i$ be $\mathfrak{k}$-dominant, that is, $(\gamma,\lambda)\ge 0$ for all $\gamma\in \Phi^+(\mathfrak{k})$.  Equivalently we have that $a_1\ge \cdots\ge a_n$.   The following 
observations will be used without explicit mention:
(a) if $\epsilon_i+\epsilon_j\in \Phi(\mathfrak{u}_\lambda\cap\mathfrak{p}_+)$, then $\epsilon_p+\epsilon_q\in 
\Phi(\mathfrak{u}_\lambda\cap\mathfrak{p}_+)$ for $1\le p\le i, 1\le q\le j$;   
(b) if $-(\epsilon_i+\epsilon_j)\in \Phi(\mathfrak{u}_\lambda\cap\mathfrak{p}_-)$, then $-\epsilon_p-\epsilon_q\in \Phi(\mathfrak{u}_\lambda\cap\mathfrak{p}_-)$ for $i\le p\le n, j\le q\le n$.  

Our aim is to establish the following result.

\begin{theorem}\label{ci}  Let $G=\Sp(n,\mathbb{R})$. (i) Suppose that $n\ne 4$.   There exists a unique (up to equivalence) 
unitary representation $\mathcal{A}_\f{q}$ of  $\emph{Sp}(n,\r)$ having Hodge type $(n-1,n-1)$.  Thus $N(n-1)=1$ and 
$N(r)=0$ for $1\le r\le n-2$.  Moreover if $R_\pm(\mathfrak{q}_\lambda)=n-1$ with $\lambda$ being $\mathfrak{k}$-dominant, then $\lambda=r(\epsilon_1-\epsilon_n)$ for some $r>0$.  (ii) When $n=4$, we have $N(3)=2$ and $N(2)=N(1)=0$.
\end{theorem}
\begin{proof}
When $n=3$ it is easily verified that $N(1)=0$ and that $N(2)=1$ corresponding to $\lambda=\epsilon_1-\epsilon_3$.  

Let $n=4$.  We have $R_\pm(\mathfrak{q}_\lambda)=R_\pm(\mathfrak{q}_\mu)=2$ when $\lambda=\epsilon_1-\epsilon_4$ and $\mu=\epsilon_1+\epsilon_2-\epsilon_3-\epsilon_4$.  It is readily 
checked that $|\Phi(\mathfrak{u}_\lambda\cap \mathfrak{p})|=4\lambda$, $|\Phi(\mathfrak{u}_\mu\cap\mathfrak{p})|=3\mu$ and 
so $\mathcal{A}_{\mathfrak{q}_\lambda}$ and $\mathcal{A}_{\mathfrak{q}_\mu}$ are inequivalent representations, resulting in  $N(3)\ge 2$.   It is easy to see that $N(3)\le 2$ and so equality must hold.   Again it is trivial to verify that $N(2)=N(1)=0$.  

Now assume that $n\ge 5$.
Let $\lambda=\sum c_i\epsilon_i\ne 0$ be $\mathfrak{k}$-dominant.  Then $c_1\ge \cdots\ge c_n$.
Suppose that $0<R_+(\f{q}_\lambda)=R_-(\f{q}_\lambda)\le n-1$.    We will 
show that $\lambda$ is a positive multiple of $\epsilon_1-\epsilon_n$.

If $(\lambda, \epsilon_1+\epsilon_n)>0,$  then $(\lambda, \epsilon_1+\epsilon_i)>0$ for $1\le i\le n$ and so $R_+(\mathfrak{q}_\lambda)\ge n$, a contradiction.  Similarly $(\lambda, \epsilon_1+\epsilon_n)<0$ contradicts the hypothesis that 
$R_-(\mathfrak{q}_\lambda)<n$.   So $c_1+c_n=(\lambda,\epsilon_1+\epsilon_n)=0$.

Let $a\le n$ be the largest positive integer such that $c_j=c_n$ for $n-a+1\le j\le n$ and let $b\le n$ be the largest positive integer so that $c_j=c_1$ for $1\le j\le b$.  Since $c_1>c_n$ we see that $1\le a,b\le n-1$; also $a+b\le n$.

We have 
$(\epsilon_i-\epsilon_j,\lambda)=0$ for all $i,j \le b$ and $i, j>n-a$. Since $c_{n-a}>c_{n-a+1},~c_b>c_{b+1}$,  we have, 
for $i\le j$ with $i>b,  j>n-a$, 
$(\epsilon_i+\epsilon_j,\lambda )=c_i+c_j<c_1+c_n=0$ 
and similarly for  $i\le j$ with $i\le b,  j\le n-a$,  we have  
$(\epsilon_i+\epsilon_j,\lambda)=c_i+c_j>c_1+c_{n}=0$.
Thus $R_+(\f{q})\ge b(n+1-a-b)+{b\choose2}$ and $R_-(\f{q})\ge a(n+1-a-b)+{a\choose2}$. There are 
three cases to consider.\\

{\it Case} (i):  Let $a=1$. 
Using the above estimates we get $R_+(\f{q}_\lambda)\ge n-1$. If $b>1$, we have $R_+(\f{q}_\lambda)\ge 2n-3$ implying that 
$n-1\ge 2n-3$ or $n\le 2$. This is contrary to our hypothesis. Hence we must have $b=1$ and so 
$R_-(\mathfrak{q})=n-1$.  This implies that $(2\epsilon_i,\lambda)=0$ for $2\le i\le n-1$.  As $c_1+c_n=0$, 
we obtain that $\lambda=r(\epsilon_1-\epsilon_n)$ for some $r>0$.  

{\it Case} (ii):  Let $b=1$.  This is analogous to Case (i) leading to the same conclusion.

{\it Case} (iii):  Suppose that $a, b\ge 2$.   
We have $(\lambda, \epsilon_i+\epsilon_j)\ge (\lambda,\epsilon_1+\epsilon_{n-a})>0$ for 
$i\le 2, 1\le j\le n-a$.  Hence we obtain that $R_+(\mathfrak{q}_\lambda)\ge 2n-2a-1$.    
As $R_+(\mathfrak{q}_\lambda)\le n-1$, we conclude that $2(n-a)-1\le n-1$ or $a\ge n/2$.  Similarly $b\ge n/2$.
Since $a+b\le n$, and since $a, b$ are both integers we conclude that $a=b=n/2$.   In particular $n$ cannot be odd.  
Write $n=2m$,  so that $a=b=m$.  Using the estimate $n-1\ge R_\pm(\mathfrak{q})\ge b(n+1-a-b)+{b\choose 2}$ we obtain that $2m-1\ge m+{m\choose 2}$ implying $m\le 2$. This is a contradiction as $n\ge 5$.  

This completes the proof.
\end{proof}

\begin{remark}\label{ciyq}
{\em 
(i) Let $R_\pm(\mathfrak{q}_\lambda)=n-1$.  By the above result, for $n\ne 4$, $\lambda=r(\epsilon_1-\epsilon_n)$ for some $r>0$ and so $\Phi(\mathfrak{l}_\lambda)$ consists of the roots $\{\pm 2\epsilon_i\mid 2\le i\le n-1\}\cup \{\pm(\epsilon_i\pm\epsilon_j)
\mid 2\le i<j\le n-1\}\cup \{\pm (\epsilon_1+\epsilon_n)\}$.  Hence $\mathfrak{l}\cong \mathfrak{sp}(n-2, \mathbb{C})\oplus
\mathfrak{sl}(2,\mathbb{C})\oplus \mathbb{C}H_\lambda$.    The roots of the summand $\mathfrak{sl}(2,\mathbb{C})$ 
are the non-compact roots $\pm(\epsilon_1+\epsilon_n)$ and the 
summand $\mathfrak{sp}(n-2,\mathbb{C})$ is generated by the set of simple roots $\epsilon_i-\epsilon_{i+1}, 2\le i\le n-3,~2\epsilon_{n-1}$.   It follows that $\mathfrak{l}_{\lambda,0}$ is isomorphic to 
$\mathfrak{sp}(n-2,\mathbb{R})\oplus \mathfrak{su}(1,1)\oplus \mathbb{R}iH_\lambda$.     
We conclude that 
the connected Lie group $L\subset G$ with Lie algebra $\mathfrak{l}_{\lambda,0}$ is locally isomorphic to $\Sp(n-2,\mathbb{R})\times \SU(1,1)\times \mathbb{S}^1$.  
We further note that $\mathfrak{q}_\lambda$ does not satisfy the hypothesis of \cite[Prop. 6.1]{li}.   

The compact dual of the symmetric space $L/(K\cap L)$ is the space $Y_\mathfrak{q}\cong  \textrm{Sp}(n-2)/\textrm{U}(n-2) \times \mathbb{S}^2$. 
}
\end{remark}

\subsection*{Type DIII}
Let $G=\SO^*(2n)$; thus $K\cong \textrm{U}(n)$.  We have $\Phi=\{\pm(\epsilon_i\pm\epsilon_j)\mid 1\le i<j\le n\}$ where 
the set of simple roots is $\{\psi_i:=\epsilon_i-\epsilon_{i+1}\mid 1\le i<n\}\cup\{\psi_n:=\epsilon_{n-1}+\epsilon_n\}$; the 
non-compact simple root being $\psi_n$. 
The set of non-compact positive roots is $\Phi^+_\n=\{\epsilon_i+\epsilon_j\mid 1\le i<j\le n\}$.  
A weight $\lambda=\sum_j a_j\epsilon_j\in i\mathfrak{t}_0^*$  is $\mathfrak{k}$-dominant
if and only if $a_1\ge a_2\ge \cdots \ge a_n$.   

When $\lambda=\epsilon_1-\epsilon_n$, we have 
$\Phi(\mathfrak{u}_\lambda\cap\mathfrak{p}_+)=\{\epsilon_1+\epsilon_j\mid 1< j\le n-1\}$ and 
$\Phi(\mathfrak{u}_\lambda\cap\mathfrak{p}_-)=\{-\epsilon_j-\epsilon_n\mid 1<j\le n-1\}$.  
Therefore $R_\pm(\mathfrak{q}_{\lambda})=n-2$.  

More generally when $\lambda=\sum_{1\le i\le n}a_i\epsilon_i$ is $\mathfrak{k}$-dominant, we 
have the following properties, as in Type CI:  
(a) if $\epsilon_i+\epsilon_j\in \Phi(\mathfrak{u}_\lambda\cap\mathfrak{p}_+)$, then $\epsilon_p+\epsilon_q\in 
\Phi(\mathfrak{u}_\lambda\cap\mathfrak{p}_+)$ for $1\le p\le i, 1\le q\le j$;   
(b) if $-(\epsilon_i+\epsilon_j)\in \Phi(\mathfrak{u}_\lambda\cap\mathfrak{p}_-)$, then $-\epsilon_p-\epsilon_q\in 
\Phi(\mathfrak{u}_\lambda\cap\mathfrak{p}_-)$ for $i\le p\le n, j\le q\le n$.  

We have the following theorem, which was 
first established in \cite{mondal} and can also be found in \cite[Prop. 3.7]{mondal-sankaran}.

\begin{theorem}\label{d3}
Let $G=\SO^*(2n)$ and let $n\ge 9$.  Then there exists a unique irreducible unitary representation $\mathcal{A}_\mathfrak{q}$ (up to unitary equivalence) having Hodge type $(n-2,n-2)$.  Moreover, if $R_\pm(\mathfrak{q}_\lambda)=n-2$, then 
$\lambda=r(\epsilon_1-\epsilon_n)$ for some $r>0$.  
Thus $N(n-2)=1$.  Also $N(r)=0$ for $1\le r\le n-3$ and for $r=n-1$. 
\end{theorem}
\begin{proof} 
Suppose that $\lambda=\sum_{1\le i\le n} c_i\epsilon_i$ is $\mathfrak{k}$-dominant so that $c_1\ge \cdots\ge c_n$.  Assume that $1\le R_+(\mathfrak{q}_\lambda)=R_-(\mathfrak{q}_\lambda)\le n-1$.  

We claim that $c_1+c_n=0$.  To obtain a contradiction, assume that $(\lambda,\epsilon_1+\epsilon_n)=c_1+c_n>0$. 
Then $\epsilon_1+\epsilon_j\in \Phi(\mathfrak{u}_\lambda\cap \mathfrak{p}_+)$ for $1<j\le n$.   As $R_+(\mathfrak{q}_\lambda)\le n-1$, equality must hold and no other positive root 
$\epsilon_i+\epsilon_j, 2\le i<j\le n$ can be in $\Phi(\mathfrak{u}_\lambda\cap \mathfrak{p}_+)$.  
Hence $c_i+c_j\le 0$ 
for $2\le i<j\le n$.  If $c_2+c_{n-1}<0,$ then $c_i+c_j<0~\forall 2\le i< j, n-1\le j\le n$ and so $R_-(\mathfrak{q}_\lambda)\ge 2n-5.$
This implies that $2n-5\le n-1$, i.e., $n\le 4$ contrary to our hypothesis.  So $c_2+c_{n-1}=0$.   
This implies that $c_2=-c_j, 3\le j\le n-1$.  If $c_2>0$, then $c_3<0$ and we see that $c_i+c_j\in\Phi( \mathfrak{q}_\lambda
\cap\mathfrak{p}_-)~\forall 3\le i<j\le n$.  This yields $n-1\ge R_-(\mathfrak{q}_\lambda)\ge {n-2\choose 2}$, a contradiction as $n\ge 6$.   
Therefore we must have $c_i=0$ for $2\le i\le n-1$.  Now $c_1+c_{n}>0$ implies that 
$R_-(\mathfrak{q}_\lambda)=n-2<R_+
(\mathfrak{q}_\lambda)$, contrary to our hypothesis that $R_+(\mathfrak{q}_\lambda)=R_-(\mathfrak{q}_\lambda)$.  
Similarly we rule out $c_1+c_n<0$ by considering $\lambda'=\sum -c_{n+1-i}\epsilon_i$ which is also $\mathfrak{k}$-dominant.
Thus we are forced to conclude that $c_1+c_n=0$. 

As in proof of Theorem \ref{ci}, let $a=\#\{2\le i\le n\mid c_i=c_n\}, 
b=\#\{1\le i<n\mid c_i=c_1\}$.  We have $a,b\ge 1, a+b\le n$.   There are three cases to consider.

{\it Case} (i):  Let $a=1$.  Then $c_1+c_j>0~\forall 2\le j\le n-1$ and so $R_+(\mathfrak{q}_\lambda)\ge n-2.$ 
 If $b>1$, then (as in the proof of Theorem \ref{ci}) we obtain the lower bound $R_+(\mathfrak{q})\ge 2n-5$. 
As $R_+(\mathfrak{q})\le n-1$, this is impossible if $n\ge 5$.   So $b=1.$

Since $R_\pm (\mathfrak{q}_\lambda)\le  n-1$, we must have 
$c_2+c_4=0$ and   
also $c_3+c_4=0, c_3+c_5=0, c_4+c_5=0$ using $n\ge 7$.  Therefore 
$c_2=c_3=0=c_4=c_5$.  Hence 
we must have $R_+(\mathfrak{q}_\lambda)  \le n-2$.  By what has been shown already, $R_\pm(\mathfrak{q}_\lambda)=n-2$.
It follows that $c_j=0~\forall 2\le j\le n-1$ and so 
$\lambda=c_1(\epsilon_1-\epsilon_n)$ in this case.   

{\it Case} (ii):  Let $b=1$.  This is similar to the above case. 

It remains to consider the case $a,b>1$. 

{\it Case} (iii):  Let $a,b\ge 2$.  We will show that, under the hypotheses of the theorem, this leads to a contradiction.
This part of the proof is similar to that in the proof of Theorem \ref{ci}.  Then $n-1\ge R_+(\mathfrak{q}_\lambda)\ge 2n-2a-3$ and so $a\ge (n-2)/2$ and similarly $b\ge (n-2)/2$.   Also we have the estimate 
$R_+(\mathfrak{q}_\lambda)\ge a(n-a-b)+{a\choose 2}$.  Writing $m=\lfloor n/2 \rfloor$, if $n$ is odd then either $\{a,b\}=\{m\}$ or $\{m,m+1\}$. In either case $R_\pm(\mathfrak{q}_\lambda)\le n-1$, implies $2m\ge (m^2+m)/2$.  Hence $m\le 3$ or $n\le 7$ 
a contradiction as we assumed $n\ge 9$.  
Finally, let $n=2m\ge 10$ so that $\{a,b\}=\{m-1\}$ or $\{m-1,m\}$ or $\{m\}$ or $\{m-1,m+1\}$. In all cases we get the inequality
$2m-1\ge m(m-1)/2$.  This implies that $m\le 4$ which is a contradiction.  This completes the proof.
\end{proof}

\begin{remark}

When $4\le n\le 8$, there are more possibilities for the $\theta$-stable parabolic subalgebras $\mathfrak{q}$ with $R_+(\mathfrak{q}),  R_-(\mathfrak{q})\le n-1$.  
The following is the complete list of such `exceptional' $\theta$-stable parabolic subalgebras:
In all these cases $R_+(\mathfrak{q})=R_-(\mathfrak{q})$.\\    
$\underline{n=8}:$  The only exceptional $\mathfrak{q}$ corresponds to $\lambda=\epsilon_1+\epsilon_2+\epsilon_3+\epsilon_4-(\epsilon_5+\epsilon_6+\epsilon_7+\epsilon_8)$.  We have $R_\pm(\mathfrak{q})=6$.\\
$\underline{n=7}:$  The only exceptional $\mathfrak{q}$ corresponds to  $\lambda=\epsilon_1+\epsilon_2+\epsilon_3-\epsilon_5-\epsilon_6-\epsilon_7$, where $R_\pm(\mathfrak{q})=6$.\\ 
$\underline{n=6}:$  There three exceptions corresponding to $\lambda=\epsilon_1+\epsilon_2+\epsilon_3-\epsilon_4-\epsilon_5-\epsilon_6$ in which case $R_\pm(\mathfrak{q})=3$, $\lambda=\epsilon_1+\epsilon_2-\epsilon_5-\epsilon_6$ where $R_\pm(\mathfrak{q})=5$ and $\lambda=2(\epsilon_1-\epsilon_6)+(\epsilon_2+\epsilon_3)-(\epsilon_4+\epsilon_5)$ with $R_\pm(\mathfrak{q})=5$.\\
$\underline{n=5}:$  There are four exceptional cases 
corresponding to $\lambda=\epsilon_1+\epsilon_2-\epsilon_4-\epsilon_5$ with $R_\pm=3, \lambda=2(\epsilon_1-\epsilon_5)+\epsilon_2-\epsilon_4$ with $R_\pm(\mathfrak{q})=4$, $\lambda=3\epsilon_1+\epsilon_2-\epsilon_3-\epsilon_4-\epsilon_5$ with $R_\pm(\mathfrak{q})=4$ and $\lambda=2\epsilon_1+\epsilon_2+\epsilon_3-\epsilon_4-3\epsilon_5$ with $R_\pm(\mathfrak{q})=4$.  \\
$\underline{n=4}:$  There are three exceptional cases corresponds to $\lambda=\epsilon_1+\epsilon_2-\epsilon_3-\epsilon_4$ with $R_\pm(\mathfrak{q})=1$, $\lambda=2\epsilon_1-\epsilon_3-\epsilon_4$ with $R_\pm(\mathfrak{q})=3$ and $\lambda=\epsilon_1+\epsilon_2-2\epsilon_4$ with $R_\pm(\mathfrak{q})=3$.
\end{remark}

\begin{remark} \label{d3yq}
{\em Let $n\ge 9$.  
Let $R_\pm(\mathfrak{q}_\lambda)=n-2$ where $\lambda$ is $\mathfrak{k}$-dominant.   By the above theorem we have $\lambda=r(\epsilon_1-\epsilon_n)$ for some $r>0$. 
It follows that $\Phi(\mathfrak{l}_\lambda)=\{\pm(\epsilon_i\pm \epsilon_j)\mid 2\le i<j\le n-1\}\cup 
\{\pm(\epsilon_1+\epsilon_n)\}$.  The set of simple 
roots for the positive system $\Phi(\mathfrak{l}_\lambda)\cap \Phi^+$ is $\{\epsilon_i-\epsilon_{i+1}\mid 2\le i\le n-2\}
\cup \{\epsilon_{n-2}+\epsilon_{n-1}\}\cup \{\epsilon_1+\epsilon_n\}$.   We note that the root spaces corresponding 
to $\pm(\epsilon_1+\epsilon_n)$ spans a copy of $\mathfrak{sl}(2,\mathbb{C})$.  As the root 
$\epsilon_1+\epsilon_n$ is orthogonal to the remaining simple roots of $\mathfrak{l}_\lambda$, whose root spaces generate a subalgebra isomorphic to $\mathfrak{so}(2n-4,\mathbb{C})$ we have $[\mathfrak{l}_\lambda,\mathfrak{l}_\lambda] 
\cong \mathfrak{so}(2n-4,\mathbb{C})\oplus \mathfrak{sl}(2,\mathbb{C})$.    
The element 
$H_{\epsilon_1-\epsilon_n}$ spans the centre of $\mathfrak{l}_\lambda$.  Since the only simple non-compact 
roots of $\mathfrak{l}_\lambda$ are $\epsilon_{n-2}+\epsilon_{n-1}, \epsilon_1+\epsilon_n$ we 
see that $\mathfrak{l}_{\lambda,0}$ 
is isomorphic to $\mathfrak{so}^*(2n-4)\oplus \mathfrak{su}(1,1) \oplus\mathbb{R}iH_{\epsilon_1-\epsilon_n}$.  
Therefore the Lie group $L\subset \SO^*(2n)$ corresponding to $\mathfrak{l}_{\lambda,0}$ is locally isomorphic 
to $\SO^*(2n-4)\times \SU(1,1)\times \mathbb{S}^1$.  In particular, $\mathfrak{q}_\lambda$ does not satisfy the 
hypothesis of \cite[Prop. 6.1]{li}.

The compact dual $Y_\mathfrak{q}$ of the non-compact symmetric space $L/(K\cap L)$ equals 
$\SO(2n-4)/\textrm{U}(n-2)\times \mathbb{S}^2$.  
}
\end{remark}


\subsection*{Type EIII}  Let $G$ be a linear connected Lie group with Lie algebra 
the real form $\mathfrak{e}_{6, (-14)}$ of $\mathfrak{e}_6$ let $K$ be a maximal compact subgroup of $G$. 
Then $K$ is  
locally isomorphic to $\SO(2)\times \SO(10)$.    The vector space $i\mathfrak{t}_0^*$ may be realised as 
subspace of the Euclidean space 
$\mathbb{R}^8$ which is orthogonal to the space spanned by the vectors 
$\epsilon_6+\epsilon_8, \epsilon_7+\epsilon_8$.  (See \cite[Planche V]{bourbaki}.)
We let $\epsilon_0=\epsilon_8-\epsilon_7-\epsilon_6\in i\mathfrak{t}_0^*$. Then $\epsilon_i, 0\le i\le 5,$ is a basis for 
$i\mathfrak{t}_0^*$.    
 The simple roots of $\mathfrak{g}$ are $\psi_1=(1/2)(\epsilon_0+\epsilon_1-\epsilon_2-\epsilon_3-
\epsilon_4-\epsilon_5), 
\psi_2=\epsilon_1+\epsilon_2, \psi_3=\epsilon_2-\epsilon_1, \psi_4=\epsilon_3-\epsilon_2,\psi_5=\epsilon_4-\epsilon_3, 
\psi_6=\epsilon_5-\epsilon_4$.   The non-compact simple root is $\psi_1$.  The set of positive non-compact roots equals 
$\Phi^+_n=\{(1/2)(\epsilon_0+\sum_{1\le i\le 5}(-1)^{s_i}\epsilon_i)\mid s_i=0,1,~ \sum_{1\le i\le 5}s_i\equiv 0\mod 2\}. 
$  Also the set of positive compact roots equals 
$\Phi^+_\mathfrak{k}=\{(\epsilon_j\pm\epsilon_i)\mid 1\le i<j\le 5\}$. 

Denote the highest root $(1/2)(\epsilon_0+\sum_{1\le j\le 5}\epsilon_j)$  by $\alpha_0$.   
Let $\lambda=\sum_{0\le i\le 5}a_i\epsilon_i\in i\mathfrak{t}_0^*$.  
It is readily seen that $\lambda$ is $\mathfrak{k}$-dominant if and only if $-a_2\le a_1\le a_2\le a_3\le a_4\le a_5$.

\begin{figure}[h]
\centering
\begin{tikzpicture}
\draw (0,0) circle [radius=.1cm];
\draw (1,0) circle [radius=.1cm];
\draw (2,0) circle [radius=.1cm];
\draw (3,0) circle [radius=.1cm];
\draw[fill] (4,0) circle [radius=.1cm];
\draw (2,1) circle [radius=.1cm];
\draw (2,2) circle [radius=.1cm];

\node [below] at (0,0) {$\psi_6$};
\node [below] at (1,0) {$\psi_5$};
\node [below] at (2,0) {$\psi_4$};
\node [below] at (3,0) {$\psi_3$};
\node [below] at (4,0) {$\psi_1$};
\node [left] at (2,1) {$\psi_2$};
\node[left] at (2,2){$-\alpha_0$};

\draw (.1,0)--(.9,0);
\draw (1.1,0)--(1.9,0);
\draw (2.1,0)--(2.9,0);
\draw (3.1,0)--(3.9,0);
\draw (2,.1)--(2,.9);
\draw[dashed] (2,1.1)--(2,1.9);

\end{tikzpicture}
\caption{Extended Dynkin diagram of EIII with non-compact root $\psi_1$.}
\label{dynkine6}
\end{figure}
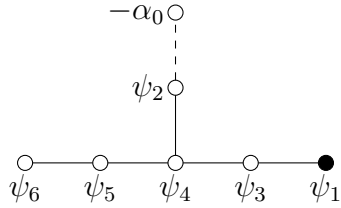

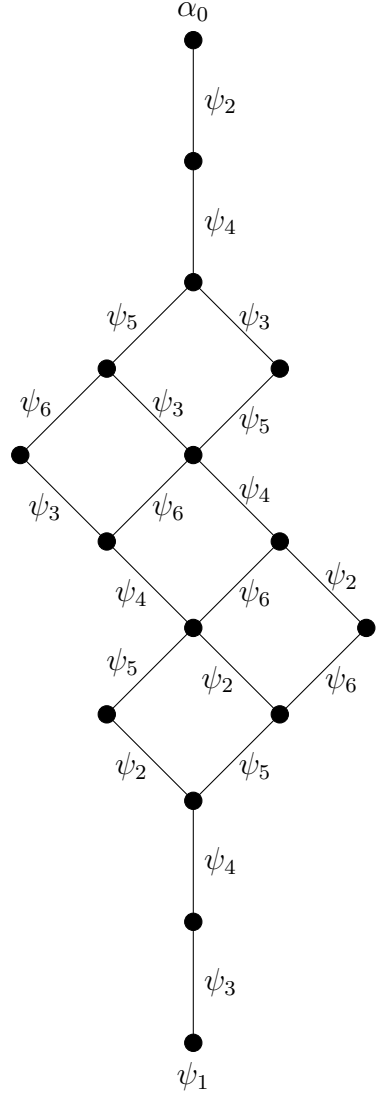
\begin{figure}[h]
\begin{tikzpicture}[scale=2.3]
\draw[fill] (0,.6) circle [radius=.05cm];
\draw[fill] (0,1.3) circle [radius=.05cm];
\draw[fill] (0,2) circle [radius=.05cm];
\draw[fill] (0,3) circle [radius=.05cm];
\draw[fill] (0,4) circle [radius=.05cm];
\draw[fill] (0,5) circle [radius=.05cm];
\draw[fill] (0,5.7) circle [radius=.05cm];
\draw[fill] (0,6.4) circle [radius=.05cm];
\draw[fill] (.5,2.5) circle [radius=.05cm];
\draw[fill] (.5,3.5) circle [radius=.05cm];
\draw[fill] (.5,4.5) circle [radius=.05cm];
\draw[fill] (-.5,2.5) circle [radius=.05cm];
\draw[fill] (-.5,3.5) circle [radius=.05cm];
\draw[fill] (-.5,4.5) circle [radius=.05cm];
\draw[fill] (-1,4) circle [radius=.05cm];
\draw[fill] (1,3) circle [radius=.05cm];
\draw (0,.6)--(0,1.3);
\draw (0,1.3)--(0,2);
\draw (0,2)--(-0.5,2.5);
\draw (0,2)--(.5,2.5);
\draw (-0.5,2.5)--(0,3);
\draw (.5,2.5)--(0,3);
\draw (.5,2.5)--(1,3);
\draw (0,3)--(.5,3.5);
\draw (1,3)--(.5,3.5);
\draw (0,3)--(-.5,3.5);
\draw (-.5,3.5)--(0,4);
\draw (.5,3.5)--(0,4);
\draw (-.5,3.5)--(-1,4);
\draw (-1,4)--(-.5,4.5);
\draw (0,4)--(-.5,4.5);
\draw (0,4)--(.5,4.5);
\draw (-.5,4.5)--(0,5);
\draw (.5,4.5)--(0,5);
\draw (0,5)--(0,5.7);
\draw (0,5.7)--(0,6.4);
\node [right] at (0,.95) {$\psi_3$};
\node [right] at (0,1.65) {$\psi_4$};
\node [right] at (.2,2.2) {$\psi_5$};
\node [right] at (.7,2.7) {$\psi_6$};
\node [right] at (.7,3.3) {$\psi_2$};
\node [right] at (.2,3.8) {$\psi_4$};
\node [right] at (.2,4.2) {$\psi_5$};
\node [right] at (.2,4.8) {$\psi_3$};
\node [right] at (0,5.35) {$\psi_4$};
\node [right] at (0,6.05) {$\psi_2$};
\node [left] at (.3,2.7) {$\psi_2$};
\node [right] at (.2,3.2) {$\psi_6$};
\node [right] at (-.3,3.7) {$\psi_6$};
\node [right] at (-.3,4.3) {$\psi_3$};
\node [left] at (-.2,2.2) {$\psi_2$};
\node [left] at (-.25,2.8) {$\psi_5$};
\node [left] at (-.2,3.2) {$\psi_4$};
\node [left] at (-.7,3.7) {$\psi_3$};
\node [left] at (-.75,4.3) {$\psi_6$};
\node [left] at (-.25,4.8) {$\psi_5$};
\node [below] at (0,.55) {$\psi_1$};
\node [above] at (0,6.45) {$\alpha_0$};
\end{tikzpicture}
\caption{Hasse diagram of the positive non-compact roots of EIII.}
\label{hassee6}
\end{figure}

\begin{proposition}  Let $G$ be a linear connected Lie group with Lie algebra $\mathfrak{e}_{6, (-14)}$.  Then  
$N(r)=0$ if $1\le r\le 3$, $N(4)=1$ and $N(5)=1$ and $N(6)=2$.
\end{proposition}

The proof will make repeated use of Remark \ref{rootsofu}(iii) without explicit reference.
It may be helpful to refer to Figure \ref{hassee6}. 

\begin{proof}
Suppose $1\le R_\pm(\mathfrak{q}_\lambda)\le 3$ with $\lambda$ being $\mathfrak{k}$-dominant. Then $(\lambda,\xi)=0$ for all roots $\xi\in\Phi^+_n$ such that $\psi_1+\psi_3+\psi_4<\xi<\alpha_0-\psi_2-\psi_4$. Then the roots $\psi_i, 2\le i\le 6$ are in the linear span of such roots. Thus $(\lambda,\psi_i)=0$ for all $2\le i\le 6$ which implies $\lambda=t\varpi_1$ for some $t\in\mathbb{R}$. This is impossible since $R_+(\mathfrak{q}_\lambda)=R_-(\mathfrak{q}_\lambda)\ge 1$.

Set $\beta:=\alpha-\psi_2-\psi_4-\psi_3-\psi_5$. When $R_\pm(\mathfrak{q}_\lambda)=4$, we must have 
$\Phi(\mathfrak{u}_\lambda\cap\mathfrak{p}_+)=\{\xi\in\Phi^+_\n\mid \xi\ge \beta+\psi_5\}$ or 
$\Phi(\mathfrak{u}_\lambda\cap\mathfrak{p}_+)=\{\xi\in\Phi^+_\n\mid \xi\ge\beta+\psi_3\}$. In the latter case $\beta,\beta-\psi_6\in\Phi(\mathfrak{l}_\lambda\cap\mathfrak{p}_+)$, which implies that $(\lambda,\psi_6)=0$. 
Therefore $(\lambda,\beta)=(\lambda,\beta-\psi_6)>0$ and so $R_\pm(\mathfrak{q}_\lambda)\ge 5$, a contradiction. 
So we must have $\Phi(\mathfrak{u}_\lambda\cap\mathfrak{p}_+)=\{\xi\in\Phi^+_\n\mid \xi\ge \beta+\psi_5\}$. Similarly 
$\Phi(\mathfrak{u}_\lambda\cap\mathfrak{p}_-)
=\{\xi\in\Phi^-_n\mid -\xi\le \psi_1+\psi_3+\psi_4+\psi_2\}$. It is easily verified that $\lambda=\varpi_5-\varpi_1$ yields $R_\pm(\mathfrak{q}_\lambda)=4$ and hence $N(4)=1$.

If $R_\pm(\mathfrak{q}_\lambda)=5$, then the only possibilities for  $\Phi(\mathfrak{u}_\lambda\cap\mathfrak{p}_+)$ are 
$\{\xi\mid \xi\ge \beta+\psi_3-\psi_6\}$
or $\{\xi\mid \xi>\beta\}$. In the latter case $\beta-\psi_6,\beta-\psi_6+\psi_3\in\Phi(\mathfrak{l}_\lambda\cap\mathfrak{p}_+)$, which implies that $(\lambda,\psi_3)=0$.   Therefore $(\lambda,\beta)=(\lambda,\beta+\psi_3)>0$ and so $R_\pm(\mathfrak{q}_\lambda)\ge 6$, a contradiction.  So we must have $\Phi(\mathfrak{u}_\lambda\cap\mathfrak{p}_+)=\{\xi\mid \xi\ge \beta+\psi_3-\psi_6\}$. Similarly $\Phi(\mathfrak{u}_\lambda\cap\mathfrak{p}_-)=\{\xi\in\Phi^-_\n\mid -\xi\le \psi_1+\psi_3+\psi_4+\psi_5+\psi_6\}$. It is easy to see that $\lambda=\varpi_2+\varpi_3-2\varpi_1$ yields $R_\pm(\mathfrak{q}_\lambda)=5$ and hence $N(5)=1$.

Suppose $R_\pm(\mathfrak{q}_\lambda)=6$ (where $\lambda$ is $\mathfrak{k}$-dominant), then there are only two possibilities for the set $\Phi(\mathfrak{u}_\lambda\cap \mathfrak{p}_+)$, namely,  
$\{\gamma\in \Phi^+_\n\mid \gamma\ge\beta\}=:A_1$ or $\{\gamma\in \Phi^+_\n\mid 
\gamma\ge \beta+\psi_3-\psi_6\}\cup \{\beta+\psi_5\}=:A_2$.   
In the former case we have $(\lambda,\psi_6)>0$ and in the latter we have $(\lambda,\psi_6)=0$. Arguing likewise with $\mathfrak{u}_\lambda\cap \mathfrak{p}_-$, we see that there are again two possibilities for 
$\Phi(\mathfrak{u}_\lambda\cap \mathfrak{p}_-)$ say $B_1, B_2$.   Exactly one of them, say  $B_1$, implies that $(\lambda,\psi_6)\ne 0$ and the other $B_2$, implies the vanishing of $(\lambda,\psi_6)$. It follows that only $A_i$ can be paired with $B_i$ for $i=1,2$, that is, $\Phi(\mathfrak{u}_\lambda\cap\mathfrak{p}_+)=A_i$ 
if and only if $\Phi(\mathfrak{u}_\lambda\cap \mathfrak{p}_-)=B_i, i=1,2$.  These two possibilities occur by choosing $\lambda=\varpi_4+\varpi_6-2\varpi_1$ and $\varpi_2+\varpi_3+\varpi_5-3\varpi_1$. Hence $N(6)=2$.  
\end{proof}

We tabulate in Table \ref{EIII} the weight $\lambda$ such that, writing $\mathfrak{q}=\mathfrak{q}_\lambda$, we have $1\le R_+(\mathfrak{q})=R_-(\mathfrak{q})\le 6$.  We also describe the 
corresponding compact Hermitian symmetric space $Y_\mathfrak{q}$ and its Euler characteristic.   We 
omit the detailed calculation that leads to the description of $Y_\mathfrak{q}$; it may be worked out 
easily using Proposition \ref{hermitian}.  

\begin{table}[h]
\centering
\begin{tabular}{|c|c|c|c|}
 \hline 
 $\lambda$ & $R_\pm(\mathfrak{q})$ & $Y_{\mathfrak{q}}$ &$\chi(Y_{\mathfrak{q}})$\\
 \hline
 $\varpi_5-\varpi_1$ & $4$&$G_2(\mathbb{C}^6)$&$15$\\
 \hline 
 $\varpi_2+\varpi_3-2\varpi_1$ & $5$&$\SO(8)/U(4)$ &$8$\\
 \hline 
 $\varpi_4+\varpi_6-2\varpi_1$ & $6$ & $\mathbb{S}^2\times \mathbb{S}^2$&$4$\\
 \hline
 $\varpi_2+\varpi_3+\varpi_5-3\varpi_1$ & $6$&$G_2(\mathbb{C}^4)$&$6$\\
\hline
\end{tabular}
\caption{The $\theta$-stable parabolic subalgebras $\mathfrak{q}=\mathfrak{q}_\lambda$ of type EIII for which 
$1\le R_\pm(\mathfrak{q}_\lambda)=r\le 6$, 
the symmetric spaces $Y_{\mathfrak{q}}$ and their Euler characteristics.}
\label{EIII}
\end{table}

\subsection*{Type EVII}
Let $G$ be a linear connected Lie group with $\mathfrak{g}_0$ isomorphic to $\mathfrak{e}_{7,(-25)}$.  
Then $K$ is locally isomorphic to the compact group $E_6\times \SO(2)$.   
As with the case of type EIII, it is convenient to set $\epsilon_0:=\epsilon_8-\epsilon_7-\epsilon_6\in \mathbb{R}^8$ and regard  
$i\mathfrak{t}^*_0\subset \mathbb{R}^8$ as the subspace orthogonal to $\epsilon_7+\epsilon_8$.   See \cite[Planche VI]{bourbaki}.
The simple roots are $\psi_1=(1/2)(\epsilon_0-\epsilon_5-\epsilon_4-\epsilon_3-\epsilon_2+\epsilon_1), 
\psi_2=\epsilon_1+\epsilon_2,\psi_3=\epsilon_2-\epsilon_1,\psi_4=\epsilon_3-\epsilon_2,\psi_5=\epsilon_4-\epsilon_3, \psi_6=\epsilon_5-\epsilon_4, \psi_7=\epsilon_6-\epsilon_5$.   The non-compact simple root is $\psi_7=\epsilon_6-\epsilon_5$. 

The set $\Phi^+_\n$ of non-compact positive roots equals 
$\{(1/2)(\epsilon_8-\epsilon_7+\epsilon_6+\sum_{1\le j\le 5} (-1)^{s_j}\epsilon_j)\mid \sum_{1\le j\le 6}(-1)^{s_j}\equiv 1 \mod 2 \}\cup\{\epsilon_6\pm\epsilon_i\mid 1\le i\le 5\}
\cup \{\epsilon_8-\epsilon_7\}$.    
The highest root is $\alpha_0:=2\psi_1+2\psi_2+3\psi_3+4\psi_4+3\psi_5+2\psi_6+\psi_7=\varpi_1
=\epsilon_8-\epsilon_7$.

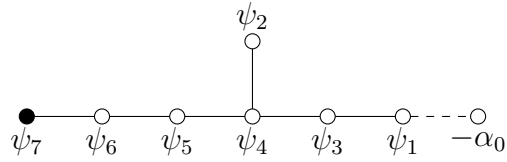
\begin{figure}
\centering
\begin{tikzpicture}
\draw[fill] (-1,0) circle [radius=.1cm];
\draw (0,0) circle [radius=.1cm];
\draw (1,0) circle [radius=.1cm];
\draw (2,0) circle [radius=.1cm];
\draw (3,0) circle [radius=.1cm];
\draw (4,0) circle [radius=.1cm];
\draw (5,0) circle [radius=.1cm];
\draw (2,1) circle [radius=.1cm];

\node [below] at (-1,0) {$\psi_7$};
\node [below] at (0,0) {$\psi_6$};
\node [below] at (1,0) {$\psi_5$};
\node [below] at (2,0) {$\psi_4$};
\node [below] at (3,0) {$\psi_3$};
\node [below] at (4,0) {$\psi_1$};
\node [above] at (2,1) {$\psi_2$};
\node[below] at (5,0) {$-\alpha_0$};

\draw (-.9,0)--(-.1,0);
\draw (.1,0)--(.9,0);
\draw (1.1,0)--(1.9,0);
\draw (2.1,0)--(2.9,0);
\draw (3.1,0)--(3.9,0);
\draw[dashed] (4.1,0)--(4.9,0);
\draw (2,.1)--(2,.9);

\end{tikzpicture}
\caption{Extended Dynkin diagram of EVII with non-compact root $\psi_7$.}
\label{dynkine7}
\end{figure}

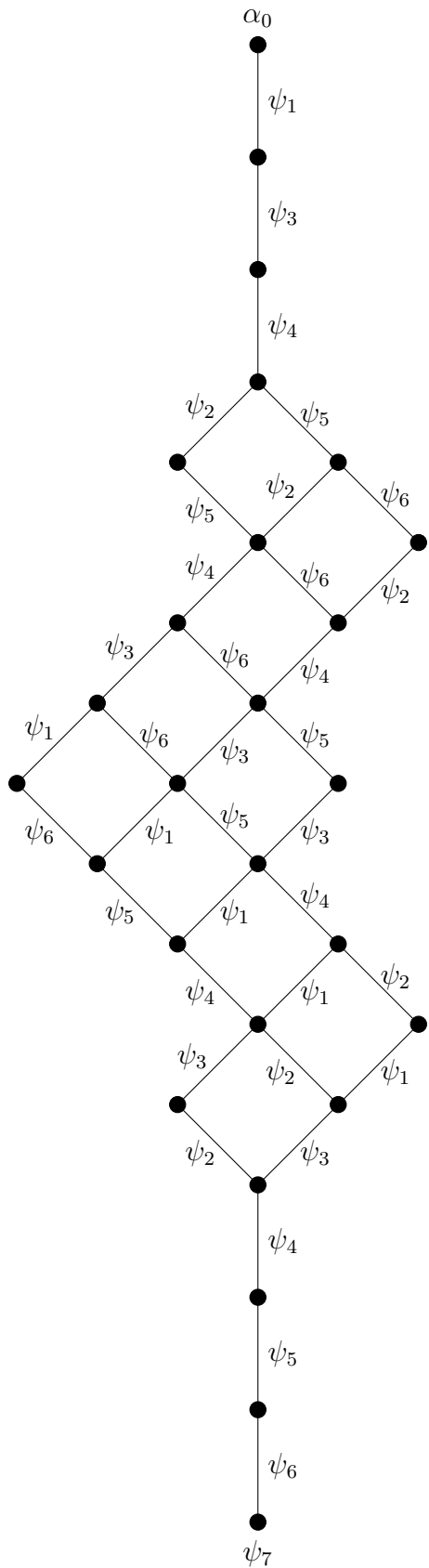
\begin{figure}
\centering
\begin{tikzpicture}[scale=2.3]
\draw[fill] (0,-.1) circle [radius=.05cm];
\draw[fill] (0,.6) circle [radius=.05cm];
\draw[fill] (0,1.3) circle [radius=.05cm];
\draw[fill] (0,2) circle [radius=.05cm];
\draw[fill] (0,3) circle [radius=.05cm];
\draw[fill] (0,4) circle [radius=.05cm];
\draw[fill] (0,5) circle [radius=.05cm];
\draw[fill] (0,6) circle [radius=.05cm];
\draw[fill] (0,7) circle [radius=.05cm];
\draw[fill] (0,7.7) circle [radius=.05cm];
\draw[fill] (0,8.4) circle [radius=.05cm];
\draw[fill] (0,9.1) circle [radius=.05cm];
\draw[fill] (.5,2.5) circle [radius=.05cm];
\draw[fill] (.5,3.5) circle [radius=.05cm];
\draw[fill] (.5,4.5) circle [radius=.05cm];
\draw[fill] (.5,5.5) circle [radius=.05cm];
\draw[fill] (.5,6.5) circle [radius=.05cm];
\draw[fill] (-.5,2.5) circle [radius=.05cm];
\draw[fill] (-.5,3.5) circle [radius=.05cm];
\draw[fill] (-.5,4.5) circle [radius=.05cm];
\draw[fill] (-.5,5.5) circle [radius=.05cm];
\draw[fill] (-.5,6.5) circle [radius=.05cm];
\draw[fill] (-1,4) circle [radius=.05cm];
\draw[fill] (-1,5) circle [radius=.05cm];
\draw[fill] (-1.5,4.5) circle [radius=.05cm];
\draw[fill] (1,3) circle [radius=.05cm];
\draw[fill] (1,6) circle [radius=.05cm];

\draw (0,-.1)--(0,.6);
\draw (0,.6)--(0,1.3);
\draw (0,1.3)--(0,2);
\draw (0,2)--(-0.5,2.5);
\draw (0,2)--(.5,2.5);
\draw (-.5,2.5)--(0,3);
\draw (.5,2.5)--(-1.5,4.5);
\draw (.5,2.5)--(1,3);
\draw (0,3)--(.5,3.5);
\draw (1,3)--(-1,5);
\draw (-1,4)--(1,6);
\draw (-1.5,4.5)--(.5,6.5);
\draw (-.5,3.5)--(.5,4.5);
\draw (.5,4.5)--(-.5,5.5);
\draw (.5,5.5)--(-.5,6.5);
\draw (1,6)--(0,7);
\draw (-.5,6.5)--(0,7);
\draw (0,7)--(0,9.1);

\node [right] at (0,.25) {$\psi_6$};
\node [right] at (0,.95) {$\psi_5$};
\node [right] at (0,1.65) {$\psi_4$};

\node [right] at (.2,2.2) {$\psi_3$};
\node [right] at (.7,2.7) {$\psi_1$};
\node [right] at (.7,3.3) {$\psi_2$};
\node [right] at (.2,3.8) {$\psi_4$};
\node [right] at (.2,4.2) {$\psi_3$};
\node [right] at (.2,4.8) {$\psi_5$};

\node [right] at (0,7.35) {$\psi_4$};
\node [right] at (0,8.05) {$\psi_3$};
\node [right] at (0,8.75) {$\psi_1$};

\node [left] at (.3,2.7) {$\psi_2$};
\node [right] at (.2,3.2) {$\psi_1$};
\node [right] at (-.3,3.7) {$\psi_1$};
\node [right] at (-.3,4.3) {$\psi_5$};
\node [right] at (-.8,4.8) {$\psi_6$};
\node [right] at (-.3,5.3) {$\psi_6$};
\node [right] at (.2,5.8) {$\psi_6$};
\node [right] at (.7,6.3) {$\psi_6$};
\node [right] at (.2,6.8) {$\psi_5$};
\node [left] at (-.2,6.2) {$\psi_5$};

\node [left] at (-.2,2.2) {$\psi_2$};
\node [left] at (-.25,2.8) {$\psi_3$};
\node [left] at (-.2,3.2) {$\psi_4$};
\node [left] at (-.7,3.7) {$\psi_5$};
\node [left] at (-1.2,4.2) {$\psi_6$};
\node [left] at (-.45,4.2) {$\psi_1$};
\node [right] at (-.3,4.7) {$\psi_3$};
\node [right] at (.2,5.2) {$\psi_4$};
\node [right] at (.7,5.7) {$\psi_2$};
\node [left] at (-1.2,4.85) {$\psi_1$};
\node [left] at (-.7,5.35) {$\psi_3$};
\node [left] at (-.2,5.85) {$\psi_4$};
\node [left] at (.3,6.35) {$\psi_2$};
\node [left] at (-.2,6.85) {$\psi_2$};

\node [below] at (0,-.15) {$\psi_7$};
\node [above] at (0,9.15) {$\alpha_0$};

\end{tikzpicture}
\caption{Hasse diagram of the positive non-compact roots of EVII.}
\label{hassee7}
\end{figure}

\begin{proposition} Let $G$ be a linear connected Lie group with Lie algebra $\mathfrak{e}_{7, (-25)}$.  Then 
 $N(r)=0$ for $1\le r\le 5, N(6)=1, N(7)=0, N(8)=0, N(9)=2, N(10)=1, N(11)=2$.
 \end{proposition}
 
 We will make repeated use of Remark \ref{rootsofu}(iii) without explicit mention.  
 
\begin{proof}
First suppose that $1\le R_\pm(\mathfrak{q}_\lambda)\le 7$ with $\lambda$ being $\mathfrak{k}$-dominant.  
Then $(\lambda,\xi)=0$ for all roots $\xi\in \Phi^+_n$ such that $\psi_7+\psi_6+\psi_5+\psi_4+\psi_3+\psi_2<\xi<\alpha_0-\psi_1-\psi_3-\psi_4-\psi_5-\psi_2$.  It is readily seen that roots 
$\psi_1, \psi_i, 3\le i\le 6,$ are in the linear span of such roots.  Set $\beta:=\alpha_0-\psi_1-\psi_3-\psi_4-\psi_5-\psi_2
\in \Phi^+_\n$.

 If $R_\pm(\mathfrak{q}_\lambda)\le 5$, then, by the same argument $(\lambda,\psi_2)$ also vanishes and so $\lambda=t\varpi_7$ for some 
$t\in \mathbb{R}$.  This is impossible when $R_+(\mathfrak{q}_\lambda)=R_-(\mathfrak{q}_\lambda)\ge 1$. 

When $R_\pm(\mathfrak{q}_\lambda)=6$, we must have $\Phi(\mathfrak{u}_\lambda\cap\mathfrak{p}_+)
=\{\xi\in \Phi^+_\n\mid \xi\ge \beta-\psi_6+\psi_2\}$ or 
$\Phi(\mathfrak{u}_\lambda\cap\mathfrak{p}_+)=\{\xi\mid \xi>\beta\}$.   In the latter case, $\beta+\psi_2-\psi_6, \beta-\psi_6\notin 
\Phi(\mathfrak{u}_\lambda\cap \mathfrak{p})$.  This implies that $(\lambda,\psi_2)=0$.  Therefore $\beta\in \Phi(\mathfrak{u}_\lambda\cap \mathfrak{p}_+)$ since $\beta+\psi_2\in \Phi(\mathfrak{u}_\lambda\cap \mathfrak{p}_+)$ and so $R_+(\mathfrak{q}_\lambda)\ge 7$, a contradiction.  So we must have $\Phi(\mathfrak{u}_\lambda\cap\mathfrak{p}_+)
=\{\xi\in \Phi^+_\n\mid \xi\ge \beta-\psi_6+\psi_2\}$.  Similarly, $\Phi(\mathfrak{u}_\lambda\cap\mathfrak{p}_-)=\{\xi\in
\Phi^-_\n\mid 
-\xi\le \psi_7+\psi_6+\psi_5+\psi_4+\psi_3+\psi_1\}$.  
It is easily verified that $\lambda=\varpi_2-\varpi_7$ yields $R_\pm(\mathfrak{q}_\lambda)=6$ and hence 
$N(6)=1$.  

Also when $R_\pm(\mathfrak{q}_\lambda)=7,$ the only possibilities for $\Phi(\mathfrak{u}_\lambda\cap \mathfrak{p}_+)
$ are $\{\xi\in \Phi^+_\n\mid \xi\ge \beta\}, \{\xi\in \Phi^+_\n\mid \xi\ge\beta-\psi_6+\psi_2\}\cup\{\beta+\psi_5\}$.  In former case, using the observation 
that $(\lambda,\psi_4)=0, (\lambda,\psi_6)=0$, we see that $\beta-\psi_4, \beta-\psi_6$ are in $\Phi(\mathfrak{u}_\lambda\cap \mathfrak{p}_+)$, a contradiction. In the latter case we have $(\lambda,\psi_5)=0$ which implies that $\beta\in\Phi(\mathfrak{u}_\lambda\cap\mathfrak{p}_+)$, again a contradiction. This proves that $N(7)=0$.  

Suppose that $R_\pm(\mathfrak{q}_\lambda)= 8$ where $\lambda$ is $\mathfrak{k}$-dominant.  
Then $\Phi(\mathfrak{u}_\lambda\cap \mathfrak{p}_+)$ equals 
$\{\xi\in \Phi^+_\n\mid \xi\ge \beta\}\cup\{\beta-\psi_4\}=:A$ or $\{\xi\in \Phi^+_\n\mid \xi\ge \beta\}
\cup\{\beta+\psi_2-\psi_6\}=:B$.    In the case $\Phi(\mathfrak{u}_\lambda\cap \mathfrak{p}_+)=A$, we have $(\lambda,\beta-\psi_6)=0$ as $(\lambda, \psi_6)=0$, a contradiction as $\beta-\psi_6\notin A$.  
 Now suppose that $\Phi(\mathfrak{q}_\lambda
\cap\mathfrak{p}_+)=B$.   Then $(\lambda, \psi_4)=0$, which implies that $\beta-\psi_4\in\Phi(\mathfrak{u}_\lambda\cap \mathfrak{p}_+),$ a contradiction.  Thus the claim that $N(8)=0$ is established.  

Next we turn to $N(9)$. Let $\lambda=\varpi_2+\varpi_4-3\varpi_7$. Then a straightforward verification shows that $\Phi(\mathfrak{u}_\lambda\cap\mathfrak{p}_+)=\{\xi\in \Phi^+_\n\mid \xi\ge \beta-\psi_6\}, \Phi(\mathfrak{u}_\lambda\cap\mathfrak{p}_-)=\{\xi\in 
\Phi^-_\n\mid -\xi\le \psi_7+\psi_6+\psi_5+\psi_4+\psi_3+\psi_2+\psi_1\}$ and so 
$R_\pm(\mathfrak{q}_{\lambda_i})=9$.   

Now  let  $\mu=\varpi_1+\varpi_6-2\varpi_7$.  Again by a direct verification $R_\pm(\mathfrak{q}_\mu)=9$.  
In fact we have $\Phi(\mathfrak{u}_\mu\cap\mathfrak{p}_+)=
\{\xi\in \Phi^+_\n\mid \xi\ge \alpha_0-\psi_1-2\psi_3-2\psi_4-\psi_2
-\psi_5\}$ and $\Phi(\mathfrak{u}_\mu\cap\mathfrak{p}_-)=\{\xi\in\Phi^-_\n\mid -\xi\le \psi_7+\psi_6+2\psi_5+2\psi_4
+\psi_3+\psi_2\}$.

The representations $\mathcal{A}_{\mathfrak{q}_\lambda}, \mathcal{A}_{\mathfrak{q}_\mu}$ are inequivalent 
since, $|\Phi(\mathfrak{u}_\lambda\cap \mathfrak{p})|\ne |\Phi(\mathfrak{u}_\mu\cap \mathfrak{p})|$ (as seen by 
comparing the coefficient of $\psi_6$ on both sides).  Therefore $N(9)\ge 2$.

We claim that $N(9)=2$. It suffices to show that there does not exist a $\mathfrak{k}$-dominant weight $\nu$ such 
that $\Phi(\mathfrak{u}_\nu\cap \mathfrak{p}_+)=\{\xi\in \Phi^+_\n\mid \xi \ge \beta -\psi_4\}
\cup \{\beta+\psi_2-\psi_6\}$ or 
$\Phi(\mathfrak{u}_\nu\cap \mathfrak{p}_-)=\{\xi\in \Phi^-_\n\mid   -\xi\le \psi_7+\psi_6+\psi_5+2\psi_4+\psi_3+\psi_2\}
\cup \{-(\psi_7+\psi_6+\psi_5+\psi_4+\psi_3+\psi_1)\}$.

Indeed, suppose that $\nu$ is such that $\Phi(\mathfrak{u}_\nu\cap \mathfrak{p}_+)=\{\xi\in \Phi^+_\textrm{n}\mid 
\xi \ge \beta -\psi_4\}\cup \{\beta+\psi_2-\psi_6\}$ 
Then $\beta-\psi_6-\psi_4, \beta-\psi_6-\psi_4-\psi_3
\in \Phi(\mathfrak{l}_\nu\cap \mathfrak{p}_+)$.  Hence $(\nu, \psi_3)=0$.  
 It follows that $(\nu,\beta-\psi_4-\psi_3)=(\nu,\beta-\psi_4)>0$.  Hence $\beta-\psi_4-\psi_3\in \Phi(\mathfrak{u}_\nu\cap\mathfrak{p}_+)$,  contrary to our hypothesis.  The possibility for $\Phi(\mathfrak{u}_\lambda\cap\mathfrak{p}_-)$ 
 is also likewise eliminated. 
 Thus $N(9)=2$.

Now suppose that $\lambda$ is a $\mathfrak{k}$-dominant weight such that $R_\pm(\mathfrak{q}_\lambda)=10$.   
Cardinality consideration shows that $\beta-\psi_4\in \Phi(\mathfrak{u}_\lambda\cap \mathfrak{p}_+)$. 
We claim that $\beta-\psi_6\notin 
\Phi(\mathfrak{u}_\lambda\cap \mathfrak{p}_+)$.  For, otherwise, 
arguing as before we see that $(\lambda,\psi_6)=0$, which implies that  
$(\lambda, \beta-\psi_4-\psi_6)=(\lambda,\beta-\psi_4)>0$ and hence $\beta-\psi_4-\psi_6\in 
\Phi(\mathfrak{u}_\lambda\cap \mathfrak{p}_+)$.   A simple cardinality argument then shows that $R(\mathfrak{q}_\lambda)\ge 11$, a contradiction.   So we must have $\beta-\psi_6\notin \Phi(\mathfrak{u}_\lambda\cap \mathfrak{p}_+)$ and $(\lambda,\psi_6)>0$.   This implies that $\beta-\psi_4-\psi_3
\in \Phi(\mathfrak{u}_\lambda\cap \mathfrak{p}_+)$.  Thus there are two possibilities: either 
$\Phi(\mathfrak{u}_\lambda\cap
\mathfrak{p}_+)$ equals $ \{\xi\in \Phi^+_\n\mid \xi\ge \beta-\psi_4-\psi_3-\psi_1 \}$ or $\{\xi\in \Phi^+_\n\mid \xi
\ge \beta-\psi_4-\psi_3 \}\cup \{\beta+\psi_2-\psi_6\}$.   Similarly there are two possibilities for 
$\Phi(\mathfrak{u}_\lambda\cap \mathfrak{p}_-)$.  Of the four possible combinations, 
 three are eliminated as in the determination of $N(9)$.  It turns out that when 
$\lambda=\varpi_1+\varpi_2+\varpi_6-3\varpi_7$, we have $R_\pm(\mathfrak{q}_\lambda)=10$.  Hence 
$N(10)=1$ as asserted.

 Next we show that $N(11)=2$.   Let $\lambda$ be a $\mathfrak{k}$-dominant weight with $R_\pm(\mathfrak{q}_\lambda)
 =11$.  
As in the above case, we must have $\beta-\psi_4-\psi_3-\psi_1\notin \Phi(\mathfrak{u}_\lambda\cap \mathfrak{p}_+)$. 
Hence $\beta-\psi_6, \beta-\psi_4\in\Phi(\mathfrak{u}_\lambda\cap \mathfrak{p}_+)$.  There are 
two possibilities: either $\beta-\psi_4-\psi_3\in \Phi(\mathfrak{u}_\lambda\cap \mathfrak{p}_+)$ or $\beta-\psi_4-\psi_6
\in \Phi(\mathfrak{u}_\lambda\cap \mathfrak{p}_+)$.  Accordingly, $\Phi(\mathfrak{u}_\lambda\cap \mathfrak{p}_+)
=\{\xi\in \Phi^+_\n\mid \xi\ge \beta-\psi_4-\psi_3\}\cup \{\xi\in \Phi^+_\n\mid \xi\ge \beta-\psi_6\}$, in which case $(\lambda,\psi_6)\ne 0$, or 
$\Phi(\mathfrak{u}_\lambda\cap \mathfrak{p}_+)=\{\xi
\mid \xi\ge \beta-\psi_6-\psi_4 \}$ in which case $(\lambda,\psi_6)=0$.  Similarly $\Phi(\mathfrak{u}_\lambda\cap \mathfrak{p}_-)=\{\xi\in \Phi^-_\n\mid -\xi\le 
\psi_7+\psi_6+2\psi_5+2\psi_4+\psi_2+\psi_3\}\cup \{\xi\in\Phi^-_\n\mid -\xi\le\psi_7+\psi_6+\psi_5+\psi_4+\psi_3+\psi_2+\psi_1\}$ in which case $(\lambda, \psi_6)>0$, or $\{\xi\in\Phi^-_\n\mid -\xi\le\psi_7+\psi_6+\psi_5+2\psi_4+\psi_3+\psi_2+\psi_1\}$ in which case $(\lambda,\psi_6)=0$.  Thus, of the four combinations for the pair 
$\Phi(\mathfrak{u}_\lambda\cap \mathfrak{p}_+), \Phi(\mathfrak{u}_\lambda\cap \mathfrak{p}_-)$ only two are possible. 
Thus $N(11)\le 2$.  Both possibilities do occur as can be seen by choosing $\lambda= \varpi_1+\varpi_4+\varpi_6-4\varpi_7$ and 
$\varpi_3+\varpi_5-3\varpi_7$. 
\end{proof}
 
We tabulate in Table \ref{EVII} the weights $\lambda$ corresponding to $\theta$-stable parabolic subalgebras $\mathfrak{q}=\mathfrak{q}_\lambda$ 
with $1\le R_+(\mathfrak{q})=R_-(\mathfrak{q})\le 11.$  We also describe the compact Hermitian symmetric 
space $Y_{\mathfrak{q}}$ dual to $[L,L]/(K\cap [L,L])$ which is determined using Proposition \ref{hermitian}.

\begin{table}[ht]
\centering
\begin{tabular}{|c|c|c|c|}
\hline
$\lambda$ & $R_\pm(\mathfrak{q})$ & $Y_\mathfrak{q}$ & $\chi(Y_\mathfrak{q})$\\
\hline
$\varpi_2-\varpi_7$ & $6$ & $\SO(12)/U(6)$& $32$ \\
\hline 
$\varpi_2+\varpi_4-3\varpi_7$ & $9$& $G_3(\mathbb{C}^6)$ & $20$  \\
\hline
$\varpi_1+\varpi_6-2\varpi_7$ & $9$& $\mathbb{S}^2\times \SO(10)/U(5)$ & $32$\\
\hline
$\varpi_1+\varpi_2+\varpi_6-3\varpi_7$ & $10$ &  $\mathbb{S}^2\times \SO(8)/U(4)$  & $16$ \\
\hline
$\varpi_1+\varpi_4+\varpi_6-4\varpi_7$& $11$ &  $\mathbb{S}^2\times G_2(\mathbb{C}^4)$&$12$ \\
\hline
$\varpi_3+\varpi_5-3\varpi_7$ & $11$ & $\mathbb{S}^2\times G_2(\mathbb{C}^4)$ & $12$\\
\hline
\end{tabular}
~\\
\caption{The $\theta$-stable parabolic subalgebras $\mathfrak{q}=\mathfrak{q}_\lambda$ of type EVII that satisfy
$1\le R_\pm(\mathfrak{q})=r\le 11$, the symmetric 
spaces $Y_\mathfrak{q}$ and their Euler characteristics.}
\label{EVII}
\end{table}

\section{Proofs of Theorems \ref{main1} and \ref{main2}}\label{proofs}

Recall the set $\mathfrak{Q}$ of equivalence classes of $\theta$-stable parabolic subalgebras $\mathfrak{q}$ of $\mathfrak{g}_0$ where $\mathfrak{q}\sim\mathfrak{q}'$ if the unitary representations of $G$, 
$\mathcal{A}_\mathfrak{q}$ and $\mathcal{A}_{\mathfrak{q}'}$, are equivalent.   Let us denote the equivalence class of 
$\mathfrak{q}$ by $[\mathfrak{q}]$.  We shall denote by $\mathfrak{Q}^0$ the set consisting of $[\mathfrak{q}]$ such that 
$R_+(\mathfrak{q})=R_-(\mathfrak{q})$.  

The equivalence class determined by $\mathfrak{g}$ consists only of 
$\mathfrak{g}$ and the corresponding irreducible representation is trivial.  When $X=G/K$ is a Hermitian symmetric space 
and $\Gamma$ is a uniform lattice, the Matsushima isomorphism yields the following isomorphism:  
For any $p\ge 0$, 
\[H^{p,p}(X_\Gamma;\mathbb{C})\cong \bigoplus_{[\mathfrak{q}]\in \mathfrak{Q}^0}H^{p,p}(\mathfrak{g},K;A_{\mathfrak{q},K}).\]

Let $V\subset X_\Gamma$ be a closed analytic cycle, not necessarily a submanifold.  It is well-known that, in view of the fact that $X_\Gamma$ is a compact 
K\"ahler manifold (in fact even a projective variety by a result of Kodaira \cite{kodaira}), $V$ determines a fundamental homology 
class $\mu_V$ whose Poincar\'e dual $[V]$ is a non-zero class in $H^{p,p}(X_\Gamma;\mathbb{C})$ where $p=\textrm{codim}_XV$ is the complex codimension.  (See \cite{gh}.) 
Write $[V]=\sum_{[\mathfrak{q}]\in \mathfrak{Q}^0} [V]_{[\mathfrak{q}]}$ where $[V]_{[\mathfrak{q}]}\in H^{p,p}(\mathfrak{g},K;
A_{\mathfrak{q},K})$.   If $[V]$ is in the image of the Matsushima map $H^*(X_\textrm{u};\mathbb{C})\to H^*(X_\Gamma;\mathbb{C})$, then $[V]_{[\mathfrak{q}]}=0$ whenever $\mathfrak{q}\ne \mathfrak{g}$.

Recall the family of lattices $\mathcal{L}(G)$ defined in \S\ref{commutinginvolutions}. 
When $\Gamma\in \mathcal{L}(G)$, we can find a finite index subgroup $\Lambda\subset \Gamma$ and an 
involution $\sigma:G\to G$ such that (i) $\sigma(\Lambda)=\Lambda$, (ii) the special cycle $C(\sigma,\Lambda)\subset 
X_\Lambda$ is complex analytic, and (iii) the Poincar\'e dual $[C(\sigma,\Lambda)]$ of $C(\sigma, \Lambda)$ 
is not in the image of the Matsushima homomorphism $H^*(X_\textrm{u};\mathbb{C})\to H^*(X_\Lambda;\mathbb{C})$.  
In fact we may (and do) choose $\sigma$ so that $\textrm{codim}_{X_\Lambda} C(\sigma,\Lambda)$ is equal to $c(X)$.  
See \S\ref{codimension}.  It follows that, taking $V$ in the above to be equal to $C(\sigma,\Lambda)$, we obtain $[C(\sigma,\Lambda)]_{[\mathfrak{q}]}\ne 0$ for some $[\mathfrak{q}]\in \mathfrak{Q}^0\setminus\{[\mathfrak{g}]\}$ 
by \cite[Theorem 2.1]{mr}.
In view of our hypothesis on $G$, 
there exists a {\it unique} class $[\mathfrak{q}_0]\in \mathfrak{Q}^0$ such that $1\le R_\pm(\mathfrak{q})\le c(X)$, 
namely the one with $R_\pm(\mathfrak{q}_0)=r(\mathfrak{g}_0)$.  It follows that 
$[C(\sigma,\Lambda)]_{[\mathfrak{q}_0]}\in H^*(\mathfrak{g},K;A_{\mathfrak{q}_0,K})$ is non-zero.  
Therefore $\mathcal{A}_{\mathfrak{q}_0}$ occurs in $L^2(\Lambda\backslash G)$ with non-zero multiplicity $m(\mathfrak{q}_0,\Lambda)$.  

Now we take $V$ to be the image of $C(\sigma,\Lambda)$ under the (finite) covering projection 
$\pi: X_\Lambda\to X_\Gamma$.  
Note that this projection is holomorphic.  Then $V$ is a complex analytic submanifold and hence its Poincar\'e dual 
$[V]$ in $H^*(X_\Gamma;\mathbb{C})$ is non-zero. If $[V]$ could be represented by a $G$-invariant form, 
then so would $[C(\sigma,\Lambda)]$.  It follows that $[V]$ is not in the image of the Matsushima homomorphism 
and we are led to the conclusion that $[V]_{[\mathfrak{q}_0]}\ne 0$.   Hence, as before, $m(\mathfrak{q}_0;\Gamma)\ne 0$.  
This completes the proof of Theorem \ref{main1}.

Note that, when there are possibly more than one element $[\mathfrak{q}]\in \mathfrak{Q}^0$ with $1\le R_\pm(\mathfrak{q})
\le c(X)$,  the above argument is still applicable, but leads to the weaker conclusion that $[V]_{[\mathfrak{q}]}\ne 0$ for 
at least one such $[\mathfrak{q}]$.  

Finally to complete the proof of Theorem \ref{main2},  we observe that, for any $\Gamma\in \mathcal{L}(G)$,  there exists a 
$[\mathfrak{q}]\in \mathfrak{Q}^0$ with $r_0(\mathfrak{g})\le R_\pm(\mathfrak{q})\le c(X)$ with $m:=m(\mathfrak{q},\Gamma)>0$.
The corresponding representation $\mathcal{A}_{\mathfrak{q}}$ contributes $H^*(\mathfrak{g},K;A_{\mathfrak{q},K})^{\oplus m}$ 
to the cohomology of $X_\Gamma$.   Recall that, writing $r=R_+(\mathfrak{q})$, we have $H^{p,p}(\mathfrak{g},K;A_{\mathfrak{q},K})\cong H^{p-r,p-r}(Y_\mathfrak{q};\mathbb{C})$ where 
$Y_\mathfrak{q}$ is Hermitian symmetric of (complex) dimension $\dim_\mathbb{C} X-2r$ 
by Proposition \ref{hermitian}.   Since 
$H^{p,p}(\mathfrak{g},K;A_{\mathfrak{q},K})\cong H^{p-r, p-r}(Y_\mathfrak{q};\mathbb{C})$ is 
non-zero for $r\le p\le \dim_\mathbb{C} X-r$, 
and since $r\le c(X)$, this completes the proof.

The irreducible representations $\mathcal{A}_\mathfrak{q}$ whose occurrence with non-zero multiplicity in 
$L^2(\Gamma\backslash G), \Gamma\in \mathcal{L}(G)$ is asserted by Theorem \ref{main1} are listed in Table \ref{lambda} below, in terms of the $\mathfrak{k}$-dominant weight $\lambda$ such that $\mathfrak{q}=\mathfrak{q}_\lambda$.
This is based on the classification results obtained in \S\ref{pptypeparabolics}.

\begin{table}[h]
\centering
\begin{tabular}{|c|c|c|c|}
\hline
Type & $G$ & $\lambda$ & $R_\pm(\mathfrak{q}_\lambda)=r(\mathfrak{g}_0)$\\
\hline
AIII & $\SU(p,q)$& $\epsilon_{p+1}-\epsilon_{p+q}$&$p$\\
&$p<p-1, q\ge 5 $&&\\
\hline
BDI & $\SO_0(2,p),p\ge 3$ & $\epsilon_2$& $1$\\
\hline
CI & $\Sp(n,\mathbb{R}),n\neq 4$ & $\epsilon_1-\epsilon_n$& $n-1$\\
\hline
DIII & $\SO^*(2n), n\ge 9$ & $\epsilon_1-\epsilon_n$& $n-2$\\
\hline
\end{tabular}
\caption{List of $\mathcal{A}_{\mathfrak{q}_\lambda}$ with $R_\pm(\mathfrak{q}_\lambda)=r(\mathfrak{g}_0)$.}
\label{lambda}
\end{table}

\begin{remark}{\em
When $G=\SU(p,q), q=p+1$, there are three irreducible representations $\mathcal{A}_\lambda$ 
with $R_\pm(\mathfrak{q})=p=c(X)$. In view of this, arguing as above, from the non-vanishing of $[V]\in H^{p,p}(X_\Gamma;\mathbb{C})$ 
we can only infer that (at least) one of the three components $[V]_{[\mathfrak{q}]}$  
is non-zero.    But we are unable to decide whether a specific component is non-zero.   For this reason 
we obtain only a weaker conclusion that for any $\Gamma\in \mathcal{L}(G)$, for one of the 
representations $\mathcal{A}_{\mathfrak{q}_{\lambda}}, \lambda\in \{\epsilon_{p+1}-\epsilon_{p+q},p\epsilon_1+(p+1)\epsilon_{p+1}-\epsilon_0 ,\epsilon_0-p\epsilon_p-(p+1)\epsilon_{p+1}   \}$, (where $\epsilon_0=\sum_{1\le i\le p+q}\epsilon_i$),  the multiplicity $m(\mathfrak{q}_\lambda,\Gamma) \ne 0$.  

The same remark applies to $G=\SU(p,p)$ and we obtain that 
$m(\mathfrak{q}_\lambda,\Gamma)\ne 0$ for at least one of $\lambda\in \{\epsilon_{p+1}-\epsilon_{2p},\epsilon_1-\epsilon_p,  p\epsilon_1+p\epsilon_{p+1}-\epsilon_0, \epsilon_0-p\epsilon_p-p\epsilon_{p+q}\}$.   

Analogously, when $G$ is an 
exceptional Lie groups with Lie algebra $\mathfrak{e}_{6, (-14)}$ or $\mathfrak{e}_{7, (-25)}$ we see that 
if $\Gamma\in \mathcal{L}(G)$, then 
one of the representations $\mathcal{A}_\mathfrak{q}$ occurs with non-zero multiplicity $m(\mathfrak{q},\Gamma)$ 
in $L^{2}(\Gamma\backslash G)$ where $\mathfrak{q}$ is as in Table \ref{EIII} and Table \ref{EVII}.}
\end{remark}

Fix $\Gamma\in \mathcal{L}(G)$.  
Let $\mathfrak{q}$ be a $\theta$-stable parabolic subalgebra with $ R_+(\mathfrak{q})=R_-(\mathfrak{q})$ and let  
$Y_\mathfrak{q}$ be the Hermitian symmetric space associated to $\mathfrak{q}$ (see Proposition \ref{hermitian}).
In view of Theorem \ref{main1}, the above remark, and the Matsushima isomorphism, we obtain a monomorphism 
$H^{s,s}(Y_{\mathfrak{q}};\mathbb{C})\to H^{r,r}(X_\Gamma;\mathbb{C})$ (with $s=r-R_+(\mathfrak{q})$) for 
some $[\mathfrak{q}]\in \mathfrak{Q}^0$ such that $r(\mathfrak{g}_0)\le R_\pm(\mathfrak{q})\le c(X)$. 
This yields the following result which is stronger than Theorem \ref{main2}.   

\begin{theorem}  Let $\Gamma\in \mathcal{L}(G)$.  Then,  
for some $[\mathfrak{q}]\in \mathfrak{Q}^0$ with $r(\mathfrak{g}_0)\le R_\pm(\mathfrak{q})=:r_0\le c(X)$  and 
for every integer $r$ such that $R_+(\mathfrak{q})\le r\le \dim X_\Gamma-R_+(\mathfrak{q})$, 
there exist a vector subspace of dimension $b_{2(r-r_0)}(Y_\mathfrak{q})$, the $2(r-r_0)^{th}$ Betti number of $Y_\mathfrak{q}$, contained 
in $H^{r,r}(X_\Gamma;\mathbb{C})$ 
whose non-zero elements are not in the image of the Matsushima homomorphism $H^*(X_\textrm{u};\mathbb{C}) 
\to H^*(X_\Gamma;\mathbb{C})$. \hfill $\Box$ 
\end{theorem}

We tabulate in Table \ref{euler} the spaces $Y_\mathfrak{q}$ and their Euler characteristics in the cases when the hypotheses of 
Theorem \ref{main1} hold. (Thus $[\mathfrak{q}]\in \mathfrak{Q}_0$ is the unique one with $r(\mathfrak{g}_0)
=R_\pm(\mathfrak{q})\le c(X)$ given in Table \ref{lambda}.)    Recall that $Y_\mathfrak{q}$ have been determined in each case 
in Remarks \ref{supqyq}, \ref{bd1yq}, \ref{ciyq}, and \ref{d3yq}. 
We have used the following notations: 
$G_r(\mathbb{C}^n)$ is the Grassmannian $\textrm{U}(n)/(\textrm{U}(r)\times \textrm{U}(n-r))$
and $Q_p$ is the quadric $\SO(p+2)/\SO(2)\times \SO(p)$ which is also the real oriented 
Grassmann manifold $\tilde{G}_2(\mathbb{R}^n)$.  The Euler characteristics of compact globally 
Hermitian symmetric 
spaces are well-known and are given by the formula $\chi(M/H)=\#W(M,T)/\#W(H,T)$ where $H$ is a 
closed connected subgroup of a compact connected Lie group $M$ and 
$T\subset H$ is a maximal torus of $M$.   

\begin{table}[H]
\centering
\begin{tabular}
{|c|c|c|c|}
\hline
Type &$G$& $Y_{\mathfrak{q}}$ & $\chi(Y_\mathfrak{q})$\\
\hline
AIII & $\SU(p,q)$ & $G_p(\mathbb{C}^{p+q-2})$ &${p+q-2\choose p}$\\
& $p<q-1, q\ge 5$&  &\\
\hline
BDI & $\SO_0(2,p), p\ge 3$ & $Q_{p-2}$ & $2\lfloor p/2\rfloor$\\
\hline
CI & $\Sp(n,\mathbb{R}), n\ne 4$ & $\Sp(n-2)/\textrm{U}(n-2) \times \textrm{P}^1(\mathbb{C})$ & $2^{n-1}$\\
\hline 
DIII & $\SO^*(2n), n\ge 9$ & $\SO(2n-4)/\textrm{U}(n-1) \times \textrm{P}^1(\mathbb{C})$ & $2^{n-2}$\\
\hline
\end{tabular}
\caption{The spaces $Y_\mathfrak{q}$, where $R_\pm(\mathfrak{q})=c(X)$, and their Euler characteristics.}
\label{euler}
\end{table}

\noindent
{\bf Acknowledgements.}  
We thank T. N. Venkataramana for pointing out to us the paper of Li \cite{li}. We thank Sravan Danda for writing a programme in Python for finding out the least dimensional or least codimensional geometric cycle in EIII and EVII cases. Most of this work was carried out when Arghya Mondal was at the Institute of Mathematical Sciences, Chennai. 
Research of both the authors were partially supported by the Department of Atomic Energy, Government of India,
under a XII Plan Project.

\end{document}